\documentclass[oneside,sn-mathphys,Numbered]{sn-jnl}

\usepackage{graphicx}%
\usepackage{multirow}%
\usepackage{amssymb,amsfonts}%
\usepackage{amsthm}%
\usepackage{mathrsfs}%
\usepackage[title]{appendix}%
\usepackage{xcolor}%
\usepackage{textcomp}%
\usepackage{manyfoot}%
\usepackage{booktabs}%
\usepackage{algorithm}%
\usepackage{algorithmicx}%
\usepackage{algpseudocode}%
\usepackage{stfloats}
\usepackage{mathtools}
\usepackage{enumitem}
\usepackage{dsfont}
\usepackage{calc}
\usepackage{listings}%
\usepackage{comment}
\newif\ifdoubleblind
\doubleblindfalse 

\ifdoubleblind
  \makeatletter
  \def\blindmask#1{{\nullfont#1}}
  
   \def\fnm#1{\blindmask{#1}}
   \def\sur#1{\blindmask{#1}}
   \def\orgdiv#1{\blindmask{#1}}
   \def\orgname#1{\blindmask{#1}}
  
  \renewcommand{\email}[1]{\def\authemail{}}
  
  \def\au@and{}
  \makeatother

  \newcommand{\blindblock}[1]{%
    \phantom{\parbox{\linewidth}{\noindent#1}}%
  }
\else
  \newcommand{\blindblock}[1]{#1}
\fi

\allowdisplaybreaks

\newtheorem{theorem}{Theorem}
\newtheorem{proposition}[theorem]{Proposition}%
\newtheorem{lemma}[theorem]{Lemma}%

\newtheoremstyle{boldremark}%
  {6pt}
  {6pt}
  {\normalfont}
  {}
  {\bfseries}
  {.}
  {1em}
  {}

\theoremstyle{boldremark}
\newtheorem{remark}{Remark}
\newtheorem{assumption}{Assumption}%

\newtheorem{definition}{Definition}%

\newcommand{\dk}[1]{\textbf{\textcolor{red}{ {DK: #1} }}}

\newcommand\blfootnote[1]{%
  \begingroup
  \renewcommand\thefootnote{}\footnote{#1}%
  \addtocounter{footnote}{-1}%
  \endgroup
}

\allowdisplaybreaks

\makeatletter
\newenvironment{tagequation}[1]{%
  \refstepcounter{equation}
  \def\@currentlabel{\theequation}
  \def\theequation{#1}
  \begin{equation}
}{%
  \end{equation}
  \addtocounter{equation}{-1}
}
\makeatother

\colorlet{lightgray}{blue!5.5}
\colorlet{orange}{orange!10}

\begin{document}

\title[Strong Duality in Risk-Constrained Nonconvex Functional Programming]
{Strong Duality in Risk-Constrained

Nonconvex Functional Programming\blfootnote{\blindblock{This work has been supported by a Microsoft gift, and by the US NSF under Grant CCF 2242215. This work was developed  while S. Pougkakiotis was a postdoctoral researcher with the Department of ECE, Yale University.}
\par \noindent \blindblock{$^\dagger$Corresponding author.}}}


\author{\fnm{Dionysis} \sur{Kalogerias}$^\dagger$}
\email{dionysis.kalogerias@yale.edu}

\author{\fnm{Spyridon} \sur{Pougkakiotis}$^\star$}\email{spyridon.pougkakiotis@kcl.ac.uk}

\affil{\orgdiv{$^\dagger$Department of ECE}---\orgname{Yale 
 University, USA}}
 \affil{\orgdiv{$^\star$Department of Mathematics}---\orgname{King's College London, UK}}


\abstract{We show that a wide class of risk-constrained nonconvex functional optimization problems exhibit strong duality, regardless of nonconvexity. We develop two novel results under distinct sets of assumptions, establishing strong duality over both decomposable policies (matching and extending prior work in the risk neutral case) and nondecomposable policies with structure (e.g., continuity or smoothness), including certain universal finite-dimensional (fixed depth/width) neural network parametrizations as special cases (improving established results in the risk-neutral setting as well). Our results hold for constrained models featuring arbitrary convex risk functionals on $\mathcal{L}_p,p\in[1,\infty]$. 
We further discuss reductions and generalizations of our base model, and establish its necessity/minimality by presenting explicitly constructed counterexamples, rigorously refuting the manifestation of strong duality for certain, potentially simpler optimization models. Lastly, we discuss applications in wireless systems resource allocation, and supervised constrained learning. Our core proof technique appears to be new and relies on a non-trivial application of the Kingman-Robertson generalization of Lyapunov’s convexity theorem for vector measures taking values in infinite-dimensional topological spaces.

\vspace{-5pt}
}




\maketitle

\vspace{-22pt}
\section{Introduction and Problem Setting} \label{sec:Introduction}
\vspace{-2pt}

On an arbitrary base probability space $(\Omega,\mathscr{F},\mu)$, consider an element $\boldsymbol{S}\triangleq(\boldsymbol{H},\boldsymbol{U}):\Omega\rightarrow \mathcal{S}\triangleq {\cal H} \hspace{-1bp}\times {\cal U}$ with $\mathcal{H}\triangleq\mathbb{R}^{N_{\boldsymbol{H}}}$ and $\mathcal{U}\triangleq\mathbb{R}^{N_{\boldsymbol{U}}}$,
modeling an observable random phenomenon, which we would like to optimally
handle in a certain sense by making appropriate decisions. The induced Borel pushforwards of $\boldsymbol{S}$ and $\boldsymbol{H}$ are denoted as $\mathrm{P}_{\boldsymbol{S}}:\mathscr{B}({\cal S})\rightarrow[0,1]$ and $\mathrm{P}:\mathscr{B}({\cal H})\rightarrow[0,1]$, respectively.
We are interested in
risk-constrained nonconvex functional programs canonically formulated
as

\vspace{-3pt}
\begin{tagequation}{RCP}
\boxed{\begin{array}{rl}
\infty>\mathsf{P}^{*}=\:\underset{\boldsymbol{x},\boldsymbol{p}(\cdot)}{\mathrm{maximize}} & g^{o}(\boldsymbol{x})\\
\mathrm{subject\,to} & \boldsymbol{x}\le-\boldsymbol{\rho}(-
[\boldsymbol{F}\boldsymbol{p}]_{\boldsymbol{f}}(\boldsymbol{H}))\\
 & \boldsymbol{g}(\boldsymbol{x})\ge{\bf 0}\\
 & (\boldsymbol{x},\boldsymbol{p})\in{\cal X}\times\Pi
\end{array},}\label{eq:Base}
\end{tagequation}
\vspace{4pt}

\noindent where $g^{o}:\mathbb{R}^{N}\rightarrow\mathbb{R}$
and $\boldsymbol{g}:\mathbb{R}^{N}\rightarrow\mathbb{R}^{N_{\boldsymbol{g}}}$
are concave utility functions; $\mathcal{X}\subseteq \mathbb{R}^{N}$ is a convex set; $\boldsymbol{p}:{\cal S}\rightarrow\mathbb{R}^{N_{\boldsymbol{p}}}$
is a measurable allocation policy on observables $(\boldsymbol{H},\boldsymbol{U})$ and restricted to an  admissible set $\Pi$;
the \textit{``smooth operator"}
$[\boldsymbol{F}\hspace{1pt}\cdot\hspace{1pt}]_{\boldsymbol{f}}:\Pi \rightarrow {\cal L}_{p}(\mathrm{P},\mathbb{R}^{N}), p\in[1,\infty]$ is defined as a conditional expectation 
$[\boldsymbol{F}\boldsymbol{p}]_{\boldsymbol{f}}(\cdot)\triangleq \mathbb{E}\{
\boldsymbol{f}(\boldsymbol{p}(\cdot,U),\cdot,U)|\cdot\}$, where
$\boldsymbol{f}:\mathbb{R}^{N_{\boldsymbol{p}}}\times{\cal S}\rightarrow \mathbb{R}^{N}$
is an arbitrarily nonconvex \textit{instantaneous} performance 
score/reward, measuring the quality of a policy $\boldsymbol{p}$ at each realization $(\boldsymbol{H},\boldsymbol{U})$ and such that $\boldsymbol{f}(\boldsymbol{p}(\cdot),\cdot)\in{\cal L}_{p}(\mathrm{P}_{\boldsymbol{S}},\mathbb{R}^{N})$
for all $\boldsymbol{p} \in \Pi$; and lastly, $\boldsymbol{\rho}:{\cal L}_{p}(\mathrm{P},\mathbb{R}^{N})\rightarrow \mathbb{R}^{N}$ 
is a finite-valued \textit{vector} risk measure following the standardized convention (imposed for simplicity) that, for
each $i\in\mathbb{N}_{N}^{+}$,
\vspace{-4.5pt}
\[
\rho_{i}(\boldsymbol{Z})=\rho_{i}(Z_{i}),\quad\text{for all }\boldsymbol{Z}\in{\cal L}_{p}(\mathrm{P},\mathbb{R}^{N}),
\]
\vspace{-15.5pt}

\noindent and where $\rho_i$ is \textit{any} convex\footnote{To avoid confusion throughout, a risk measure is called convex \textit{if and only if} it is a convex functional of its argument; note that this is in contrast to, e.g., \cite[Definition 6.4]{ShapiroLectures_2ND}, where a risk measure is called convex if it satisfies additional conditions to mere convexity, namely, monotonicity and translation equivariance.} (and real-valued) functional on ${\cal L}_{p}(\mathrm{P},\mathbb{R})$.

Our assumptions support a very broad class of functionals including all classical risk measures such as the Conditional Value-at-Risk (CVaR), mean-variance, mean-standard deviation, mean-(upper-semi)deviations (including non-monotone cases) \citep{Ogryczak1999,Ogryczak2002}, distributionally robust functionals such as $f$-divergences (taken on $\mathcal{L}_\infty(\mathrm{P},\mathbb{R})$), as well as standard and large well-known subclasses with common features such as coherent risk measures \cite{Artzner1999}, spectral risk measures, distortion risk measures \cite{Balbas2009}, and beyond (see, e.g., \citep[Chapters 6 and 7]{ShapiroLectures_2ND}). A notable risk measure which is a nonconvex functional of its argument is the VaR, excluded from our treatment.

Back to \eqref{eq:Base}, we call $\boldsymbol{x}$ a \textit{performance risk vector}, since it is used to optimally bound the performance of a policy $\boldsymbol{p}$ evaluated using the composite functional $-\boldsymbol{\rho}(-
[\boldsymbol{F}\boldsymbol{p}]_{\boldsymbol{f}}(\boldsymbol{H}))$ within the corresponding constraint. Observe that the functional inequality constraints of \eqref{eq:Base} couple the finite/infinite-dimensional variables of the problem. As indicated in \eqref{eq:Base} and also mentioned above, performance risks are further restricted
to a finite-dimensional set ${\cal X}$ and policies are further restricted to an infinite-dimensional set
$\Pi$. More specific assumptions on the structure
of \eqref{eq:Base} enabling the development of the results presented
in this paper will be discussed in due course.

Problem \eqref{eq:Base} manifests itself in a variety of interesting applications.
As a canonical example coming from wireless systems engineering \cite{Ribeiro2012}, if $i\in\mathbb{N}_{N}^{+}$ refers to a user of a
wireless system, then the $i$-th entry of $\boldsymbol{\rho}$ 
may evaluate the risk associated with the service experienced by the $i$-th
user; hereafter, $\boldsymbol{f}$ will also be called the \textit{service
function}. However, with the exception of very special scenarios,
the service experienced by the $i$-th user \textit{can} (and, in
general, will) be dependent on the services experienced by the rest
of the users in the system. This is due to each service $f_{i}$ being
a function of the coupling variables $\boldsymbol{p}(\boldsymbol{S})$ ---the resources---
and $\boldsymbol{S}=(\boldsymbol{H},\boldsymbol{U})$ ---the uncertainty, hereafter referred to as a \textit{split state}---, which are common to all users in the system.
This general structure adheres to, among others, a large variety of resource
allocation problems encountered in practice (e.g., see Section \ref{subsec:Contributions} below, as well as Section \ref{subsec: resource allocation examples}). 

Temporarily and to motivate our setting, let us make the relatively mild additional assumption that each (real-valued) component of $\boldsymbol{\rho}$ is lower semicontinuous with respect to the corresponding weak topology on $\mathcal{L}_p(\mathrm{P},\mathbb{R})$ (e.g., $\mathrm{weak}$ topology for $p\in[1,\infty)$, $\mathrm{weak}^{*}$ topology for $p=\infty$, depending on the corresponding set of ``test" functions, relative to the canonical Hölder pairing).
Then, every $\rho_i$ is proper and admits the dual representation (by Fenchel-Moreau) 
\begin{equation}\nonumber
    {\mathbb{R}\ni\rho_i (Z) = \sup_{\zeta \in \mathrm{dom}(\rho_i^*)}
    \langle \zeta, Z \rangle - \rho_i^*(\zeta),\quad Z \in \mathcal{L}_p(\mathrm{P},\mathbb{R}),}
\end{equation}
where the functional $\rho_i^* : \mathcal{L}_q(\mathrm{P},\mathbb{R})\rightarrow (-\infty, \infty]$ is the convex conjugate of $\rho_i$ (Fenchel-Legendre transform), and $\langle \zeta, Z \rangle \triangleq \int \zeta Z  \mathrm{d}\mathrm{P}$ is the corresponding pairing. Under this provisioning, it is easily seen that our initial functional program may be re-expressed in the form
\vspace{-3.5pt}
\begin{tagequation}{RCP-E}
\boxed{\begin{array}{rl}
\underset{\boldsymbol{x},\boldsymbol{p}(\cdot)}{\mathrm{maximize}} & g^{o}(\boldsymbol{x})\\
\mathrm{subject\,to} & \boldsymbol{x}\le\inf_{\boldsymbol{\zeta}\in\mathds{A}^{S}}\mathbb{E}\{\boldsymbol{\zeta}(\boldsymbol{H})\odot[\boldsymbol{F}\boldsymbol{p}]_{\boldsymbol{f}}(\boldsymbol{H})\} + \boldsymbol{\rho}^*(\boldsymbol{\zeta}(\boldsymbol{H}))\\
 & \boldsymbol{g}(\boldsymbol{x})\ge{\bf 0}\\
 & (\boldsymbol{x},\boldsymbol{p})\in{\cal X}\times\Pi
\end{array},}\label{eq:Base_Problem_D}
\end{tagequation}

\noindent where ``$\odot$'' denotes the Hadamard product, the infimum is
understood in a component-wise manner, the service vector uncertainty
set (i.e., the vector risk envelope) $\mathds{A}^{S}$ is defined as
the Cartesian product $\mathds{A}^{S}\triangleq \mathrm{dom}(\rho^*_1)
\times\mathrm{dom}(\rho^*_2)
\times\ldots\times
\mathrm{dom}(\rho^*_N)$, and where ${\rho}_i^*(\boldsymbol{\zeta})={\rho}_i^*({\zeta_i}),i \in \mathbb{N}^+_N$, by convention.
\par In light of problem \eqref{eq:Base_Problem_D}, one may observe that the risk constraint $\boldsymbol{x} \leq -\boldsymbol{\rho}(-[\boldsymbol{F}\boldsymbol{p}]_{\boldsymbol{f}}(\boldsymbol{H}))$ (appearing in \eqref{eq:Base}) is equivalent to enforcing an uncountable multitude of constraints $\boldsymbol{x} \leq \mathbb{E}\{\boldsymbol{\zeta}(\boldsymbol{H}) \odot [\boldsymbol{F}\boldsymbol{p}]_{\boldsymbol{f}}(\boldsymbol{H})\} + \boldsymbol{\rho}^*(\boldsymbol{\zeta}(\boldsymbol{H}))$, for all $\boldsymbol{\zeta} \in \mathbb{A}^S$. Thus, in passing, it may be worth noting that (\ref{eq:Base_Problem_D})
is equivalent to the semi-infinite functional program 
\begin{tagequation}{RCP-I}
\begin{array}{rl}
\underset{\boldsymbol{x},\boldsymbol{p}(\cdot)}{\mathrm{maximize}} & g^{o}(\boldsymbol{x})\\
\mathrm{subject\,to} & \boldsymbol{x}\le
\mathbb{E}\{\boldsymbol{\zeta}(\boldsymbol{H})\odot[\boldsymbol{F}\boldsymbol{p}]_{\boldsymbol{f}}(\boldsymbol{H})\} + \boldsymbol{\rho}^*(\boldsymbol{\zeta}(\boldsymbol{H})),\quad\forall\boldsymbol{\zeta}\in\mathds{A}^{S}\\
 & \boldsymbol{g}(\boldsymbol{x})\ge{\bf 0}\\
 & (\boldsymbol{x},\boldsymbol{p})\in{\cal X}\times\Pi
\end{array},\label{eq:Base_Problem_I}
\end{tagequation}

\noindent illustrating the generality of the risk-constrained problem \eqref{eq:Base},
as compared with its risk-neutral counterpart, which may be obtained by
replacing the risk operator $\boldsymbol{\rho}$ by a finite-dimensional vector
of linear functionals (expectations).

Before proceeding, it is also worth discussing the operational role of the smooth operator  $[\boldsymbol{F}\hspace{1pt}\cdot\hspace{1pt}]_{\boldsymbol{f}}$ in problem \eqref{eq:Base}, which might look somewhat unnatural at first sight. 
Specifically, the presence of the conditional expectation in the definition of $[\boldsymbol{F}\hspace{1pt}\cdot\hspace{1pt}]_{\boldsymbol{f}}$ might indicate that we are somehow optimizing in ``sequence" or ``stages". This is not quite accurate though, and we can see this via the dual representations of the risk measures involved. Indeed, for any generic random variable $Z\in{\cal L}_{p}(\mathrm{P}_{\boldsymbol{S}},\mathbb{R})$
and given any (weakly) lower semicontinuous risk measure $\rho$ on ${\cal L}_{p}(\mathrm{P},\mathbb{R})$,
Fubini implies that 
\begin{align*}
\rho(\mathbb{E}\{Z(\boldsymbol{H},\boldsymbol{U})|\boldsymbol{H}\}) & =\sup_{\zeta\in\mathrm{dom}(\rho^{*})}\mathbb{E}\{\zeta(\boldsymbol{H})\mathbb{E}\{Z(\boldsymbol{H},\boldsymbol{U})|\boldsymbol{H}\}\}-\rho^{*}(\zeta(\boldsymbol{H}))\\
 & =\sup_{\zeta\in\mathrm{dom}(\rho^{*})}\mathbb{E}\{\zeta(\boldsymbol{H})Z(\boldsymbol{H},\boldsymbol{U})\}-\rho^{*}(\zeta(\boldsymbol{H}))\\
 & =\sup_{\zeta\in\mathrm{dom}(\rho^{*})}\mathbb{E}\{\zeta(\boldsymbol{H},\boldsymbol{U})Z(\boldsymbol{H},\boldsymbol{U})\}-\rho^{*}(\zeta(\boldsymbol{H},\boldsymbol{U}))\triangleq\rho'(Z(\boldsymbol{H},\boldsymbol{U})),
\end{align*}
with the extension of $\rho^{*}$ on the product space ${\cal L}_{q}(\mathrm{P}_{\boldsymbol{S}},\mathbb{R})$
taken as
\[
\rho^{*}(\zeta(\boldsymbol{H},\boldsymbol{U}))=\rho^{*}(\zeta(\boldsymbol{H}))\cdot\mathds{1}_{{\cal L}_{q}(\mathrm{P},\mathbb{R})}(\zeta(\boldsymbol{H},\boldsymbol{U}))+\infty\cdot\mathds{1}_{{\cal L}_{q}(\mathrm{P}_{\boldsymbol{S}},\mathbb{R})\setminus{\cal L}_{q}(\mathrm{P},\mathbb{R})}(\zeta(\boldsymbol{H},\boldsymbol{U})),
\]
where we may view ${\cal L}_{q}(\mathrm{P},\mathbb{R})$ as a subspace
of ${\cal L}_{q}(\mathrm{P}_{\boldsymbol{S}},\mathbb{R})$, in the
sense that elements on the former may be viewed as elements on the
latter which are almost everywhere constant in the second part of
the split state $(\boldsymbol{H},\boldsymbol{U})$. Thus, we verify
equivalence with another convex and weakly lower semicontinuous risk measure $\rho'$ on
${\cal L}_{p}(\mathrm{P}_{\boldsymbol{S}},\mathbb{R})$, which implements
\textit{partial adversarial (i.e., worst-case) reweighting} of $Z(\boldsymbol{H},\boldsymbol{U})$
restricted to the first part of the split state $(\boldsymbol{H},\boldsymbol{U})$.
In the perhaps more familiar setting where $\rho$ is a coherent risk
measure to start with, it is easy to see that $\rho'$ admits a cleaner
\textit{partial distributionally robust} representation
\[
\rho'(Z(\boldsymbol{H},\boldsymbol{U}))=\sup_{\frac{\mathrm{dQ}}{\mathrm{dP}}\in\mathds{A}}\mathbb{E}_{(\boldsymbol{h},\boldsymbol{u})\sim\mathrm{Q}\times\mathrm{P}_{\boldsymbol{U}|\boldsymbol{H}}}\{Z(\boldsymbol{h},\boldsymbol{u})\}.
\]
Note that, unless $\boldsymbol{H}$ and $\boldsymbol{U}$ are statistically
independent, Bayes (disintegration) implies that each choice of $\mathrm{Q}$
also results in a different marginal for $\boldsymbol{U}$ as compared
with that induced by $\mathrm{P}_{\boldsymbol{S}}$, despite invariance
of the reference conditional distribution $\mathrm{P}_{\boldsymbol{U}|\boldsymbol{H}}$.


As it turns out, including operator  $[\boldsymbol{F}\hspace{1pt}\cdot\hspace{1pt}]_{\boldsymbol{f}}$ in the definition of \eqref{eq:Base} is in fact indispensable for the purposes in this work (i.e., for recovering strong duality relations), as we discuss briefly in the next subsection, and in full detail in later sections.

\subsection{Origins, Related Literature and Contributions}\label{subsec:Contributions}

We will be mainly invested in investigating \textit{Lagrangian duality relations} for the infinite-dimensional constrained program \eqref{eq:Base} ---in particular, whether strong duality holds or not--- under general assumptions, notably without enforcing any additional structure on the integrable service vector function $\boldsymbol{f}(\cdot,\boldsymbol{S})$ (or by enforcing merely continuity in some cases, as discussed later). Recall that strong duality is a very desirable property for both theoretical and practical purposes and ensures in particular existence of optimal dual variables (Lagrange multipliers), in addition to a null duality gap \cite{Boyd2004,Ruszczynski2006b,Bertsekas2009}; see Section \ref{sec: Lagrangian duality} for a gentle exposition focusing on our purposes herein.


Although we are mainly interested in the infinite-dimensional setting, it may be worth noting that strong (Lagrangian) duality in finite-dimensional nonconvex optimization is typically unlikely to hold. While there are special (and highly structured) cases in which strong duality is present in finite-dimensional optimization (e.g., see, among others, \cite{CHIEU2020441,geoffrion_1974,hoffman_kruskal_1956,Lemaire1998}), a general treatment is, to the best of our knowledge, not possible without the utilization of additional, often unusual or hard-to-verify conditions. As we argue herein, in the context of (stochastic) functional optimization (as in the case of \eqref{eq:Base}), general (and practical) conditions on the underlying probability measure and the policy space involved reveal surprising structure, enabling a comprehensive study of strong duality under a very generic and versatile technical framework.

Before we discuss related prior work and existing results, it would be helpful to first present evidence justifying the usefulness and relevance of \eqref{eq:Base} in applications by mainly capitalizing on the standard \textit{risk-neutral} case, i.e., when $\boldsymbol{\rho}$ is a vector of expectations. Such evidence will hopefully provide at least some preliminary motivation for studying \eqref{eq:Base} in the general risk-constrained setting we consider in this paper. Indeed, due to the modular structure of \eqref{eq:Base}, and beyond its applicability in numerous settings in resource allocation for networking and wireless communications \cite{Tassiulas1992,Neely2005,Georgiadis2006a,Sidiropoulos2006,WeiYu2006,Bazerque2007, Luo2008,Neely2010, Ribeiro2010, Ribeiro2010a,  Ribeiro2012, Hu2012, He2014, Naderializadeh2014, Eisen2019, Eisen2019b, Kalogerias2020e}, instances of \eqref{eq:Base} appear in areas including nonlinear sparse functional programming \cite{Chamon2020a} with applications such as nonlinear line spectrum estimation and robust functional data analysis \cite{Chamon2020a}, nonconvex constrained machine learning \cite{Chamon2021}, the so-called ``fluid problem"
in dynamic resource allocation \cite{Balseiro2023} (with various applications such as network dynamic pricing, network revenue management, dynamic bidding, online matching, and order fulfillment; see \cite{Balseiro2023} and the references therein), and wireless control \cite{Gatsis2015,Gatsis2018}, to name a few. In all those instances of the constrained problem \eqref{eq:Base}, nonconvexity is abundant, and it should be expected that a rigorous dual-domain characterization of \eqref{eq:Base}, although highly desirable, is a challenging task.  

Dual-domain analysis and strong duality results concerning the \emph{risk-neutral} counterpart of \eqref{eq:Base} just discussed, although part of a much older and longer story (see \cite{BerliocchiLasry1973} as well as the discussion next paragraph), have a relatively short history in their most contemporary form, starting, to the best of our knowledge, in $2008$ with the seminal paper \cite{Luo2008} by Luo and Zhang in the context of wireless spectrum management, which was in turn based on an earlier preliminary discovery reported by Yu and Lui in \cite{WeiYu2006}. Since then, and after subsequent developments by Ribeiro and Giannakis in \cite{Ribeiro2010,Ribeiro2012}, there has been a flurry of research activity positioned around \eqref{eq:Base} and its various applications, in particular informing an extensive literature on dual-domain methods for obtaining optimal resource policies in the context of networking and wireless systems, including recent advances in model-free learning for wireless communications ---see, e.g., references in the previous paragraph---, such as approximate strong duality relations for finite-dimensional versions of \eqref{eq:Base} exploiting convolutional smoothing and universal policy parameterizations \cite{Eisen2019, Kalogerias2020e}. Strong duality of \eqref{eq:Base} in the risk-neutral setting has also been studied in \cite{Chamon2020a}, and more recently leveraged in \cite{Chamon2021} to explore a notion of (dual) PAC learnability in the context of constrained machine learning.

More specifically and from a technical standpoint, and following a long tradition of classical results in functional optimization problems arising in calculus of variations and mathematical economics (see, e.g.,  \cite{AubinEkeland1976,AumannPerles1965}, \cite[Appendix I.4: The Lyapunov Effect]{EkelandTemam1999} and the review presented in \cite{BerliocchiLasry1973}), articles \cite{Luo2008} and subsequently \cite{Ribeiro2010,Ribeiro2012} established the remarkable fact that, under a set of standard and generic assumptions (see our Assumption \ref{assu:Assumption} in Section \ref{subsec: decomposable}) ---the most important of which being the 
\textit{nonatomicity} of the reference Borel measure $\mathrm{P}_{\boldsymbol{S}}$--- the risk-neutral counterpart of \eqref{eq:Base} exhibits strong duality (also implying a zero duality gap) regardless of the nature of the service $\boldsymbol{f}(\cdot,\boldsymbol{S})$; the latter may in fact be \textit{arbitrary}, in particular (component-wise) nonconcave or even discontinuous. In a sense (see, e.g., \citep[Theorem 1 and its proof]{Chamon2020a}), a nonatomic distribution $\mathrm{P}_{\boldsymbol{S}}$ ``patches the holes" of the \textit{image space set}  of \eqref{eq:Base} \cite{Giannessi2005}, thus preserving its convexity; this follows from a careful application of the celebrated theorem of A. A. Lyapunov \cite[Corollary IX.1.6]{Diestel1977} (see also \cite{Lyapunov_original} for the original paper written by Lyapunov and \cite{AUMANN19651} for a related result due to Aumann), on the convexity of the range of finite-dimensional nonatomic vector measures. The guaranteed convexity of the image space set of \eqref{eq:Base} can then be utilized to establish strong duality ---itself a consequence of convex geometry \cite{FloresBazan2013}--- of \eqref{eq:Base} in the risk-neutral case via a standard and almost elementary application of the supporting hyperplane theorem.
The latter argument is a special part of another longer story on strong duality relations in general conic programming in infinite dimensions, for instance as considered in \cite{FloresBazan2013}; see also the books \cite{Luenberger1968,Giannessi2005}.

Evidently, if its risk-neutral version is at all valuable, then so is the substantially more general risk-constrained problem \eqref{eq:Base}, which is the focus of this work. This is because most, if not all, problems discussed above admit at least one interpretable and practically justifiable risk-averse reformulation (by replacing expectations with appropriate and meaningful risk measures, at the very least). For a concrete example, we refer the reader to a series of recent articles \cite{Yaylali2023, Yaylali2024,Yaylali2025,Yaylali2026}, where \eqref{eq:Base} is considered in the ``canonical" context of wireless systems resource allocation,  albeit under a classical and stylized convex programming setting (thus ensuring strong Lagrangian duality under standard constraint qualifications). For example, in \cite{Yaylali2023}, the standard expectation is replaced by the Conditional Value-at-Risk (CVaR) \citep{Rockafellar1997,ShapiroLectures_2ND}, provably resulting in optimal resource policies ensuring robust and reliable system performance, and operationally desirable quality-of-service.  


Despite their relevance, though, risk-constrained policy search programs such as in the form of \eqref{eq:Base} are not currently as well-studied as their risk-neutral versions. We believe that this is natural and mainly due to the fact that, specifically in the nonconvex setting, dual-domain analysis and properties of \eqref{eq:Base}, which are crucial for further developments such as the design of efficient algorithms for tackling such problems, are currently absent and essentially unexplored. Our work in this paper is exactly on initiating an effort for illuminating and rigorously addressing those issues. 

\textbf{{Contributions:}} We primarily show that, perhaps surprisingly, \textit{the risk-constrained functional problem \eqref{eq:Base} exhibits strong duality under almost exactly the same assumptions as in the risk-neutral case} (i.e., Assumption \ref{assu:Assumption} in Section \ref{subsec: decomposable}); 
the only difference is a demand for nonatomicity of the conditional distribution $\mathrm{P}_{\boldsymbol{U}|\boldsymbol{H}}$, rather than nonatomicity of the full distribution $\mathrm{P}_{\boldsymbol{S}}$ (which is the case in the risk-neutral setting), with the two conditions being equivalent if additionally $\mathrm{P}_{\boldsymbol{S}}$ admits a Lebesgue density. No further assumptions are required,
and the result holds for numerous risk measures, namely all (finite-valued) convex functionals on the space $\mathcal{L}_p,p\in[1,\infty]$, where the operator $\boldsymbol{\rho}$ may be comprised of different such risk measures and/or combinations of those in each dimension (for a common $p$). 
Further, despite the slight (partial) in-agreement on the assumptions, our result reduces exactly to that in the risk-neutral case, unconditionally (Section \ref{Counter}).

\par A core assumption in both the aforementioned and all previous results in the related literature (for the risk-neutral case; see discussion above) is the imposition of the property of \textit{decomposability} on the policy search space $\Pi$ in \eqref{eq:Base} (Assumption \ref{assu:Assumption}, condition 5). Decomposability is a well-established notion (e.g., see \cite{Rockafellar1971Integrals, HIAI1977149}), and instrumental in establishing  a wide array of results, ranging from strong duality (e.g., see \cite{Ribeiro2012,PennanenPerkkio2024}), to the well-known interchangeability principle of integration (e.g., see \cite[Theorem 14.60]{Rockafellar2009VarAn}). Nevertheless, an issue with decomposability is that it mostly holds for policy spaces with arbitrarily complex structure (e.g., all $\mathcal{L}_p$ spaces, for $p\in[1,\infty]$, or Orlicz spaces), potentially hindering interpretability and tractability of the corresponding solutions of \eqref{eq:Base}, even in the presence of strong duality. Importantly, decomposability is not compatible with structured policies (e.g., continuous, smooth, or those lying in Reproducing Kernel Hilbert Spaces (RKHSs)), which are desirable in certain applications (e.g., (constrained) learning). 

To circumvent these shortcomings, we relax decomposability by introducing the notion of a \textit{densely decomposable set} (see Definition \ref{def: densely decomposable set} and related discussion), enabling the direct imposition of structure (e.g., continuity or smoothness) on the feasible policies. Within this setting, we establish another novel (even in the risk-neutral case) strong duality result for \eqref{eq:Base} under an alternative set of conditions (i.e., Assumption \ref{assu:Assumption2} in Section \ref{subsec: densely decomposable sets}), while still allowing for arbitrarily nonconvex behavior of the service function $\boldsymbol{f}$. A  remarkable byproduct of this result is verification of strong duality in certain \textit{finite-dimensional}, neural network-induced, parametrized instances of \eqref{eq:Base}.

To the best of our knowledge, the general technique we develop to establish strong duality of \eqref{eq:Base} is new. It is inspired from exploiting risk duality in \eqref{eq:Base} so that its risk constraints ---a finite-dimensional vector of nonlinear functionals--- can be \textit{lifted to a vector of affine functionals of uncountably infinite dimension}; see equivalent problem (\ref{eq:Base_Problem_D}) (under lower semicontinuity). The hope is for such a collection of affine functionals to be
easier to handle than the risk measures they represent, by leveraging elements of functional analysis, specifically vector measure theory \cite{Diestel1977}. 

However, in contrast to the risk-neutral setting, application of Lyapunov's convexity theorem would not be possible in the case of (\ref{eq:Base_Problem_D}), since the corresponding vector measure construction ---as in \cite{Luo2008,Ribeiro2010,Ribeiro2012}--- is infinite-dimensional. For such vector measures, the conclusions Lyapunov's convexity theorem provably may not hold in general; see, e.g., counter-examples in \cite[Section IX. 1]{Diestel1977}. To deal with this fundamental difficulty, we leverage another well-known and powerful result, namely, the \textit{Kingman-Robertson convexity theorem} \cite{KingmanRobertson} (later further extended by Knowles \cite{SIAMCon:Knowles}), which massively generalizes Lyapunov’s theorem by providing necessary \textit{and} sufficient conditions for range convexity and compactness of density-admitting vector measures (in a certain precise sense \cite{KingmanRobertson}) which take values in product topological spaces of even uncountably infinite dimension.

The conditions identified by the Kingman-Robertson theorem are related with their notion of a \textit{thin set of functions} (essentially an infinite-dimensional kernel-related condition), and a significant part of our technique is on verifying such thinness of the set of derivatives associated with our uncountably infinite-dimensional vector measure construction (Section \ref{sec:Proof-of-Theorems}), motivated by the fact that the functional constraints of problem (\ref{eq:Base_Problem_D}) are all affine. This enables us to leverage the Kingman-Robertson theorem, even \textit{bypassing} the need for utilization of any dual representation of the risk measures involved. In the decomposable case (Assumption \ref{assu:Assumption}, Section \ref{subsec: decomposable}) we prove that image space set associated with \eqref{eq:Base} (the latter denoted as $\mathcal{C}$) is convex, providing a clean (though nontrivial) extension to the corresponding technique employed for the risk-neutral case. This conclusion then suffices to establish strong duality of \eqref{eq:Base}, as a consequence of convex geometry.
Alternatively, in the nondecomposable case (Assumption \ref{assu:Assumption2}, Section \ref{subsec: densely decomposable sets}), we leverage our notion of a densely decomposable set to show that the (norm) \textit{closure} of the image space set $\mathcal{C}$ is convex. In this case, we establish strong duality via an extension to the classical hyperplane separation argument (used, e.g., in the decomposable case), by exploiting convexity of $\mathrm{cl}(\mathcal{C})$ only \cite{Bertsekas2009}.

A natural question we faced after the development of the results above is whether the presence of the smooth operator $[\boldsymbol{F}\hspace{1pt}\cdot\hspace{1pt}]_{\boldsymbol{f}}$ is necessary while posing \eqref{eq:Base} (see discussion before this subsection for its operational properties), and whether one could assert strong duality by removing operator $[\boldsymbol{F}\hspace{1pt}\cdot\hspace{1pt}]_{\boldsymbol{f}}$ in two intuitively reasonable cases: that of deterministic policies (i.e., as in \eqref{eq:Base}, but where risk-aversion is over the full split state $\boldsymbol{S}$), or that of ``randomized policies", in which $\boldsymbol{U}$ enters exclusively the policy $\boldsymbol{p}(\boldsymbol{H},\cdot)$ as ``pure noise" (and nowhere else in \eqref{eq:Base}).
We rigorously answer this question \textit{to the negative for both cases} (perhaps surprisingly at first), by providing an explicitly constructed mean-CVaR instance of \eqref{eq:Base} satisfying our conditions (say decomposable case, Assumption \ref{assu:Assumption}) and applicable to both cases simultaneously, which exhibits \textit{a non-zero duality gap} available in closed form as a function of problem-related quantities (Section \ref{Counter}). This counterexample rigorously proves the necessity of the presence of the integral operator $[\boldsymbol{F}\hspace{1pt}\cdot\hspace{1pt}]_{\boldsymbol{f}}$ in \eqref{eq:Base} for strong duality to manifest, within our assumption system. We also discuss some further extensions of our theory (Section \ref{Counter}), including basic variations of \eqref{eq:Base} in which risk measures might appear in the objective as well, or where linear inequalities might be substituted by appropriate convex conic inequalities. 



\par Lastly, we demonstrate the versatility of our base model and assumption system by looking at two application areas. We first discuss the natural fit of \eqref{eq:Base} in the context of risk-constrained resource allocation in wireless communication systems, and demonstrate the compatibility of \eqref{eq:Base} for two relevant models, namely, a
multiple access interference channel \cite{Kalogerias2020e,Eisen2019}, and a
frequency division broadcast channel \cite{Ribeiro2012}, previously studied in the risk-neutral setting.
Subsequently, we consider the setting of functional risk-constrained supervised learning, in which the associated loss functions are allowed to be nonconvex \cite{Chamon2021}. By exploiting our result for densely decomposable policy spaces and while remaining compatible with existing literature, we show that strong duality is guaranteed for a wide class of relevant nonconvex risk-constrained learning problems. Further, we recover and improve existing results on risk-neutral learning obtained in \cite{Chamon2021} under relaxed assumptions, and enable a unified treatment of risk-constrained regression and classification tasks. 

As a remark in passing, we would like to mention that while Lyapunov's convexity theorem has many profound applications in various areas with obvious practical interest, this does not (yet) seem to be the case for its infinite-dimensional extensions, notably by Kingman and Robertson \cite{KingmanRobertson}, Knowles \cite{SIAMCon:Knowles}, Uhl \cite{Uhl1969}, and Kadets and Schechtman \cite{Kadets1992}. 
Therefore, we believe that the recognition of the usefulness of the Kingman-Robertson theorem in this paper, as a core technical ingredient in proving strong duality of \eqref{eq:Base} in the ---hopefully practically relevant--- context of risk-constrained policy optimization  under a standard Borel space setup, is quite interesting, just some
fifty-eight ($58$) years after its formulation (the Kingman-Robertson theorem was published in $1968$ \cite{KingmanRobertson}).

\subsection{Structure and Notation}

\par Paper structure is outlined as follows. 
In Section \ref{sec: Lagrangian duality}, we briefly introduce standard Lagrangian duality in the context of program \eqref{eq:Base}. In Section \ref{sec: Main results}, we state and discuss our assumptions as well as the main results of this work ---establishing strong duality of \eqref{eq:Base}---, which we proceed to prove in Section \ref{sec:Proof-of-Theorems}. Reduction to the risk-neutral setting, our counterexamples and extensions to our main results are discussed in Section \ref{Counter}. Then, in Section \ref{sec:Applications}, we elaborate on two relevant and general application settings that can be seen as special cases of \eqref{eq:Base}, namely risk-constrained wireless systems resource allocation and risk-constrained learning, and finally conclude in Section \ref{sec: Conclusions}.

\textbf{Notation ---also applicable above---:} Calligraphic
letters (such as ${\cal A}$), or plain capital letters
(such as $A$) will denote finite-dimensional sets/spaces, such as
Euclidean spaces. Double stroke letters (such as $\mathds{A}$) will
denote infinite-dimensional sets/spaces, such as Banach spaces. Math
script letters (such as $\mathscr{A}$) will denote $\sigma$-algebras.
Boldsymbol letters (such as $\boldsymbol{A}$ or $\boldsymbol{a}$)
will denote (possibly random) vectors. The space of $p$-integrable functions
from a measurable space $(\Omega,\mathscr{F})$ equipped with a finite
measure $\mu:\mathscr{F}\rightarrow\mathbb{R}_{+}$ to a space
$\mathds{A}$, with standard notation ${\cal L}_{p}(\Omega,\mathscr{F},\mu;\mathds{A})$
\citep{Diestel1977,ShapiroLectures_2ND}, is abbreviated as ${\cal L}_{p}(\mu,\mathds{A})$,
as it is also common practice \citep{Diestel1977}. The rest of the
notation is standard.

\section{Lagrangian Duality} \label{sec: Lagrangian duality}

A celebrated approach for dealing with the explicit inequality
constraints of the risk-aware problem \eqref{eq:Base} is by exploiting
\textit{Lagrangian duality}, which has been proven fundamental in analyzing
and efficiently solving constrained convex optimization problems;
see, e.g., \citep{Boyd2004,Ruszczynski2006b,Bertsekas2009}. Note,
however, that since the services $\boldsymbol{f}(\cdot,\boldsymbol{S})$ appearing in problem \eqref{eq:Base} are nonconcave in general (i.e., with respect to the first argument corresponding to the policy $\boldsymbol{p}$), standard results in Lagrangian duality for convex optimization
do not automatically apply (e.g., see \cite[Chapters 4 and 6]{Bertsekas2016Nonlinear}). On top of that, one has to incorporate
the structural complexity of the vector risk measure $\boldsymbol{\rho}$,
which is a nonlinear operator. As a result, fundamental
properties of expectation (being the most trivial risk measure) which
enable an elegant and straightforward analysis, such as linearity,
do not hold for $\boldsymbol{\rho}$.

The \textit{Lagrangian function} $\mathsf{L}:\mathbb{R}^{N}\times\Pi\times\mathbb{R}^{N_{\boldsymbol{g}}}\times\mathbb{R}^{N}\rightarrow\mathbb{R}$
associated with the risk-constrained problem \eqref{eq:Base} is defined
by scalarizing its constraints as
\[
\mathsf{L}(\boldsymbol{x},\boldsymbol{p},\boldsymbol{\lambda})\triangleq g^{o}(\boldsymbol{x})+\langle\boldsymbol{\lambda}_{\boldsymbol{g}},\boldsymbol{g}(\boldsymbol{x})\rangle+\big\langle\boldsymbol{\lambda}_{\boldsymbol{\rho}},-\boldsymbol{\rho}(-[\boldsymbol{F}\boldsymbol{p}]_{\boldsymbol{f}}(\boldsymbol{H}))-\boldsymbol{x}\big\rangle,
\]
where $\boldsymbol{\lambda}\equiv(\boldsymbol{\lambda}_{\boldsymbol{g}},\boldsymbol{\lambda}_{\boldsymbol{\rho}})\in\mathbb{R}^{N_{\boldsymbol{g}}}\times \mathbb{R}^{N}$
are dual multipliers associated with the constraints of \eqref{eq:Base}.
Then the \textit{dual function} $\mathsf{D}:\mathbb{R}^{N_{\boldsymbol{g}}}\times \mathbb{R}^{N} \rightarrow[-\infty,\infty]$
is defined as 
$
\mathsf{D}(\boldsymbol{\lambda})\triangleq\sup_{(\boldsymbol{x},\boldsymbol{p})\in{\cal X}\times\Pi}\mathsf{L}(\boldsymbol{x},\boldsymbol{p},\boldsymbol{\lambda}).
$
If the optimal value of problem \eqref{eq:Base} is $\mathsf{P}^{*}\in\mathbb{R}$,
it is then understood that $\mathsf{P}^{*}\le\mathsf{D}$ on the
positive orthant (i.e., for $\boldsymbol{\lambda}\ge{\bf 0}$), and
thus it is most reasonable to consider the \textit{dual problem $\inf_{\boldsymbol{\lambda}\ge{\bf 0}}\mathsf{D}(\boldsymbol{\lambda})$},
which is always convex and whose optimal value 
$
\mathsf{D}^{*}\triangleq\inf_{\boldsymbol{\lambda}\ge{\bf 0}}\mathsf{D}(\boldsymbol{\lambda})\in(-\infty,\infty]
$
serves as the tightest over-estimate of the optimal value of problem
\eqref{eq:Base}, $\mathsf{P}^{*}$, when knowing only $\mathsf{D}$.
Then, a basic question is whether
we can essentially replace an original constrained problem with its
dual, in the sense that $\mathsf{P}^{*}\equiv\mathsf{D}^{*}$; in
such a case, we say that the problem has \textit{zero duality gap}.
Referring to \eqref{eq:Base}, this would imply that one could instead look at the minimax problem
\begin{equation}\nonumber
\inf_{\boldsymbol{\lambda}\ge{\bf 0}}\mathsf{D}(\boldsymbol{\lambda})\equiv\inf_{\boldsymbol{\lambda}\ge{\bf 0}}\sup_{(\boldsymbol{x},\boldsymbol{p})\in{\cal X}\times\Pi}\mathsf{L}(\boldsymbol{x},\boldsymbol{p},\boldsymbol{\lambda}),
\end{equation}
whose optimal value is $\mathsf{D}^{*}$. If, in addition, $\mathsf{D}^{*}$ is
attained, i.e., an optimal multiplier vector exists, then \eqref{eq:Base}
is said to exhibit \textit{strong duality}. A zero duality gap also
implies that problem \eqref{eq:Base} satisfies the \textit{saddle
point property }(whether a saddle point exists or not), which is expressed
as
\[
\mathsf{D}^{*}=\inf_{\boldsymbol{\lambda}\ge{\bf 0}}\sup_{(\boldsymbol{x},\boldsymbol{p})\in{\cal X}\times\Pi}\mathsf{L}(\boldsymbol{x},\boldsymbol{p},\boldsymbol{\lambda})=\sup_{(\boldsymbol{x},\boldsymbol{p})\in{\cal X}\times\Pi}\inf_{\boldsymbol{\lambda}\ge{\bf 0}}\mathsf{L}(\boldsymbol{x},\boldsymbol{p},\boldsymbol{\lambda})=\mathsf{P}^{*}.
\]

Zero duality gaps are desirable: Because there is a finite number
of constraints, the dual function is finite-dimensional even though
the original functional problem \eqref{eq:Base} is infinite-dimensional.
Additionally, for every choice of the dual variable $\boldsymbol{\lambda}\ge{\bf 0}$
(and therefore for any optimal multiplier vector), joint maximization
of the Lagrangian $\mathsf{L}(\boldsymbol{x},\boldsymbol{p},\boldsymbol{\lambda})$
over the pair $(\boldsymbol{x},\boldsymbol{p})\in{\cal X}\times\Pi$
is separable. We thus see that duality transforms a constrained problem
into an unconstrained problem in a principled way and,
provided that the original constrained stochastic program (in our
case \eqref{eq:Base}) exhibits zero duality gap or strong duality, presents a general
methodological approach to study it. While zero duality gaps are
a common and fundamental characteristic of problems in convex 
constrained optimization (under appropriate regularity conditions),
proving zero duality gaps in nonconvex problems such as \eqref{eq:Base}
is a much more delicate and challenging task.

\section{Main Results} \label{sec: Main results}

In this section, we present the main contributions of this work. Specifically, in Section \ref{subsec: decomposable},  we provide minimal conditions ensuring that \eqref{eq:Base} exhibits strong duality notably without any assumption on the structure on the service function $\boldsymbol{f}$ (besides integrability of $\boldsymbol{f}(\boldsymbol{p}(\boldsymbol{S}),\boldsymbol{S})$), relying in particular on the decomposability of the policy space $\Pi$ and the nonatomicity of the conditional distribution $\mathrm{P}_{\boldsymbol{U}|\boldsymbol{H}}$ (cf. Assumption \ref{assu:Assumption} below).

A potential issue with the imposition of decomposability is that the corresponding search spaces are typically large, in the sense that they contain policies with arbitrarily complex structure (e.g., well-beyond continuity or smoothness). Consequently, such spaces may not be appropriate for constrained problems demanding structured decisions (i.e., continuous or smooth), such as those often arising in modern machine learning applications (e.g., RKHSs). To this end (Section \ref{subsec: densely decomposable sets}), we relax decomposability by introducing the notion of a \emph{densely decomposable set}, enabling the native imposition of structure (e.g., continuity or smoothness) on the feasible policies while still allowing for the establishment of strong duality of \eqref{eq:Base}, however under a mild pointwise continuity assumption on the (still arbitrarily nonconvex) service function $\boldsymbol{f}$ (cf. Assumption \ref{assu:Assumption2}). Remarkably, we show that under this setting one can ascertain strong duality for \eqref{eq:Base} even when optimizing over certain \textit{finite-dimensional}, neural network-induced, parametrized policy spaces with a universal behavior. 
\setcounter{subsection}{-1}
\subsection{WarmUp: Convexification via Randomization (``Convex-Analytic" Forms)}\label{warmup}
Before proceeding with the rigorous formulation of
our results, we would like to demonstrate the advantages of the presence
of the smooth operator in \eqref{eq:Base}, restricted in a special
``easy'' setting where the observable $\boldsymbol{U}$ of the split
state $\boldsymbol{S}=(\boldsymbol{H},\boldsymbol{U})$ acts as \textit{pure
policy randomization}. For simplicity in this subsection, we set $N^{\boldsymbol{U}}=1$,
and we take $\boldsymbol{U}\equiv U$ as a standard uniform random
variable essentially supported on $[0,1]$ (equivalently, Lebesgue
measure on $[0,1]$), independent of the observables $\boldsymbol{H}$.
In the case of pure policy randomization, $U$ enters only through
the policy $\boldsymbol{p}$ (as general action noise), resulting
in the problem
\begin{equation}
\begin{array}{rl}
\underset{\boldsymbol{x},\boldsymbol{p}(\cdot)}{\mathrm{maximize}} & g^{o}(\boldsymbol{x})\\
\mathrm{subject\,to} & \boldsymbol{x}\le-\boldsymbol{\rho}(-\mathbb{E}_{U}\{\boldsymbol{f}(\boldsymbol{p}(\boldsymbol{H},U),\boldsymbol{H})\})\\
 & \boldsymbol{g}(\boldsymbol{x})\ge{\bf 0}\\
 & (\boldsymbol{x},\boldsymbol{p}(\cdot))\in{\cal X}\times\Gamma(\cdot),\quad\mathrm{P}_{\boldsymbol{S}}\text{-a.e.}
\end{array},\tag{{PureR}}\label{eq:PureR}
\end{equation}
where we have explicitly chosen a decomposable set $\Pi=\{\boldsymbol{p}:\boldsymbol{S}\rightarrow\mathbb{R}^{N_{\boldsymbol{p}}}\hspace{1bp}|\hspace{1bp}\boldsymbol{p}(\boldsymbol{s})\hspace{-2bp}\in\Gamma(\boldsymbol{s}),\mathrm{P}_{\boldsymbol{S}}\text{ -a.e.}\}$,
with $\Gamma:\boldsymbol{S}\rightrightarrows\mathbb{R}^{N_{\boldsymbol{p}}}$
being a measurable multifunction \textit{constant in the second argument},
i.e., $\Gamma(\boldsymbol{s})=\Gamma(\boldsymbol{h}),\boldsymbol{s}\in{\cal S}$.
Note that, aside from constancy, this is already a quite general choice:
By Hiai and Umegaki \cite{HIAI1977149}, we know that a nonempty
closed subset of ${\cal L}_{p<\infty}(\mathrm{P},\mathbb{R}^{N_{\boldsymbol{p}}})$
is decomposable if and only if it coincides with the set containing
the measurable selections of some closed-valued measurable multifunction
(which set is of course decomposable). 
Now, for every policy $\boldsymbol{p}$, we may write
\begin{align*}
\mathbb{E}_{U}\{\boldsymbol{f}(\boldsymbol{p}(\boldsymbol{H},U),\boldsymbol{H})\} & =\int^{1}_{0}\boldsymbol{f}(\boldsymbol{p}(\boldsymbol{H},u),\boldsymbol{H})\mathrm{d}u
=\int\boldsymbol{f}(\boldsymbol{u},\boldsymbol{H})\mathrm{d}\kappa_{\boldsymbol{p}}(\boldsymbol{u}|\boldsymbol{H})\triangleq\langle\boldsymbol{f}(\cdot,\boldsymbol{H}),\kappa_{\boldsymbol{p}}(\cdot|\boldsymbol{H})\rangle,
\end{align*}
which confirms that every possible $\boldsymbol{p}(\boldsymbol{H},\cdot)$
(randomized in $U$) induces a stochastic kernel $\kappa_{\boldsymbol{p}}(\cdot|\boldsymbol{H})$;
in other words, $\kappa_{\boldsymbol{p}}(\cdot|\boldsymbol{H})$ is
precisely (a version of) the conditional distribution of $\boldsymbol{p}$
relative to $\boldsymbol{H}$. What is more, because $U$ is nonatomic,
it is well-known that \textit{the pushforward map that takes $\boldsymbol{p}(\boldsymbol{H},\cdot)$
to $\kappa_{\boldsymbol{p}}(\cdot|\boldsymbol{H})$ is onto (surjective)}.
This fact allows us to recast problem (\ref{eq:PureR}) into the form
\begin{equation}
{\begin{array}{rl}
\underset{\boldsymbol{x},\kappa}{\mathrm{maximize}} & g^{o}(\boldsymbol{x})\\
\mathrm{subject\,to} & \boldsymbol{x}\le-\boldsymbol{\rho}(-\langle\boldsymbol{f}(\cdot,\boldsymbol{H}),\kappa(\cdot|\boldsymbol{H})\rangle)\\
 & \boldsymbol{g}(\boldsymbol{x})\ge{\bf 0}\\
 & \boldsymbol{x}\in{\cal X}
 ,\,
\kappa(\Gamma(\cdot)|\cdot)=1,\,\,\mathrm{P}\text{-a.e.}
\end{array}.}\tag{{KPureR}}\label{eq:KPureR}
\end{equation}
Problems (\ref{eq:PureR}) and (\ref{eq:KPureR}) are equivalent,
even though the feasible set of (\ref{eq:PureR}) is generally larger
(i.e., contains more elements). However, (\ref{eq:KPureR}) is a perfect
instance of convex programming under mere component-wise convexity
of the risk operator $\boldsymbol{\rho}$; just observe that $\boldsymbol{\rho}$
evaluates a linear functional of the Markov kernel $\kappa(\cdot|\boldsymbol{H})$,
whereas the support constraint related to the multifunction $\Gamma$
is also convex (even though the images of $\Gamma$ might be nonconvex
sets).
Consequently, assuming Slater's condition,
strong duality is established for both problems immediately via standard
hyperplane separation. We thus see that at least in the case of pure
randomization discussed above, strong duality follows as a trivial
fact from standard convex programming. 

But the setting we are studying in this paper is by far the most general:
We place \textit{absolutely no assumptions on the ways either the
function $\boldsymbol{f}$ or the decomposable set $\Pi$ can depend
on the joint split state} $\boldsymbol{S}=(\boldsymbol{H},\boldsymbol{U})$.
Both variables of the tuple $(\boldsymbol{H},\boldsymbol{U})$ can
enter $\boldsymbol{f}(\boldsymbol{p}(\cdot,\bullet),\cdot,\bullet)$
and $\Pi$ in a completely arbitrary manner, as long $\boldsymbol{f}$
remains integrable and $\Pi$ remains decomposable. 

As discussed in the next sections, we establish that strong duality
indeed holds regardless of problem-level convexity or, to state it
more concretely, regardless of whether \eqref{eq:Base} can be reduced
to a convex Markov kernel program or not, and one can assert strong
duality for \eqref{eq:Base} totally without the need to check for
convexity of some equivalent kernel problem. In fact, we reveal a
far more general mechanism for certifying strong duality, that applies
seamlessly even in cases where \eqref{eq:Base} cannot be convexified
(at problem level) via the introduction of kernel representations.
More specifically, our analysis shows that \textit{convexification
indeed happens (as perhaps expected), but at the higher ``gestalt''
of image space sets} (and not necessarily at problem level) \cite{FloresBazan2013, Giannessi2005, Bertsekas2009}, as introduced later in Section \ref{sec:Proof-of-Theorems}. We believe that this
is the right way of looking at the problem of certifying strong duality,
since an image space set is by construction the universal object of
interest in relation to the geometry of strong duality relations \cite{FloresBazan2013, Giannessi2005, Bertsekas2009} (again, see Section \ref{sec:Proof-of-Theorems}).


\subsection{Decomposable Policies} \label{subsec: decomposable}

Quite surprisingly, the first set of conditions enabling the dual-domain analysis
of problem \eqref{eq:Base} that we develop is standard and in fact \textit{almost} identical to the set of conditions employed in 
the relevant literature on the corresponding risk-neutral related problems; see, e.g., \citep{Luo2008,Ribeiro2012,Chamon2020a,Chamon2021} and also Section \ref{subsec:Contributions}. We further discuss this in more detail below and later in Section \ref{Counter}.

\noindent %
\noindent\fcolorbox{black}{lightgray}{\begin{minipage}[t]{1\columnwidth - 2\fboxsep - 2\fboxrule}%
\begin{assumption}[\textbf{Problem Setting I}]
\label{assu:Assumption}The following conditions are in effect:
\begin{enumerate}
\item The set $\mathcal{X}\subseteq \mathbb{R}^N$ of performance risk vectors is convex.
\item The utilities $g^{o}$ and $\boldsymbol{g}$ are concave on $\mathcal{X}$.
\item The Borel conditional measure $\mathrm{P}_{\boldsymbol{U}|\boldsymbol{H}}$ is nonatomic $\mathrm{P}$-a.e.\begingroup
  \renewcommand\thempfootnote{\ddag}%
  \footnote{Recall that a measure $\mathrm{P}$ is nonatomic if for any event $E$ with
$\mathrm{P}(E)>0$, an event $E'\subseteq E$ exists such that $\mathrm{P}(E)>\mathrm{P}(E')>0$.}%
\endgroup
.
\item The operator $\boldsymbol{\rho}:{\cal L}_{p}(\mathrm{P},\mathbb{R}^{N})\rightarrow \mathbb{R}^{N}$ satisfies, for every $i\in\mathbb{N}_{N}^{+}$:
\begin{enumerate}
    \item $\rho_i(\boldsymbol{Z}) = \rho_i(Z_i)$, for all $\boldsymbol{Z}\in{\cal L}_{p}(\mathrm{P},\mathbb{R}^{N})$,
    \item $\rho_i$ is convex.
\end{enumerate}
\item The policy feasible set $\Pi$ is decomposable\begingroup
  \renewcommand\thempfootnote{\dag}%
  \footnote{We call $\Pi$ a decomposable space/set if for all $\boldsymbol{p}, \boldsymbol{p}' \in \Pi$ and any event $E \in \mathscr{B}(\mathcal{S})$, we have that $\mathds{1}_{E} \boldsymbol{p} + \mathds{1}_{\mathcal{S}\setminus E} \boldsymbol{p}' \in \Pi$. We note (see, e.g., \cite{HIAI1977149}) that this definition of decomposability is slightly different from the original definition given in \cite{Rockafellar1971Integrals}.}%
\endgroup.
\item The composition $\boldsymbol{f}(\boldsymbol{p}(\cdot),\cdot)$ is in ${\cal L}_{p}(\mathrm{P}_{\boldsymbol{S}},\mathbb{R}^{N})$ 
for all $\boldsymbol{p} \in \Pi$.
\item Problem \textnormal{\eqref{eq:Base}} satisfies Slater's condition (i.e., inequality constraints are
strictly feasible).
\end{enumerate}
\end{assumption}
\end{minipage}}\vspace{6.5bp}

A useful observation is that the conditions comprising Assumption
\ref{assu:Assumption} are fully compatible, and can be easily satisfied
in applications. For instance, the space of all
(Borel-)measurable functions (policies) from ${\cal H}$ to $\mathbb{R}^{N_{\boldsymbol{p}}}$
as well as all spaces of integrable functions ${\cal L}_{p}(\mathrm{P},\mathbb{R}^{N_{\boldsymbol{p}}}),\ p\in[1,\infty]$,
are decomposable \citep{Rockafellar2009VarAn}, and all those examples
can be possible choices for the feasible set $\Pi$. A more specific
choice of a decomposable feasible set $\Pi$, which is relevant
in applications, is the uniform box
\[
\Pi=\big\{\boldsymbol{p}:{\cal S}\rightarrow\mathbb{R}^{N_{\boldsymbol{p}}}\big|\mathrm{ess\,sup}_{\mathrm{P}_{\boldsymbol{S}}}\,\Vert\boldsymbol{p}(\cdot)\Vert_{\infty}\le \gamma\big\},
\]
where $\gamma>0$ is a fixed number, or a more refined
rectangular box reflecting the condition $\mathrm{ess\,sup}_{\mathrm{P}_{\boldsymbol{S}}}\, |p^{i}(\cdot)|\le \gamma^{i},i\in\mathbb{N}_{N_{\boldsymbol{p}}}^{+}$,
where $\gamma^{i}$ are fixed. A
very general choice for a decomposable $\Pi$ is
\[
\Pi=\big\{\boldsymbol{p}:{\cal S}\rightarrow\mathbb{R}^{N_{\boldsymbol{p}}}\big|\boldsymbol{p}(\boldsymbol{s})\in{\Gamma}(\boldsymbol{s}),\mathrm{P}_{\boldsymbol{S}}\text{-a.e.}\big\},
\]
where ${\Gamma}:{\cal S}\rightrightarrows\mathbb{R}^{N_{\boldsymbol{p}}}$ is \textit{any} closed-valued measurable multifunction (\textit{not} necessarily convex-valued). Note that, by Hiai and Umegaki \cite{HIAI1977149} (cf. our discussion in Section \ref{warmup}), we know that a nonempty and
\textit{closed} subset of ${\cal L}_{p<\infty}(\mathrm{P}_{\boldsymbol{S}},\mathbb{R}^{N_{\boldsymbol{p}}})$
is decomposable if and only if it coincides with the set of
the measurable selections of some closed-valued measurable multifunction. 
In our setting, if ${\Gamma}$ is a closed-valued multifunction which is \textit{additionally} closed \citep[Section 9.2.4, Remark 62]{ShapiroLectures_2ND}, then ${\Gamma}$ is measurable and thus admits at least one measurable
selection \citep[Theorem 9.46]{ShapiroLectures_2ND}. In other words,
every feasible $\boldsymbol{p}\in\Pi$ may be taken to be a well-behaved
Borel-measurable selection of ${\Gamma}$.

Regarding the rest of
the conditions of Assumption \ref{assu:Assumption}, all are standard
as mentioned above, except possibly for condition 3 (nonatomicity), especially when compared point-to-point with the corresponding condition enforced in risk-neutral setting \citep{Luo2008,Ribeiro2012,Chamon2020a,Chamon2021}. In particular, condition 3 demands  
nonatomicity of the conditional Borel measure $\mathrm{P}_{\boldsymbol{U}|\boldsymbol{H}}$, which implies nonatomicity of $\mathrm{P}_{\boldsymbol{S}}$ (the standard relevant condition in the risk-neutral setting), but the converse need \textit{not} be the case in general. 
However, if $\mathrm{P}_{\boldsymbol{S}}$ admits a density with respect with the Lebesgue measure, then $\mathrm{P}_{\boldsymbol{U}|\boldsymbol{H}}$ also admits a Lebesgue density and nonatomicity of the two are in fact equivalent conditions. Of course, having a Lebesgue density is a valid
assumption in numerous practical settings (in which case the whole Assumption \ref{assu:Assumption} becomes identical with the correponding conditions enforced in the risk-neutral setting).
Despite the slight mismatch, we emphasize that our main result in this section to be presented after this discussion and developed under Assumption \ref{assu:Assumption} \textit{does} recover the risk-neutral case \textit{exactly}; we discuss this interesting reduction in more detail later in Section \ref{Counter}.  

Overall, nonatomicity is crucial for our developments; if $\mathrm{P}_{\boldsymbol{U}|\boldsymbol{H}}$ contained atoms (i.e., point masses), our results would cease to hold. Nonetheless, even in this case one may artificially enforce nonatomicity by means of an integral convolution of an atomic measure $\mathrm{P}_{\boldsymbol{U}|\boldsymbol{H}}$ with, e.g., a Gaussian probability density function of arbitrarily small variance (i.e., by just adding such small-variance noise to $\boldsymbol{U}$), thus approximating the original (atomic) measure, to an arbitrary accuracy, via a new purely nonatomic measure. Thus, although crucial, the nonatomicity assumption is either met in (a plethora of) practical applications, or can be enforced by appropriately approximating the problem under consideration, with (hopefully) minimal loss of optimality.

We are now in position to state the first main result of this work, proved in Section \ref{sec:Proof-of-Theorem}.
\vspace{1.5bp}

\noindent\fcolorbox{black}{orange}{\begin{minipage}[t]{1\columnwidth - 2\fboxsep - 2\fboxrule}%
\begin{theorem}[\textbf{Strong Duality---Decomposable Case}]
 \label{thm:Main} Let Assumption \textnormal{\ref{assu:Assumption}} be in effect.
Then problem \textnormal{\eqref{eq:Base}} exhibits strong duality, i.e., $\mathsf{P}^{*}=\mathsf{D}^{*}$ and optimal dual
variables exist.
\end{theorem}
\end{minipage}}\vspace{4bp}

\begin{remark} Observe that attainment of the primal solution is neither assumed nor implied by Theorem \textnormal{\ref{thm:Main}}. The result holds irrespectively of whether the primal solution is attained or not. \hfill \qed
\end{remark}

\subsection{Nondecomposable Policies} \label{subsec: densely decomposable sets}

In this section, we extend the strong duality result given in Theorem \ref{thm:Main} for a large class of possibly nondecomposable policy spaces/sets (i.e., $\Pi$), by an appropriate adjustment of the minimal conditions given in Assumption \ref{assu:Assumption}. To this end, we first introduce the notion of a \textit{densely decomposable set}, strictly generalizing the notion of decomposability. Subsequently, we provide a rich (though non-exhaustive) list of examples of such densely decomposable sets, together with a simple recipe for the identification of such sets. We emphasize that the developments in this section are novel even for the risk-neutral counterpart of \eqref{eq:Base} (i.e., when $\boldsymbol{\rho}$ is replaced by a vector of expectations), further pushing the frontiers of nonconvex constrained stochastic programming in this case as well.
\begin{definition}[\textbf{\textit{Densely Decomposable Sets}}] \label{def: densely decomposable set}
 A set of $\mathscr{B}({\cal S})$-measurable functions $\Pi$ is
called \emph{densely decomposable relative to a norm $\Vert\cdot\Vert$}
iff for every pair of elements $\boldsymbol{p}\in\Pi$ and $\boldsymbol{p}'\in\Pi$
and every event $E\in\mathscr{B}({\cal S})$, it is true that
\begin{align*}
\forall\varepsilon>0\quad\exists\boldsymbol{p}_{E}^{\varepsilon}\in\Pi\quad\text{such that}\quad & \Vert\boldsymbol{p}_{E}^{\varepsilon}-[\mathds{1}_{E}\boldsymbol{p}+\mathds{1}_{{\cal S}\backslash E}\boldsymbol{p}']\Vert\le\varepsilon,
\end{align*}
also satisfying $\min\{\Vert\boldsymbol{p}_{E}^{\varepsilon}\Vert,\Vert\mathds{1}_{E}\boldsymbol{p}+\mathds{1}_{{\cal S}\backslash E}\boldsymbol{p}'\Vert\}<\infty$.
\end{definition}

As a consequence of the definition above, an application
of the reverse triangle inequality shows that both $\Vert\boldsymbol{p}_{E}^{\varepsilon}\Vert<\infty$
and $\Vert\mathds{1}_{E}\boldsymbol{p}+\mathds{1}_{{\cal S}\backslash E}\boldsymbol{p}'\Vert<\infty$,
for all $E\in\mathscr{B}({\cal S})$. Also note that decomposability
of the normed vector space of $\mathscr{B}({\cal S})$-measurable
functions $\Pi_{\Vert\cdot\Vert}\triangleq\{\boldsymbol{p}:{\cal S}\rightarrow\mathbb{R}^{N_{\boldsymbol{p}}}\ \vert\ \Vert\boldsymbol{p}\Vert<\infty\}\supseteq\Pi$
(to see the inclusion, choose $E={\cal S}$ above) does \textbf{\textit{not}}
follow as a result of $\Pi$ being densely decomposable;
it merely holds that for every pair $\boldsymbol{p}\in\Pi$ and $\boldsymbol{p}'\in\Pi$
and $E\in\mathscr{B}({\cal S})$, the splice $\mathds{1}_{E}\boldsymbol{p}+\mathds{1}_{{\cal S}\backslash E}\boldsymbol{p}'$
is in $\Pi_{\Vert\cdot\Vert}$. Although densely decomposable sets could be defined in relation to a particular decomposable set (see below), Definition \ref{def: densely decomposable set} reflects a minimal structure required to fully characterize them. 

Of course, decomposable sets are trivially densely decomposable
(relative to any finite norm on $\Pi$, if such norm exists). But
there are many densely decomposable sets which are not even remotely
decomposable. In fact, let $\overline{\Pi}$ be a depomposable normed
space (such as ${\cal L}_{p}(\mathrm{P}_{\boldsymbol{S}},\mathbb{R}^{N_{\boldsymbol{p}}})$, $p\in[1,\infty)$).
Then, we have the lemma that every dense subset $\Pi\subseteq\overline{\Pi}$
is densely decomposable relative to the native norm of $\overline{\Pi}$.
Indeed, if for every pair of policies $\boldsymbol{p}\in\Pi$ and
$\boldsymbol{p}'\in\Pi$ and for every $E\in\mathscr{B}({\cal S})$
we define the policy $\boldsymbol{p}_{E}=\mathds{1}_{E}\boldsymbol{p}+\mathds{1}_{{\cal S}\backslash E}\boldsymbol{p}',$
then it follows that $\boldsymbol{p}_{E}\in\overline{\Pi}$ since
$\overline{\Pi}$ is decomposable and contains both $\boldsymbol{p}$
and $\boldsymbol{p}'$ and, for every $\varepsilon>0$, we can choose
$\boldsymbol{p}_{E}^{\varepsilon}\in\Pi$ such that $\Vert\boldsymbol{p}_{E}^{\varepsilon}-\boldsymbol{p}_{E}\Vert_{\overline{\Pi}}\le\varepsilon$
since $\Pi$ is dense in $\overline{\Pi}$. In this case, we may also
set $\Pi_{\Vert\cdot\Vert}=\overline{\Pi}$. In the other way around
and for the same reasons as above, we also have the lemma that if
the closure of $\Pi$ relative to some norm $\Vert\cdot\Vert$ (finite
on $\Pi$), say $\mathrm{cl}_{\Vert\cdot\Vert}(\Pi),$ is decomposable,
then $\Pi$ must be densely decomposable (relative to $\Vert\cdot\Vert$).

Some explicit examples of densely decomposable sets (or spaces) that
are identifiable via the aforementioned construction are as follows:
\begin{itemize}[leftmargin=.40cm]
\item The vector space of $\mathbb{R}^{N_{\boldsymbol{p}}}$-valued infinitely smooth compactly
supported functions on ${\cal S}$, i.e., $\Pi={\cal C}_{c}^{\infty}({\cal S},\mathbb{R}^{N_{\boldsymbol{p}}})$,
relative to the norm of ${\cal L}_{p}(\mathrm{P}_{\boldsymbol{S}},\mathbb{R}^{N_{\boldsymbol{p}}})=\overline{\Pi}$,
for any choice of $p\in[1,\infty)$.
\item The set of continuous and bounded functions, i.e., $\Pi={\cal C}_b({\cal S},{\mathbb{R}^{N_{\boldsymbol{p}}}})\supset{\cal C}_{c}^{\infty}({\cal S},\mathbb{R}^{N_{\boldsymbol{p}}})$,
relative to the norm of ${\cal L}_{p}(\mathrm{P}_{\boldsymbol{S}},\mathbb{R}^{N_{\boldsymbol{p}}})=\overline{\Pi}$,
for any choice of $p\in[1,\infty)$.
\item The vector space of Schwarz functions, i.e., $\Pi={\mathfrak S}({\cal S},\mathbb{R}^{N_{\boldsymbol{p}}})\supset{\cal C}_{c}^{\infty}({\cal S},\mathbb{R}^{N_{\boldsymbol{p}}})$,
relative to the norm of ${\cal L}_{p}(\mathrm{P}_{\boldsymbol{S}},\mathbb{R}^{N_{\boldsymbol{p}}})=\overline{\Pi}$,
for any choice of $p\in[1,\infty)$ and for  (say) $N_{\boldsymbol{p}} = 1$.
\item Large Reproducing Kernel Hilbert Spaces (RKHS) which are dense in
${\cal L}_{p}(\mathrm{P}_{\boldsymbol{S}},\mathbb{R}^{N_{\boldsymbol{p}}})$, for some $p\in[1,\infty)$.
In this case $\Pi$ is identified by the RKHS, and ${\cal L}_{p}(\mathrm{P}_{\boldsymbol{S}},\mathbb{R}^{N_{\boldsymbol{p}}})=\overline{\Pi}$.
A standard example is the Gaussian RKHS for $N_{\boldsymbol{p}}=1$,
where again any $p\in[1,\infty)$ works (see \cite[Theorem 4.63]{steinwart2008support}).
\item The universal class of \textit{finite-dimensional} real-valued \textit{Two hidden
Layer Feedforward neural Networks (TLFN)} on compact boxes as introduced
in \cite[Theorem 2.1]{GULIYEV2018262}, i.e., $\Pi=\mathrm{TFLN}({\cal B})$, where ${\cal B}\subset{\cal S}$
is a compact box, which is dense in the vector space of (bounded)
continuous functions on ${\cal B}$ metrized through the supremum
norm, say ${\cal C}({\cal B})$. In turn, ${\cal C}({\cal B})$ is
itself dense in ${\cal L}_{p}(\mathrm{P}_{\boldsymbol{S}},\mathbb{R}),p\in[1,\infty)$
for any (probability) measure compactly supported on ${\cal B}$.
Under these circumstances, for every $\varepsilon$ and
$\boldsymbol{p}\in{\cal L}_{p}(\mathrm{P}_{\boldsymbol{S}},\mathbb{R})$, we can find
a $\boldsymbol{p}_{p}^{\varepsilon}\in{\cal C}({\cal B})$ such
that $\Vert\boldsymbol{p}-\boldsymbol{p}_{p}^{\varepsilon}\Vert_{{\cal L}_{p}}\le\varepsilon$.
Also, for the same $\varepsilon>0$ for every such $\boldsymbol{p}_{p}^{\varepsilon}\in{\cal C}({\cal B})$,
we can find $\boldsymbol{p}_{c}^{\varepsilon}\in\mathrm{TFLN}({\cal B})$,
such that
\vspace{-1pt}
\begin{align*}
\varepsilon & \ge\sup_{\boldsymbol{s}\in{\cal B}}|\boldsymbol{p}_{p}^{\varepsilon}(\boldsymbol{s})-\boldsymbol{p}_{c}^{\varepsilon}(\boldsymbol{s})|\ge\mathrm{ess\,sup}_{\mathrm{P}}|\boldsymbol{p}_{p}^{\varepsilon}-\boldsymbol{p}_{c}^{\varepsilon}|\ge\Vert\boldsymbol{p}_{p}^{\varepsilon}-\boldsymbol{p}_{c}^{\varepsilon}\Vert_{{\cal L}_{p}}.
\end{align*}
\vspace{-11pt}

\noindent Therefore, we can write
$
\Vert\boldsymbol{p}-\boldsymbol{p}_{c}^{\varepsilon}\Vert_{{\cal L}_{p}}\le\Vert\boldsymbol{p}-\boldsymbol{p}_{p}^{\varepsilon}\Vert_{{\cal L}_{p}}+\Vert\boldsymbol{p}_{p}^{\varepsilon}-\boldsymbol{p}_{c}^{\varepsilon}\Vert_{{\cal L}_{p}}\le2\varepsilon,
$
showing that $\Pi=\mathrm{TFLN}({\cal B})$ is dense in the decomposable
space ${\cal L}_{p}(\mathrm{P}_{\boldsymbol{S}},\mathbb{R})=\overline{\Pi}$.
\end{itemize}
Exploiting the notion of densely decomposable sets, another set of conditions suitable for analyzing problem \eqref{eq:Base}, however assuming a nondecomposable policy feasible set, is as follows.

\noindent %
\noindent\fcolorbox{black}{lightgray}{\begin{minipage}[t]{1\columnwidth - 2\fboxsep - 2\fboxrule}%
\begin{assumption}[\textbf{Problem Setting II}]
\label{assu:Assumption2} Together with conditions 1--4 and 7 of Assumption \textnormal{\ref{assu:Assumption}}, the following additional conditions are in effect:
\begin{enumerate}
\item [$5'$.] The policy feasible set $\Pi$ is densely-decomposable relative to a norm $\Vert \cdot \Vert$.
\item [$6'$.] The function $\boldsymbol{f}$ is such that:
\begin{enumerate}
    \item $\boldsymbol{f}(\boldsymbol{p}(\cdot),\cdot)\in{\cal L}_{p}(\mathrm{P}_{\boldsymbol{S}},\mathbb{R}^{N})$ 
for all $\boldsymbol{p} \in \Pi_{\Vert \cdot\Vert}$,
\item $-\boldsymbol{\rho}(-[\boldsymbol{F}\boldsymbol{p}]_{\boldsymbol{f}}(\boldsymbol{H}))$ is continuous relative to $\boldsymbol{p}\in\Pi_{\Vert\cdot\Vert}$.
\end{enumerate}

\end{enumerate}
\end{assumption}
\end{minipage}}
\vspace{3bp}

Out of the newly introduced conditions of Assumption \ref{assu:Assumption2},
continuity of $-\boldsymbol{\rho}(-[\boldsymbol{F}\boldsymbol{p}]_{\boldsymbol{f}}(\boldsymbol{H}))$
at $\boldsymbol{p} \in \Pi_{\Vert\cdot\Vert}$ deserves some discussion. First, let $\boldsymbol{f}(\cdot,\boldsymbol{S})$
be $\alpha$-H\"older continuous pointwise in $\boldsymbol{S}=(\boldsymbol{H},U)$, with constants $L(\boldsymbol{S})$. Let $p=1$ for simplicity, and let us take $\Pi_{\Vert\cdot\Vert}={\cal L}_{p'}(\mathrm{P}_{\boldsymbol{S}},\mathbb{R}^{N_{\boldsymbol{p}}}),p'\in[1,\infty)$.
For every pair of policies $\boldsymbol{p}\in\Pi_{\Vert\cdot\Vert}$
and $\boldsymbol{p}'\in\Pi_{\Vert\cdot\Vert}$, we may write (with
$1/q^o+1/p^o=1,$ $p^o<\infty$)
\begin{align}
\Vert\boldsymbol{f}(\boldsymbol{p}(\boldsymbol{S}),\boldsymbol{S})-\boldsymbol{f}(\boldsymbol{p}'(\boldsymbol{S}),\boldsymbol{S})\Vert_{{\cal L}_{1}} & \le\mathbb{E}\{L(\boldsymbol{S})\Vert\boldsymbol{p}(\boldsymbol{S})-\boldsymbol{p}'(\boldsymbol{S})\Vert_{2}^{\alpha}\}\nonumber \\
 & \le\Vert L(\boldsymbol{S})\Vert_{{\cal L}_{q^o}}\big(\mathbb{E}\{\Vert\boldsymbol{p}(\boldsymbol{S})-\boldsymbol{p}'(\boldsymbol{S})\Vert_{2}^{p^o\alpha}\}\big)^{1/p^o}\nonumber \\
 & =\Vert L(\boldsymbol{S})\Vert_{{\cal L}_{q^o}}\Vert\boldsymbol{p}(\boldsymbol{S})-\boldsymbol{p}'(\boldsymbol{S})\Vert_{{\cal L}_{\alpha p^o}}^{\alpha} \nonumber\\
 & \le\Vert L(\boldsymbol{S})\Vert_{{\cal L}_{q^o}}\Vert\boldsymbol{p}(\boldsymbol{S})-\boldsymbol{p}'(\boldsymbol{S})\Vert_{{\cal L}_{p'}}^{\alpha},\nonumber
\end{align}
where the last inequality is valid for 
$\alpha \in (0,1]$ with $\alpha p^o \le p'$. 
Then we are done if each $\rho_{i}$ (i.e., the $i$-th coordinate of
$\boldsymbol{\rho}$) is lower semicontinuous (thus continuous) on ${\cal L}_{p=1}(\mathrm{P},\mathbb{R})$,
for $i\in\mathbb{N}_{N}^{+}$ and
\begin{equation}\nonumber    \Vert[\boldsymbol{F}\boldsymbol{p}]_{\boldsymbol{f}}(\boldsymbol{H})-[\boldsymbol{F}\boldsymbol{p}']_{\boldsymbol{f}}(\boldsymbol{H})\Vert_{{\cal L}_{1}}\le\Vert\boldsymbol{f}(\boldsymbol{p}(\boldsymbol{S}),\boldsymbol{S})-\boldsymbol{f}(\boldsymbol{p}'(\boldsymbol{S}),\boldsymbol{S})\Vert_{{\cal L}_{1}}
\end{equation}
by Jensen. More generally, if we can show that
\begin{equation}
\Vert\boldsymbol{f}(\boldsymbol{p}(\boldsymbol{S}),\boldsymbol{S})-\boldsymbol{f}(\boldsymbol{p}'(\boldsymbol{S}),\boldsymbol{S})\Vert_{{\cal L}_{p}}\le L\Vert\boldsymbol{p}(\boldsymbol{S})-\boldsymbol{p}'(\boldsymbol{S})\Vert^{\alpha},\quad L<\infty,\quad\alpha\in(0,1],\nonumber
\end{equation}
continuity of $-\boldsymbol{\rho}(-[\boldsymbol{F}\boldsymbol{p}]_{\boldsymbol{f}}(\boldsymbol{H}))$
at every $\boldsymbol{p} \in \Pi_{\Vert\cdot\Vert}$ follows (as above). Minimally, mere pointwise ${\cal L}_{p}$-continuity
of $\boldsymbol{f}(\cdot,\boldsymbol{S})$ on $\Pi_{\Vert\cdot\Vert}$
suffices (and it is implied by the H\"older conditions above).

Our second main result now follows.
The proof is deferred to Section \ref{subsec: proof of theorem 2}.

\noindent\fcolorbox{black}{orange}{\begin{minipage}[t]{1\columnwidth - 2\fboxsep - 2\fboxrule}%
\begin{theorem}[\textbf{Strong Duality|Densely Decomposable Case}]
 \label{thm:Main2} Let Assumption \textnormal{\ref{assu:Assumption2}} be in effect.
Then problem \textnormal{\eqref{eq:Base}} exhibits strong duality.
\end{theorem}
\end{minipage}}
\begin{remark}[\textbf{On $\mathrm{TFLN}$ Parametrizations}]
    Although $\mathrm{TFLN}$s is a finite-dimensional densely decomposable class ---which, under Assumption \ref{assu:Assumption2}, ensures strong duality via Theorem \ref{thm:Main2}---, the implementation of such functions might not be particularly computationally appealing \cite{GULIYEV2018262}. In applications, one may combine our strong duality result together with the exploitation of other more standard function approximators, aiming for a trade-off between implementation/computational tractability and a possibly nonzero (but controlled) duality gap in the  optimization task.\hfill \qed
\end{remark}
\section{\label{sec:Proof-of-Theorems}Proofs}

We will work with the image space set \cite{FloresBazan2013, Giannessi2005} associated with problem \eqref{eq:Base} and defined as
\[
{\cal C}\triangleq\left\{ (\delta_{o},\boldsymbol{\delta}_{r},\boldsymbol{\delta}_{d})\left|\begin{array}{l}
g^{o}(\boldsymbol{x})\ge\delta_{o}\\
-\boldsymbol{\rho}(-[\boldsymbol{F}\boldsymbol{p}]_{\boldsymbol{f}}(\boldsymbol{H}))-\boldsymbol{x}\ge\boldsymbol{\delta}_{r}\\
\boldsymbol{g}(\boldsymbol{x})\ge\boldsymbol{\delta}_{d}
\end{array}\hspace{-1bp},\;\text{for some }(\boldsymbol{x},\boldsymbol{p})\in{\cal X}\times\Pi\right.\right\} .
\]
Following \citep{Chamon2020a,Chamon2021}, showing that ${\cal C}$
is convex and the strict feasibility of problem \eqref{eq:Base} (i.e., Slater's condition in
Assumption \ref{assu:Assumption}) would suffice to ensure strong duality of \eqref{eq:Base}, as a relatively simple consequence of the
supporting hyperplane theorem; see Section \ref{subsec:Convexity-of-C},
or \citep[Chapter 4]{Bertsekas2009}, or \citep[Theorem 1 and its proof]{Chamon2020a} 
for the details. In fact, this is a special part of a far more general
story; see, e.g., \citep{FloresBazan2013} and also \citep[Chapter 4]{Bertsekas2009}. Proving convexity of ${\cal C}$
is nontrivial in the case of \eqref{eq:Base} though, and does not
follow from existing standard analysis.
This is due to the nonlinearity of the functionals present in the
risk constraints of \eqref{eq:Base}, in sharp contrast to standard
problems considered in the literature (e.g., in \citep{Luo2008,Ribeiro2012,Chamon2020a,Chamon2021}),
where the corresponding constraints evaluate the vector $\boldsymbol{f}$
solely through linear functionals, i.e., expectations.

\textbf{{The Challenge of Risk:}} From the discussion above
it follows that ultimately we would like to prove that the set ${\cal C}$
is convex. This would mean that if $(\delta_{o},\boldsymbol{\delta}_{r},\boldsymbol{\delta}_{d})\in{\cal C}$
for $(\boldsymbol{x},\boldsymbol{p})\in{\cal X}\times\Pi$ and if
$(\delta'_{o},\boldsymbol{\delta}'_{r},\boldsymbol{\delta}'_{d})\in{\cal C}$
for $(\boldsymbol{x}',\boldsymbol{p}')\in{\cal X}\times\Pi$, then,
for every $\alpha\in[0,1]$, it should be the case that
$
\alpha(\delta_{o},\boldsymbol{\delta}_{r},\boldsymbol{\delta}_{d})+(1-\alpha)(\delta'_{o},\boldsymbol{\delta}'_{r},\boldsymbol{\delta}'_{d})\in{\cal C}
$.
In other words, we would have to show that there exists another pair
$(\boldsymbol{x}_{\alpha},\boldsymbol{p}_{\alpha})\in{\cal X}\times\Pi$,
such that
\begin{align*}
g^{o}(\boldsymbol{x}_{\alpha}) & \ge\alpha\delta_{o}+(1-\alpha)\delta'_{o},\\
-\boldsymbol{\rho}(-[\boldsymbol{F}\boldsymbol{p}_\alpha]_{\boldsymbol{f}}(\boldsymbol{H}))-\boldsymbol{x}_{\alpha} & \ge\alpha\boldsymbol{\delta}_{r}+(1-\alpha)\boldsymbol{\delta}'_{r}\quad\text{and}\\
\boldsymbol{g}(\boldsymbol{x}_{\alpha}) & \ge\alpha\boldsymbol{\delta}_{d}+(1-\alpha)\boldsymbol{\delta}'_{d}.
\end{align*}
By choosing $\boldsymbol{x}_{\alpha}=\alpha\boldsymbol{x}+(1-\alpha)\boldsymbol{x}'\in{\cal X}$
(${\cal X}$ is convex), and appealing to the concavity
of $g^{o}$ and $\boldsymbol{g}$, the proof would be complete if
we showed that, for every $\alpha\in[0,1]$, there is a policy $\boldsymbol{p}_{\alpha}\in\Pi$,
such that
\begin{tagequation}{KEY}
\boxed{
-\boldsymbol{\rho}(-[\boldsymbol{F}\boldsymbol{p}_\alpha]_{\boldsymbol{f}}(\boldsymbol{H}))\ge-\alpha\boldsymbol{\rho}(-[\boldsymbol{F}\boldsymbol{p}]_{\boldsymbol{f}}(\boldsymbol{H}))-(1-\alpha)\boldsymbol{\rho}(-[\boldsymbol{F}\boldsymbol{p}']_{\boldsymbol{f}}(\boldsymbol{H})).}\label{eq:The_Key}
\end{tagequation}

To prove Theorem \ref{thm:Main}, we indeed establish {convexity of the image set ${\cal C}$}, verifying \eqref{eq:The_Key} directly; this nontrivial fact suffices to establish strong duality of \eqref{eq:Base} under Assumption \ref{assu:Assumption}.
However, for proving Theorem \ref{thm:Main2}, what we can show is merely \textit{convexity of
the closure of ${\cal C}$} which, despite somewhat nonstandard, is still sufficient to establish strong duality. As will be discussed in due course, convexity
of \textit{$\mathrm{cl}({\cal C})$ }is combined with the fact
that the point $(\mathsf{P}^{*},{\bf 0},{\bf 0})$ cannot be in the
interior of $\mathrm{cl}({\cal C})$, enabling a separating hyperplane
argument for the pair $(\mathsf{P}^{*},{\bf 0},{\bf 0})$ and $\mathrm{cl}({\cal C})$,
which together with Slater's condition (i.e., strict feasibility) implies strong duality for
\eqref{eq:Base} under Assumption \ref{assu:Assumption2}; this step is similar to but extends the case where ${\cal C}$ is a convex set.
Still, we believe that establishing strong duality for \eqref{eq:Base} using
merely convexity of \textit{$\mathrm{cl}({\cal C})$ }is interesting,
as it reveals a weaker chain of self-contained logical
arguments sufficient to prove Theorem \ref{thm:Main2} as compared
with Theorem \ref{thm:Main} or the risk-neutral setting; see also \citep[Chapter 4]{Bertsekas2009}. For an alternative
(though non-elementary and more involved) way of establishing strong
duality for \eqref{eq:Base} ---again by exploiting convexity of \textit{$\mathrm{cl}({\cal C})$}
and Slater's condition---, one can exploit more advanced results, e.g., \citep[Theorem 3.2]{FloresBazan2013}.

In passing, note that \eqref{eq:The_Key} is trivially true and thus
convexity of ${\cal C}$ is implied in the special case where the composition $-\boldsymbol{\rho}\left(-[\boldsymbol{F}\boldsymbol{p}]_{\boldsymbol{f}}(\boldsymbol{H})\right)$
is concave in $\boldsymbol{p} \in \Pi$ and $\Pi$ \textit{is
a convex set} (instead of a (densely) decomposable one), simply by choosing
$\boldsymbol{p}_{\alpha}$ as a convex combination of $\boldsymbol{p}$
and $\boldsymbol{p}'$; this is a standard
case of convex programming. In particular, this holds whenever $\boldsymbol{\rho}$ is convex and monotone, and $\boldsymbol{f}(\cdot,\boldsymbol{s})$ is concave for all $\boldsymbol{s} \in \mathcal{S}$.  

\subsection{Preliminaries on Vector Measures of Infinite Dimension}

Suppose that $\mathds{X}$ is a possibly infinite-dimensional space, e.g., a Banach space or a product space $\mathbb{R}^\Theta$, where $\Theta$ is a general set (possibly uncountably infinite). A \textit{vector measure} on a measurable space $(\Omega,\mathscr{F})$
is a countably additive function $\boldsymbol{G}:\mathscr{F}\rightarrow\mathds{X}$. Hereafter, we let $\mathds{X}=\mathbb{R}^\Theta$ with $\Theta$ as above. We say that a vector measure $\boldsymbol{G}:\mathscr{F}\rightarrow\mathds{X}$ \textit{admits a density} whenever there is a common finite measure $\mu :\mathscr{F}\rightarrow\mathbb{R}_+$, with respect to which each component of the vector $\boldsymbol{G}(\cdot)$ is an absolutely continuous finite (signed) measure on $(\Omega, \mathscr{F})$. In other words, there exists a vector \textit{density function} $\boldsymbol{g}=[g_\theta \in \mathcal{L}_1(\mu,\mathbb{R}) ]_{\theta \in \Theta} \in \mathds{X}$ such that
\begin{equation}\nonumber
    \boldsymbol{G}(E)=\bigg[ 
    \int_E g_\theta\, \mathrm{d}\mu 
    \bigg]_{\theta \in \Theta},\quad E\in\mathscr{F}.
\end{equation}
For any measurable set $E\in\mathscr{F}$, the range of $\boldsymbol{G}$ restricted on $E$ is defined as
\begin{equation}\nonumber
    \boldsymbol{G}(\mathscr{F}_E)\triangleq \{\boldsymbol{x}\in\mathds{X}\ \vert\ \boldsymbol{x}=\boldsymbol{G}(E),\text{ for some }E\in\mathscr{F}_E\},
\end{equation}
where $\mathscr{F}_E \triangleq \{E'\in\mathscr{F}|E' \subseteq E\}$. Lastly, let us also define the following key concept of from \cite[Section 2]{KingmanRobertson}: A set of functions $\mathcal{M}\subset \mathcal{L}_1(\mu,\mathbb{R})$ is called \textit{thin} if and only if  its ``restricted kernel"
\begin{equation}\nonumber
    \mathcal{M}^\bot(E)\triangleq \bigg\{
    \phi\in \mathcal{L}_\infty(\mu,\mathbb{R})
    \ \bigg\vert\ 
    \phi \mathds{1}_{\Omega\backslash E}\equiv 0 \text{ and }
    \int \phi g\, \mathrm{d}\mu = 0, \text{ for all } g \in \mathcal{M}   
    \bigg\},\quad E\in\mathscr{F}
\end{equation}
satisfies $\mathcal{M}^\bot(E) \neq \{0\}$ for every $E\in\mathscr{F}$ with $\mu(E)>0$. In words, $\mathcal{M}$ is thin iff $\mathcal{M}^\bot(\cdot)$ is different from the zero subspace on all  events of positive measure. 

Our analysis below leverages the discussion above and will be based on the following generalization of
the celebrated convexity theorem of A. A. Lyapunov, due to Kingman and Robertson \citep[Theorem 1 and Corollary]{KingmanRobertson}, later significantly extended by Knowles \citep[Theorem 1]{SIAMCon:Knowles}. This also classical
result \textit{characterizes} the range of an infinite dimensional vector measure, by providing necessary \textit{and} sufficient conditions for convexity and compactness (in a suitable topology).
\vspace{6pt}
\begin{theorem}[\textbf{\citep{KingmanRobertson} Kingman-Robertson Characterization}]
\label{thm:Weak_Lyapunov} Let $(\Omega,\mathscr{F})$ be a measurable
space, and let $\mathds{X}\triangleq \mathbb{R}^\Theta$ endowed with its product topology (where $\Theta$ is a general set). Let $\boldsymbol{G}:\mathscr{F}\rightarrow\mathds{X}$
be a vector measure which admits a density function $\boldsymbol{g}$. Then, the ranges $\boldsymbol{G}(\mathscr{F}_E),E\in\mathscr{F}$ are each convex and compact if and only if the set of derivatives $\mathcal{M}\triangleq \{g_\theta\}_{\theta \in \Theta} \subset \mathcal{L}_1(\mu,\mathbb{R})$ is thin.
\end{theorem}

\subsection{\label{sec:Proof-of-Theorem}Proof of Theorem \ref{thm:Main}}
We split the proof of the result into two subsections, as follows. The first establishes convexity of the image space set $\mathcal{C}$ as a consequence of range convexity of a properly defined vector measure of uncountably infinite dimension, while the second establishes range convexity of the vector measure itself. Throughout this section, we let Assumption \ref{assu:Assumption} be in effect.

\subsubsection[Convexity of C]{Convexity of $\boldsymbol{\mathcal{C}}$} \label{subsec:Core}

Consider product space $\mathds{X}=\mathbb{R}^{\{1,2\}\times\mathbb{N}^{+}_{N}\times{\cal L}_{q}(\mathrm{P},\mathbb{R})}$
of uncountably infinite dimension (with its product topology). For
an arbitrary and fixed pair of policies $\boldsymbol{p}_{1}\in\Pi$
and $\boldsymbol{p}_{2}\in\Pi$, we define a \textit{vector measure}
\textit{of} \textit{universal index} $\boldsymbol{G}:\mathscr{B}({\cal S})\rightarrow\mathds{X}\triangleq\mathbb{R}^{\{1,2\}\times\mathbb{N}^{+}_{N}\times{\cal L}_{q}(\mathrm{P},\mathbb{R})}$ as
\[
\boldsymbol{G}(E)\triangleq\bigg[\int_{E}\zeta(\boldsymbol{h})f_{i}(\boldsymbol{p}_{j}(\boldsymbol{h},\boldsymbol{u}),\boldsymbol{h},\boldsymbol{u})\mathrm{d}\mathrm{P}_{\boldsymbol{S}}(\boldsymbol{h},\boldsymbol{u})\bigg]_{(j,i,\zeta)\in\{1,2\}\times\mathbb{N}^{+}_{N}\times{\cal L}_{q}(\mathrm{P},\mathbb{R})},\;E\in\mathscr{B}({\cal S}),
\]
and we also consider the coordinate-restricted vector measure $\boldsymbol{G}_{\boldsymbol{p}}:\mathscr{B}({\cal S})\rightarrow\mathbb{R}^{\mathbb{N}^{+}_{N}\times{\cal L}_{q}(\mathrm{P},\mathbb{R})}$
obtained by fixing the choice of $\boldsymbol{p}_{j}$ in $\boldsymbol{G}$
to some single particular $\boldsymbol{p}$ (e.g., either $\boldsymbol{p}_{1}$
or $\boldsymbol{p}_{2}$). The next key result is concerned with the
geometry of the range of $\boldsymbol{G}$.

\noindent\fcolorbox{black}{lightgray}{\begin{minipage}[t]{1\columnwidth - 2\fboxsep - 2\fboxrule}%
\begin{lemma}[\textbf{Universal Range Convexity}]
\label{thm:Weak_Lyapunov_G} The range of $\boldsymbol{G}$
\[
\boldsymbol{G}(\mathscr{B}({\cal S}))=\{\boldsymbol{x}\in\mathds{X}|\boldsymbol{x}=\boldsymbol{G}(E),\text{ for some }E\in\mathscr{B}({\cal S})\}
\]
 is convex and compact (in product topology).
\end{lemma}
\end{minipage}}

\medskip{}
Before providing a proof for Lemma \ref{thm:Weak_Lyapunov_G} (Section \ref{subsec: proof of vector measure range convexity}), let
us see how its conclusions, in particular convexity the range of $\boldsymbol{G}$,
are sufficient to ensure convexity of ${\cal C}$. To see this, let
us first consider the range $\boldsymbol{G}(\mathscr{B}({\cal S}))$
of $\boldsymbol{G}$. Of course $\boldsymbol{z}=\boldsymbol{G}({\cal S})$
and $\boldsymbol{z}'=\boldsymbol{G}(\emptyset)={\bf 0}$ are both
elements of $\boldsymbol{G}(\mathscr{B}({\cal S}))$. Therefore, Lemma
\ref{thm:Weak_Lyapunov_G} implies that, for every $\alpha\in[0,1]$,
the convex combination $\alpha\boldsymbol{z}+(1-\alpha)\boldsymbol{z}'\equiv\alpha\boldsymbol{z}$
lies in $\boldsymbol{G}(\mathscr{B}({\cal S}))$. In order words,
for each $\alpha\in[0,1]$, there exists a single event $E^{\alpha}\in\mathscr{B}({\cal S})$
such that \[
\boxed{\alpha\boldsymbol{z}=\boldsymbol{G}(E^{\alpha})\iff(1-\alpha)\boldsymbol{G}({\cal S})=\boldsymbol{G}((E^{\alpha})^{c}).}
\]
This implies that
$
\alpha\boldsymbol{G}_{\boldsymbol{p}_{1}}({\cal S})=\boldsymbol{G}_{\boldsymbol{p}_{1}}(E^{\alpha})
$, and symmetrically,
$
(1-\alpha)\boldsymbol{G}_{\boldsymbol{p}_{2}}({\cal S})=\boldsymbol{G}_{\boldsymbol{p}_{2}}((E^{\alpha})^{c})
$.
Now, we may define the splice
\[
\boldsymbol{p}_{\alpha}(\boldsymbol{s})=\mathds{1}_{E^{\alpha}}(\boldsymbol{s})\boldsymbol{p}_{1}(\boldsymbol{s})+\mathds{1}_{{\cal S}\backslash E^{\alpha}}(\boldsymbol{s})\boldsymbol{p}_{2}(\boldsymbol{s}), \quad \boldsymbol{s}\in \mathcal{S}.
\]
It holds that $\boldsymbol{p}_{\alpha}\in\Pi$ because $\Pi$ is decomposable.
Then, it follows that
\begin{align*}
\boldsymbol{G}_{\boldsymbol{p}_{\alpha}}({\cal S}) & =\boldsymbol{G}_{\boldsymbol{p}_{\alpha}}(E^{\alpha})+\boldsymbol{G}_{\boldsymbol{p}_{\alpha}}((E^{\alpha})^{c})=\boldsymbol{G}_{\boldsymbol{p}_{1}}(E^{\alpha})+\boldsymbol{G}_{\boldsymbol{p}_{2}}((E^{\alpha})^{c}),
\end{align*}
which is the same as
\[
\boldsymbol{G}_{\boldsymbol{p}_{\alpha}}({\cal S})=\alpha\boldsymbol{G}_{\boldsymbol{p}_{1}}({\cal S})+(1-\alpha)\boldsymbol{G}_{\boldsymbol{p}_{2}}({\cal S}).
\]
Therefore, we have that, for every $(i,\zeta)\in\mathbb{N}^{+}_{N}\times{\cal L}_{q}(\mathrm{P},\mathbb{R})$,
\begin{align*}
 & \hspace{-1bp}\hspace{-1bp}\hspace{-1bp}\hspace{-1bp}\hspace{-1bp}\hspace{-1bp}\hspace{-1bp}\hspace{-1bp}\hspace{-1bp}\hspace{-1bp}\hspace{-1bp}\hspace{-1bp}\mathbb{E}\{\zeta(\boldsymbol{H})f_{i}(\boldsymbol{p}_{\alpha}(\boldsymbol{S}),\boldsymbol{S})\}
 =\alpha\mathbb{E}\{\zeta(\boldsymbol{H})f_{i}(\boldsymbol{p}_{1}(\boldsymbol{S}),\boldsymbol{S})\}+(1-\alpha)\mathbb{E}\{\zeta(\boldsymbol{H})f_{i}(\boldsymbol{p}_{2}(\boldsymbol{S}),\boldsymbol{S})\},
\end{align*}
which by Fubini implies that
\begin{equation}\nonumber
    \mathbb{E}\{\zeta(\boldsymbol{H})[\boldsymbol{F}\boldsymbol{p}_{\alpha}]^{i}_{\boldsymbol{f}}(\boldsymbol{H})\}=\mathbb{E}\{\zeta(\boldsymbol{H})[\alpha[\boldsymbol{F}\boldsymbol{p}_{1}]^{i}_{\boldsymbol{f}}(\boldsymbol{H})+(1-\alpha)[\boldsymbol{F}\boldsymbol{p}_{2}]^{i}_{\boldsymbol{f}}(\boldsymbol{H})]\}.
\end{equation}
Since $\zeta$ can be any element of ${\cal L}_{q}(\mathrm{P},\mathbb{R})$ separating ${\cal L}_{p}(\mathrm{P},\mathbb{R})$, the display above implies that
\begin{equation}\nonumber
[\boldsymbol{F}\boldsymbol{p}_{\alpha}]_{\boldsymbol{f}}(\boldsymbol{H})=\alpha[\boldsymbol{F}\boldsymbol{p}_{1}]_{\boldsymbol{f}}(\boldsymbol{H})+(1-\alpha)[\boldsymbol{F}\boldsymbol{p}_{2}]_{\boldsymbol{f}}(\boldsymbol{H}),\quad\mathrm{P}\text{-a.e.}
\end{equation}
Therefore, it follows that
\[
-\boldsymbol{\rho}(-[\boldsymbol{F}\boldsymbol{p}_{\alpha}]_{\boldsymbol{f}}(\boldsymbol{H}))=-\boldsymbol{\rho}(-\alpha[\boldsymbol{F}\boldsymbol{p}_{1}]_{\boldsymbol{f}}(\boldsymbol{H})-(1-\alpha)[\boldsymbol{F}\boldsymbol{p}_{2}]_{\boldsymbol{f}}(\boldsymbol{H})),
\]
which from component-wise convexity of $\boldsymbol{\rho}$ yields
\[
-\boldsymbol{\rho}(-[\boldsymbol{F}\boldsymbol{p}_{\alpha}]_{\boldsymbol{f}}(\boldsymbol{H}))\ge-\alpha\boldsymbol{\rho}(-[\boldsymbol{F}\boldsymbol{p}_{1}]_{\boldsymbol{f}}(\boldsymbol{H}))-(1-\alpha)\boldsymbol{\rho}(-[\boldsymbol{F}\boldsymbol{p}_{2}]_{\boldsymbol{f}}(\boldsymbol{H})).
\]
This shows that the image space set ${\cal C}$ is convex, which
immediately establishes strong duality of problem (\ref{eq:Base})
(in light of Slater's condition, after a standard separation argument).\hfill{}$\blacksquare$

\subsection{Proof of Lemma \ref{thm:Weak_Lyapunov_G}} \label{subsec: proof of vector measure range convexity}

We work towards an application of the Kingman-Robertson theorem (cf. Lemma \ref{thm:Weak_Lyapunov_G}) on the vector measure $\boldsymbol{G}$, as defined
in the previous subsection. For brevity, define functions $F_{ij}(\cdot)\triangleq f_{i}(\boldsymbol{p}_{j}(\cdot),\cdot),(j,i)\in\{1,2\}\times\mathbb{N}^{+}_{N}$.
Because the (regular) conditional distribution (stochastic kernel/random
Borel measure) $\mathrm{P}_{\boldsymbol{U}|\boldsymbol{H}}$ is nonatomic
almost everywhere relative to $\mathrm{P}$, there exists for each
$\boldsymbol{h}\in{\cal H}$ a measure isomorphism $\Phi_{\boldsymbol{h}}:{\cal U}\rightarrow[0,1]$,
jointly measurable when seen as a function in $(\boldsymbol{h},\boldsymbol{u})$,
such that
\begin{align*}
\int_{E}\zeta(\boldsymbol{h})F_{ij}(\boldsymbol{h},\boldsymbol{u})\mathrm{d}\mathrm{P}_{\boldsymbol{S}}(\boldsymbol{h},\boldsymbol{u}) & =\int\mathds{1}_{E}(\boldsymbol{h},\boldsymbol{u})\zeta(\boldsymbol{h})F_{ij}(\boldsymbol{h},\boldsymbol{u})\mathrm{d}\mathrm{P}_{\boldsymbol{S}}(\boldsymbol{h},\boldsymbol{u})\\
 & =\int_{{\cal H}\times[0,1]}\mathds{1}_{E}(\boldsymbol{h},\Phi^{-1}_{\boldsymbol{h}}(u))\zeta(\boldsymbol{h})F_{ij}(\boldsymbol{h},\Phi^{-1}_{\boldsymbol{h}}(u))[\mathrm{d}\mathrm{P}(\boldsymbol{h})\hspace{-1bp}\hspace{-1bp}\times\hspace{-1bp}\hspace{-1bp}\mathrm{d}u]\\
 & =\int_{\{(\boldsymbol{h},u)\in{\cal H}\times[0,1]  \vert (\boldsymbol{h},\Phi^{-1}_{\boldsymbol{h}}(u))\in E\}}\zeta(\boldsymbol{h})F_{ij}(\boldsymbol{h},\Phi^{-1}_{\boldsymbol{h}}(u))[\mathrm{d}\mathrm{P}(\boldsymbol{h})\hspace{-1bp}\hspace{-1bp}\times\hspace{-1bp}\hspace{-1bp}\mathrm{d}u],
\end{align*}
where the preimage $\{\boldsymbol{h},u|(\boldsymbol{h},\Phi^{-1}_{\boldsymbol{h}}(u))\in E\}$
is in $\mathscr{B}({\cal H}\hspace{-1bp}\hspace{-0bp}\times\hspace{-0bp}\hspace{-1bp}[0,1])$.
Together with a symmetric argument (fully utilizing the isomorphism $\Phi_{\boldsymbol{h}}$), this implies that for the full range of $\boldsymbol{G}$
it holds that
\[
\boldsymbol{G}(\mathscr{B}({\cal S}))=\{\boldsymbol{x}\in\mathds{X}|\boldsymbol{x}=\widetilde{\boldsymbol{G}}(E),\text{ for some }E\in\mathscr{B}({\cal H}\hspace{-2bp}\times\hspace{-2bp}[0,1])\}=\widetilde{\boldsymbol{G}}(\mathscr{B}({\cal H}\hspace{-1bp}\hspace{-1bp}\times\hspace{-1bp}\hspace{-1bp}[0,1])),
\]
where
\[
\widetilde{\boldsymbol{G}}(E)\triangleq\bigg[\int_{E}\zeta(\boldsymbol{h})F_{ij}(\boldsymbol{h},\Phi^{-1}_{\boldsymbol{h}}(u))[\mathrm{d}\mathrm{P}(\boldsymbol{h})\hspace{-1bp}\hspace{-1bp}\times\hspace{-1bp}\hspace{-1bp}\mathrm{d}u]\bigg]_{(j,i,\zeta)\in\{1,2\}\times\mathbb{N}^{+}_{N}\times{\cal L}_{q}(\mathrm{P},\mathbb{R})},\;E\in\mathscr{B}({\cal H}\hspace{-1bp}\hspace{-1bp}\times\hspace{-1bp}\hspace{-1bp}[0,1]).
\]
Therefore, we may equivalently work on the vector measure $\widetilde{\boldsymbol{G}}$.
To this end, let us first observe that by construction $\widetilde{\boldsymbol{G}}$
admits a density function taking values in the product space $\mathds{X}$, 
which is identified by the set of Radon-Nikodym derivatives
\[
{\cal M}\triangleq\{\zeta(\bullet)F_{ij}(\bullet,\Phi^{-1}_{(\bullet)}(\cdot))\}_{(j,i,\zeta)\in\{1,2\}\times\mathbb{N}^{+}_{N}\times{\cal L}_{q}(\mathrm{P},\mathbb{R})}.
\]
Consequently, the conclusions of Theorem \ref{thm:Weak_Lyapunov} for  $\widetilde{\boldsymbol{G}}$
(and thus for $\boldsymbol{G}(\mathscr{B}(\mathcal{S}))$) will be immediate if we show that
${\cal M}$ is a thin set. Specifically, if $E\in\mathscr{B}({\cal H}\hspace{-1bp}\hspace{-1bp}\times\hspace{-1bp}\hspace{-1bp}[0,1])$
is an event of positive probability, we should find a nonzero function
$\phi\equiv\phi_{E}\in{\cal L}_{\infty}(\mathrm{P}\hspace{-1bp}\hspace{-1bp}\times\hspace{-1bp}\hspace{-1bp}\lambda,\mathbb{R})$,
vanishing outside $E$ and such that
\[
\int\phi(\boldsymbol{h},u)[\zeta(\boldsymbol{h})F_{ij}(\boldsymbol{h},\Phi^{-1}_{\boldsymbol{h}}(u))][\mathrm{d}\mathrm{P}(\boldsymbol{h})\hspace{-1bp}\hspace{-1bp}\times\hspace{-1bp}\hspace{-1bp}\mathrm{d}u]=0,\quad\forall(j,i,\zeta)\in\{1,2\}\times\mathbb{N}^{+}_{N}\times{\cal L}_{q}(\mathrm{P},\mathbb{R}).
\]
By Fubini, we know that, for every $\phi\in{\cal L}_{\infty}(\mathrm{P}\hspace{-1bp}\hspace{-1bp}\times\hspace{-1bp}\hspace{-1bp}\lambda,\mathbb{R})$
and all tuples $(j,i,\zeta)$,
\[
\int\phi(\boldsymbol{h},u)[\zeta(\boldsymbol{h})F_{ij}(\boldsymbol{h},\Phi^{-1}_{\boldsymbol{h}}(u))][\mathrm{d}\mathrm{P}(\boldsymbol{h})\hspace{-1bp}\hspace{-1bp}\times\hspace{-1bp}\hspace{-1bp}\mathrm{d}u]=\int\zeta(\boldsymbol{h})\bigg[\int^{1}_{0}\phi(\boldsymbol{h},u)F_{ij}(\boldsymbol{h},\Phi^{-1}_{\boldsymbol{h}}(u))\mathrm{d}u\bigg]\mathrm{d}\mathrm{P}(\boldsymbol{h}),
\]
where, additionally, it is easy to see that the integrand inside the
brackets on the right-hand side is in ${\cal L}_{p}(\mathrm{P},\mathbb{R})$
(by Jensen). Since $\zeta$ varies freely in ${\cal L}_{q}(\mathrm{P},\mathbb{R})$,
and ${\cal L}_{q}(\mathrm{P},\mathbb{R})$ separates points in ${\cal L}_{p}(\mathrm{P},\mathbb{R})$
(via the canonical bilinear pairing, for all possible Hölder pairs
$(p,q)$), we obtain the \textit{equivalent condition} 
\[
\int^{1}_{0}\phi(\boldsymbol{h},u)F_{ij}(\boldsymbol{h},\Phi^{-1}_{\boldsymbol{h}}(u))\mathrm{d}u=0,\quad\forall(j,i)\in\{1,2\}\times\mathbb{N}^{+}_{N},\quad\mathrm{P}\text{-a.e.}
\]
We next show that, for \textit{any} selection of functions $\{z_{i}\in{\cal L}_{1}(\mathrm{P},\mathbb{R})\}_{i\in\mathbb{N}^{+}_{N}}$
with \textit{any} $N$ finite (thus also applicable to $2N$ above),
and for every non-negligible $E$, a nonzero \textit{annihilator}
$\phi\in{\cal L}_{\infty}(\mathrm{P}\hspace{-1bp}\hspace{-1bp}\times\hspace{-1bp}\hspace{-1bp}\lambda,\mathbb{R})$
vanishing outside $E$ such that
\begin{equation}
\boxed{\int^{1}_{0}\phi(\boldsymbol{h},u)z_{i}(\boldsymbol{h},u)\mathrm{d}u=0,\quad\forall i\in\mathbb{N}^{+}_{N},\quad\mathrm{P}\text{-a.e.}}\tag{{Kernel}}\label{eq:Ortho}
\end{equation}
can be constructed explicitly.

Let us consider any set $E\in\mathscr{B}({\cal H}\times[0,1])$ of
positive joint probability. Since $\phi$ must be vanishing outside
$E$, it is must necessarily be true that
\[
\phi(\boldsymbol{h},u)\equiv\phi(\boldsymbol{h},u)\mathds{1}_{E(\boldsymbol{h})}(u),\quad(\boldsymbol{h},u)\in{\cal H}\times[0,1],
\]
where $E(\boldsymbol{h})\triangleq\{u\in[0,1]|(\boldsymbol{h},u)\in E\}\in\mathscr{B}([0,1])$
is the parameterized $\boldsymbol{h}$-section of $E$ (an elementary
measurable multifunction). Therefore, a general construction for $\phi$
reads
\[
\phi(\boldsymbol{h},u)=\phi'(\boldsymbol{h},u)\mathds{1}_{E(\boldsymbol{h})}(u),\quad(\boldsymbol{h},u)\in{\cal H}\times[0,1],
\]
where $\phi'\in{\cal L}_{\infty}(\mathrm{P}\hspace{-1bp}\hspace{-1bp}\times\hspace{-1bp}\hspace{-1bp}\lambda,\mathbb{R})$
is an appropriate non-zero function that we should choose with the
properties that we want (possibly minus the restricted support, which
is supplied by the indicator). Also observe that $E(\boldsymbol{h})$
must be of strictly positive Lebesgue measure on a strictly positive
$\mathrm{P}$-probability event; indeed, if this was not true, then
$E(\boldsymbol{h})$ would be of zero Lebesgue measure almost everywhere
relative to $\mathrm{P}$, which would imply that $E$ is a null set. 
Hereafter, let us call the \textit{maximal} such strictly positive
$\mathrm{P}$-probability event $H^{E}_{+}\in\mathscr{B}({\cal H})$.
Of course, on the complement of $H^{E}_{+}$ we have $\lambda(E(\boldsymbol{h}))=0$.
Under these provisions, we may rewrite (\ref{eq:Ortho}) as
\[
\int_{E(\boldsymbol{h})}\phi'(\boldsymbol{h},u)z_{i}(\boldsymbol{h},u)\mathrm{d}u=0,\quad\forall i\in\mathbb{N}^{+}_{N},\quad\boldsymbol{h}\in H^{E}_{+},
\]
with the only additional constraint that $\phi'\neq0$ (on $E(\boldsymbol{h}),\boldsymbol{h}\in H^{E}_{+}$)
must reside in ${\cal L}_{\infty}(\mathrm{P}\hspace{-1bp}\hspace{-1bp}\times\hspace{-1bp}\hspace{-1bp}\lambda,\mathbb{R})$.
Note that outside $H^{E}_{+}$ we can define $\phi'$ arbitrarily;
thus we may set $\phi'=0$ on that event.

At this point, we may observe that, for each $\boldsymbol{h}\in H^{E}_{+},$
we have a system of $N$ equations and \textit{one unkown of infinite
dimension} (since the Lebesgue measure is nonatomic). Therefore it
must be expected that the nullspace of this system is not trivial,
at least pointwise in $\boldsymbol{h}\in H^{E}_{+}$. The idea is
to properly parameterize $\phi$ so that we end up with a concrete
construction which will be jointly measurable and in ${\cal L}_{\infty}(\mathrm{P}\hspace{-1bp}\hspace{-1bp}\times\hspace{-1bp}\hspace{-1bp}\lambda,\mathbb{R})$.
Since we have $N$ equations, a natural choice here is to take 
\[
\phi'(\boldsymbol{h},u)\triangleq\sum^{N+1}_{j=1}c_{j}(\boldsymbol{h})\mathds{1}_{L_{j}(\boldsymbol{h})}(u),
\]
for some measurable coefficients $c_{j}$ and multifunctions $L_{j}$,
$j\in\mathbb{N}^{+}_{N+1}$, to be determined. The motivation is to
create an appropriate classical underdetermined system of $N$ equations
and at least $N+1$ uknowns, which will always exhibit a nontrivial
kernel. 

Since the Lebesgue measure is nonatomic (and thus its cdf is continuous),
we can write, for every $\boldsymbol{h}\in H^{E}_{+}$ and $i\in\mathbb{N}^{+}_{N}$,
\begin{align*}
\int_{E(\boldsymbol{h})}\phi'(\boldsymbol{h},u)z_{i}(\boldsymbol{h},u)\mathrm{d}u & =\lambda(E(\boldsymbol{h}))\int^{1}_{0}\phi'(\boldsymbol{h},u)z_{i}(\boldsymbol{h},u)\dfrac{\mathds{1}_{E(\boldsymbol{h})}(u)}{\lambda(E(\boldsymbol{h}))}\mathrm{d}u\\
 & =\lambda(E(\boldsymbol{h}))\int^{1}_{0}\phi'(\boldsymbol{h},u)z_{i}(\boldsymbol{h},u)\mathrm{d}F_{{\mathrm U}[E(\boldsymbol{h})]}(u)\\
 & =\lambda(E(\boldsymbol{h}))\int^{1}_{0}\phi'(\boldsymbol{h},F^{-1}_{{\mathrm U}[E(\boldsymbol{h})]}(u))z_{i}(\boldsymbol{h},F^{-1}_{{\mathrm U}[E(\boldsymbol{h})]}(u))\mathrm{d}u,
\end{align*}
where $F_{{\mathrm U}[E(\boldsymbol{h})]}:\mathbb{R}\rightarrow[0,1]$
is the standard uniform cdf (cf. Lebesgue measure on $[0,1]$) restricted
to the event $E(\boldsymbol{h})$ (cf. normalized Lebesgue measure on
$E(\boldsymbol{h})$), and $F^{-1}_{{\mathrm U}[E(\boldsymbol{h})]}:[0,1]\rightarrow[0,1]$
denotes \textit{any} measurable selection of the standard inverse
(quantile) multifunction of $F_{{\mathrm U}[E(\boldsymbol{h})]}$ (for
every fixed $\boldsymbol{h}$). The last equality in fact reveals \textit{how}
one can choose the multifunctions $L_{j}(\boldsymbol{h}),j\in\mathbb{N}^{+}_{N+1}$
as guaranteed nontrivial intervals (thus with positive Lebesgue measure).
By regularly partitioning the \textit{range} of $F_{{\mathrm U}[E(\boldsymbol{h})]}$
(i.e., the set $[0,1]$), we may measurably set
\[
L_{j}(\boldsymbol{h})\triangleq\bigg(F^{-1}_{{\mathrm U}[E(\boldsymbol{h})]}\bigg(\dfrac{j-1}{N+1}\bigg),F^{-1}_{{\mathrm U}[E(\boldsymbol{h})]}\bigg(\dfrac{j}{N+1}\bigg)\bigg],\quad j\in\mathbb{N}^{+}_{N+1},\quad\boldsymbol{h}\in H^{E}_{+},
\]
where we observe that the endpoints of each interval never coincide
($F_{{\mathrm U}[E(\boldsymbol{h})]}$ is continuous as well). This results
in the top-level choice
\[
\phi(\boldsymbol{h},u)=\sum^{N+1}_{j=1}c_{j}(\boldsymbol{h})\mathds{1}_{E(\boldsymbol{h})\cap L_{j}(\boldsymbol{h})}(u),\quad(\boldsymbol{h},u)\in H^{E}_{+}\times[0,1],
\]
for which we can explicitly calculate, for each fixed $\boldsymbol{h}\in H^{E}_{+}$,
\begin{align*}
\lambda(E(\boldsymbol{h})\cap L_{j}(\boldsymbol{h})) & =\lambda(E(\boldsymbol{h}))\int^{1}_{0}\mathds{1}_{L_{j}(\boldsymbol{h})}(u)\dfrac{\mathds{1}_{E(\boldsymbol{h})}(u)}{\lambda(E(\boldsymbol{h}))}\mathrm{d}u\\
 & =\lambda(E(\boldsymbol{h}))\int^{1}_{0}\mathds{1}_{L_{j}(\boldsymbol{h})}(u)\mathrm{d}F_{{\mathrm U}[E(\boldsymbol{h})]}(u)\\
 & =\lambda(E(\boldsymbol{h}))\bigg[F_{{\mathrm U}[E(\boldsymbol{h})]}\bigg(F^{-1}_{{\mathrm U}[E(\boldsymbol{h})]}\bigg(\dfrac{j}{N+1}\bigg)\hspace{-1bp}\hspace{-1bp}\hspace{-1bp}\bigg)-F_{{\mathrm U}[E(\boldsymbol{h})]}\bigg(F^{-1}_{{\mathrm U}[E(\boldsymbol{h})]}\bigg(\dfrac{j-1}{N+1}\bigg)\hspace{-1bp}\hspace{-1bp}\hspace{-1bp}\bigg)\bigg]\\
 & =\dfrac{\lambda(E(\boldsymbol{h}))}{N+1}>0, \quad j\in \mathbb{N}^+_{N+1}.
\end{align*}
Therefore, we end up with the system of $N$ equations and $N+1$
unknowns (i.e., the $c_{j}$'s) reading
\[
\sum^{N+1}_{j=1}c_{j}(\boldsymbol{h})\bigg[\int^{1}_{0}\mathds{1}_{E(\boldsymbol{h})\cap L_{j}(\boldsymbol{h})}(u)z_{i}(\boldsymbol{h},u)\mathrm{d}u\bigg]=0,\quad i\in\mathbb{N}^{+}_{N},\quad\boldsymbol{h}\in H^{E}_{+}.
\]
This shows that, if at least one $c_{j}$'s is nonzero (for each $\boldsymbol{h}\in H^{E}_{+}$),
then $\phi\neq0$ with probability at least $\lambda(E(\boldsymbol{h}))/(N+1)$
on $H^{E}_{+}$, or total probability at least that of $E$, divided
by $N+1$ (strictly positive since $E$ is not null). Of course, we
may organize the above system into a (measurable) matrix-vector form
as
\[
\boldsymbol{A}(\boldsymbol{h})\boldsymbol{c}(\boldsymbol{h})={\bf 0},\quad\boldsymbol{h}\in H^{E}_{+}.
\]
By construction, such a system always has a nontrivial solution, say
$\boldsymbol{c}^{o}(\boldsymbol{h})$, which we may normalize
to unit Euclidean norm. Note that because one can use a finite-step
algorithm in order to determine $\boldsymbol{c}^{o}(\boldsymbol{h})$
via elementary arithmetic operations (e.g., Householder QR with a
fixed deterministic pivoting and sign-convention rule), $\boldsymbol{c}^{o}(\boldsymbol{h})$
can always be selected measurably.

By replacing $c_{j}$ with $c^{o}_{j}$, we may finally determine
the measurable choice (for every non-negligible event $E$)
\[
\phi(\boldsymbol{h},u)=\sum^{N+1}_{j=1}c^{o}_{j}(\boldsymbol{h})\mathds{1}_{E(\boldsymbol{h})\cap L_{j}(\boldsymbol{h})}(u),\quad(\boldsymbol{h},u)\in H^{E}_{+}\times[0,1],
\]
while $\phi=0$ outside $H^{E}_{+}$. In any case, it follows that
$\phi=0$ outside $E$, $\phi\neq0$ with always positive joint probability,
while it easy to see that $\phi\in{\cal L}_{\infty}(\mathrm{P}\hspace{-1bp}\hspace{-1bp}\times\hspace{-1bp}\hspace{-1bp}\lambda,\mathbb{R})$,
while satisfying the annilation property (\ref{eq:Ortho}). Consequently,
we have shown that the set of Radon-Nikodym derivatives ${\cal M}$
is thin, and therefore the Kingman-Robertson theorem applies.\hfill{}$\blacksquare$

\subsection{Proof of Theorem  \ref{thm:Main2}} \label{subsec: proof of theorem 2}
Again, we split the proof into two subsections. The first establishes convexity of the image \textit{closure} $\mathrm{cl}(\mathcal{C})$, manifesting the lack of decomposability for $\Pi$, and the second establishes strong duality by leveraging convexity of $\mathrm{cl}(\mathcal{C})$ only. Throughout this section, we let Assumption \ref{assu:Assumption2} be in effect.

\subsubsection[Proof Core: Convexity of cl(C)]{Convexity of $\boldsymbol{\mathrm{cl}({\cal C})}$}

The procedure follows similar steps as in the proof to Theorem \ref{thm:Main} presented above, but with some necessary modifications in the arguments taking place right after the invocation of Lemma \ref{thm:Weak_Lyapunov_G} in Section \ref{subsec:Core}. Throughout this section, we assume that Assumption \ref{assu:Assumption2} is in effect (instead of Assumption \ref{assu:Assumption}).

Specifically, let us recall the construction of the policy $\boldsymbol{p}_\alpha$, chosen as $\boldsymbol{p}_{\alpha}(\boldsymbol{s})=\mathds{1}_{E^{\alpha}}(\boldsymbol{s})\boldsymbol{p}(\boldsymbol{s})+\mathds{1}_{{\cal S}\backslash E^{\alpha}}(\boldsymbol{s})\boldsymbol{p}'(\boldsymbol{s}),\boldsymbol{s}\in\mathcal{S}$, where both $\boldsymbol{p}$ and $\boldsymbol{p}'$ are arbitrarily chosen in the policy set $\Pi$ (and $\alpha\in[0,1]$). Note that, in this case (i.e., under Assumption \ref{assu:Assumption2}), $\Pi$ is not decomposable, and therefore each policy $\boldsymbol{p}_\alpha$ does not in general belong to $\Pi$, but to $\Pi_{\Vert\cdot\Vert}\supseteq \Pi$ ($\Pi$ is densely decomposable). 
Nevertheless, thanks to condition $6'.(a)$ of Assumption $\ref{assu:Assumption2}$, the 
vector measure $\boldsymbol{G}_{ \boldsymbol{p}_\alpha}$ is well-defined (on the larger set $\Pi_{\Vert\cdot\Vert}$). Consequently, the complete argument following the policy sequence construction in Section \ref{subsec:Core} (proof of Theorem \ref{thm:Main}) is still valid, implying that
\[
-\boldsymbol{\rho}(-[\boldsymbol{F}\boldsymbol{p}_{\alpha}]_{\boldsymbol{f}}(\boldsymbol{H}))=-\boldsymbol{\rho}(-\alpha[\boldsymbol{F}\boldsymbol{p}]_{\boldsymbol{f}}(\boldsymbol{H})-(1-\alpha)[\boldsymbol{F}\boldsymbol{p}']_{\boldsymbol{f}}(\boldsymbol{H})).
\]
We now crucially exploit the assumption of $\Pi$ being densely decomposable relative to a norm $\Vert\cdot\Vert$: For every $\delta>0$, there exists
a policy $\overline{\boldsymbol{p}}_{\alpha}^{\delta}\in\Pi$ such that
$\Vert\overline{\boldsymbol{p}}_{\alpha}^{\delta}-\boldsymbol{p}_{\alpha}\Vert\le\delta$.
Since $-\boldsymbol{\rho}(-[\boldsymbol{F}\boldsymbol{p}_{\alpha}]_{\boldsymbol{f}}(\boldsymbol{H}))$
is assumed to also be (pointwise) continuous relative to $\boldsymbol{p}$ in the topology generated by $\Vert\cdot\Vert$,
we can choose, for every $\varepsilon >0$, a $\delta(\varepsilon,\boldsymbol{p}_{\alpha})>0$
such that 
\begin{equation}\nonumber
\Vert\overline{\boldsymbol{p}}_{\alpha}^{\delta}-\boldsymbol{p}_{\alpha}\Vert\le\delta\implies
\Vert
-\boldsymbol{\rho}(-[\boldsymbol{F}\overline{\boldsymbol{p}}_{\alpha}^{\delta}]_{\boldsymbol{f}}(\boldsymbol{H}))
-
[-\boldsymbol{\rho}(-[\boldsymbol{F}\boldsymbol{p}_{\alpha}]_{\boldsymbol{f}}(\boldsymbol{H}))]
\Vert_{\infty}\le\varepsilon.
\end{equation}
Therefore, we have that for every $\varepsilon>0$,
\begin{equation}\nonumber
\Vert
-\boldsymbol{\rho}(-[\boldsymbol{F}\overline{\boldsymbol{p}}_{\alpha}^{\delta}]_{\boldsymbol{f}}(\boldsymbol{H}))
-
[-\boldsymbol{\rho}(-\alpha[\boldsymbol{F}\boldsymbol{p}]_{\boldsymbol{f}}(\boldsymbol{H})-(1-\alpha)[\boldsymbol{F}\boldsymbol{p}']_{\boldsymbol{f}}(\boldsymbol{H}))]
\Vert_{\infty}\le\varepsilon,
\end{equation}
where $\boldsymbol{p}\in\Pi$, $\boldsymbol{p}'\in\Pi$ and also $\overline{\boldsymbol{p}}_{\alpha}^{\delta}\in\Pi$. Convexity further implies that\[
-\boldsymbol{\rho}(-[\boldsymbol{F}\overline{\boldsymbol{p}}_{\alpha}^{\delta}]_{\boldsymbol{f}}(\boldsymbol{H}))\ge
-\alpha\boldsymbol{\rho}(-[\boldsymbol{F}\boldsymbol{p}]_{\boldsymbol{f}}(\boldsymbol{H}))-(1-\alpha)\boldsymbol{\rho}(-[\boldsymbol{F}\boldsymbol{p}']_{\boldsymbol{f}}(\boldsymbol{H})) -\varepsilon \bf{1}.
\] 

Let us now see how the preceding fact implies that image set \textit{closure} $\mathrm{cl}({\cal C})$
is a convex set. Let $(\delta_{o},\boldsymbol{\delta}_{r},\boldsymbol{\delta}_{d})\in{\cal C}$
for $(\boldsymbol{x},\boldsymbol{p})\in{\cal X}\times\Pi$ and let
$(\delta'_{o},\boldsymbol{\delta}'_{r},\boldsymbol{\delta}'_{d})\in{\cal C}$
for $(\boldsymbol{x}',\boldsymbol{p}')\in{\cal X}\times\Pi$. Also
choose a sequence $\{\varepsilon_{n}>0\}_{n\in\mathbb{N}}$ decreasing
to zero, say $\varepsilon_{n}\triangleq1/(n+1),n\in\mathbb{N}$. By
our discussion above, we have actually shown that, for every $\alpha\in[0,1]$,
it holds that, for every $n$,
\[
\alpha(\delta_{o},\boldsymbol{\delta}_{r}-\varepsilon_{n}{\bf 1},\boldsymbol{\delta}_{d})+(1-\alpha)(\delta'_{o},\boldsymbol{\delta}'_{r}-\varepsilon_{n}{\bf 1},\boldsymbol{\delta}'_{d})\in{\cal C}.
\]
To verify this claim, observe that, for every $n\in\mathbb{N}$, there
exists a policy $\boldsymbol{p}_{\alpha}^{\varepsilon_{n}}\in\Pi$
such that for the pair $(\boldsymbol{x}_{\alpha}=\alpha\boldsymbol{x}+(1-\alpha)\boldsymbol{x}',\boldsymbol{p}_{\alpha}^{\varepsilon_{n}})\in{\cal X}\times\Pi$
(by assumption, ${\cal X}$ is convex) it is true that
\begin{align*}
 -\boldsymbol{\rho}(-[\boldsymbol{F}\boldsymbol{p}_{\alpha}^{\varepsilon_{n}}]_{\boldsymbol{f}}(\boldsymbol{H}))-\boldsymbol{x}_{\alpha}
 &  =-\boldsymbol{\rho}(-[\boldsymbol{F}\boldsymbol{p}_{\alpha}^{\varepsilon_{n}}]_{\boldsymbol{f}}(\boldsymbol{H}))-\alpha\boldsymbol{x}-(1-\alpha)\boldsymbol{x}'\\
 & \ge-\alpha\boldsymbol{\rho}(-[\boldsymbol{F}\boldsymbol{p}]_{\boldsymbol{f}}(\boldsymbol{H}))-(1-\alpha)\boldsymbol{\rho}(-[\boldsymbol{F}\boldsymbol{p}']_{\boldsymbol{f}}(\boldsymbol{H}))-\alpha\boldsymbol{x}-(1-\alpha)\boldsymbol{x}'-\varepsilon_{n}{\bf 1}\\
 & =-\alpha\boldsymbol{\rho}(-[\boldsymbol{F}\boldsymbol{p}]_{\boldsymbol{f}}(\boldsymbol{H}))-\alpha\boldsymbol{x}-(1-\alpha)\boldsymbol{\rho}(-[\boldsymbol{F}\boldsymbol{p}']_{\boldsymbol{f}}(\boldsymbol{H}))-(1-\alpha)\boldsymbol{x}'-\varepsilon_{n}{\bf 1}
 \\ &
 \ge\alpha(\boldsymbol{\delta}_{r}-\varepsilon_{n}{\bf 1})+(1-\alpha)(\boldsymbol{\delta}'_{r}-\varepsilon_{n}{\bf 1}).
\end{align*}
This verifies our claim. But, evidently, 
\[
\lim_{n\rightarrow\infty}\alpha(\delta_{o},\boldsymbol{\delta}_{r}-\varepsilon_{n}{\bf 1},\boldsymbol{\delta}_{d})+(1-\alpha)(\delta'_{o},\boldsymbol{\delta}'_{r}-\varepsilon_{n}{\bf 1},\boldsymbol{\delta}'_{d})=\alpha(\delta_{o},\boldsymbol{\delta}_{r},\boldsymbol{\delta}_{d})+(1-\alpha)(\delta'_{o},\boldsymbol{\delta}'_{r},\boldsymbol{\delta}'_{d}),
\]
which of course implies that
$
\alpha(\delta_{o},\boldsymbol{\delta}_{r},\boldsymbol{\delta}_{d})+(1-\alpha)(\delta'_{o},\boldsymbol{\delta}'_{r},\boldsymbol{\delta}'_{d})\in\mathrm{cl}({\cal C}),
$
as a limit point of elements in ${\cal C}$. Finally, convexity of
$\mathrm{cl}({\cal C})$ follows in light of the next elementary result,
whose proof we also present for completeness.
\vspace{6pt}
\begin{proposition}
\label{prop:PROP}Suppose that a set ${\cal A}\subseteq\mathbb{R}^{N}$
has the property that all convex combinations of any two points in
${\cal A}$ belong to its closure, i.e., $\boldsymbol{x}\in{\cal A}$ and $\boldsymbol{y}\in{\cal A}$ implies that, for every $\alpha \in [0,1]$, $\alpha\boldsymbol{x}+(1-\alpha)\boldsymbol{y}\in\mathrm{cl}({\cal A})$.
Then $\mathrm{cl}({\cal A})$ is (closed) convex.
\end{proposition}
\begin{proof}[Proof of Proposition \ref{prop:PROP}]
Let $\boldsymbol{x}\in\mathrm{cl}({\cal A})$ and $\boldsymbol{y}\in\mathrm{cl}({\cal A})$.
Then, we can find sequences $\{\boldsymbol{x}_{n}\in{\cal A}\}_{n\in\mathbb{N}}$
and $\{\boldsymbol{y}_{n}\in{\cal A}\}_{n\in\mathbb{N}}$ such that
$\boldsymbol{x}=\lim_{n\rightarrow\infty}\boldsymbol{x}_{n} \text{ and } \boldsymbol{y}=\lim_{n\rightarrow\infty}\boldsymbol{y}_{n},$
where the limits are interpreted in the standard Euclidean sense.
So, given $\alpha\in[0,1]$,
\begin{align*}
\alpha\boldsymbol{x}+(1-\alpha)\boldsymbol{y} & =\alpha\lim_{n\rightarrow\infty}\boldsymbol{x}_{n}+(1-\alpha)\lim_{n\rightarrow\infty}\boldsymbol{y}_{n}
=\lim_{n\rightarrow\infty}\alpha\boldsymbol{x}_{n}+(1-\alpha)\boldsymbol{y}_{n}.
\end{align*}
However, by assumption for the set ${\cal A}$ it holds that, for
every $n\in\mathbb{N}$,
$
\alpha\boldsymbol{x}_{n}+(1-\alpha)\boldsymbol{y}_{n}\in\mathrm{cl}({\cal A}),
$
and the sequence converges; therefore it must converge in $\mathrm{cl}({\cal A})$,
and the limit is $\alpha\boldsymbol{x}+(1-\alpha)\boldsymbol{y}$.
Therefore, $\mathrm{cl}({\cal A})$ is closed convex.
\end{proof}
\noindent By an application of Proposition \ref{prop:PROP},
we obtain that $\mathrm{cl}({\cal C})$ is a (closed) convex set.\hfill{}$\blacksquare$

\subsubsection[Convexity of cl(C)
Implies Strong Duality]{Convexity of $\boldsymbol{\mathrm{cl}({\cal C})}$
Implies Strong Duality} \label{subsec:Convexity-of-C}

Let us now finish the proof of Theorem \ref{thm:Main2} by exploiting
the convexity of the closure of set $\mathcal{C}$,
together
with condition $7$ of Assumption \ref{assu:Assumption2},
namely that problem \eqref{eq:Base} satisfies Slater's condition.
Let us also recall the Lagrangian associated with problem \eqref{eq:Base},
for which we already know that $\mathsf{D}^{*}\ge\mathsf{P}^{*}$.
Strong duality of \eqref{eq:Base} will be immediate if we can show
that $\mathsf{D}^{*}\le\mathsf{P}^{*}$ (and dual solutions exist).

Our discussion extends the usual standard separation argument, but is still eventually based on an application of the supporting hyperplane theorem (see e.g.,
\citep[Proposition 1.5.1]{Bertsekas2009}), which we outline below
for completeness. Note that the same technique would be applicable
in case our initial problem \eqref{eq:Base} was originally convex.
The reason is that the only fact needed at this point is the ``top-level''
convexity of $\mathrm{cl}({\cal C})$. 
Before we assemble
everything, we need an additional result.
\vspace{5pt}
\begin{lemma}[\textbf{Point on the Shell}]
\label{lem:Shell} $(\mathsf{P}^{*},{\bf 0},{\bf 0})$ is not in
the interior of $\mathrm{cl}({\cal C})$.
\end{lemma}
\vspace{-2pt}
\begin{proof}[Proof of Lemma \ref{lem:Shell}]
 First, we observe that the point $(\mathsf{P}^{*},{\bf 0},{\bf 0})$
cannot be in the interior of ${\cal C}$, for otherwise there would
exist an $\varepsilon>0$ such that $(\mathsf{P}^{*}+\varepsilon,{\bf 0},{\bf 0})\in{\cal C}$,
contradicting $\mathsf{P}^{*}$ being the optimal value of the initial
constrained problem. In fact, if $(\delta_{o}^{+},\boldsymbol{\delta}_{r}^{+},\boldsymbol{\delta}_{d}^{+})>{\bf 0}$,
every perturbation of the form
$
(\mathsf{P}^{*},{\bf 0},{\bf 0})+(\delta_{o}^{+},\boldsymbol{\delta}_{r}^{+},\boldsymbol{\delta}_{d}^{+})=(\mathsf{P}^{*}+\delta_{o}^{+},\boldsymbol{\delta}_{r}^{+},\boldsymbol{\delta}_{d}^{+})
$
cannot be in ${\cal C}$ either; if it was, this would imply the existence
of a pair $(\boldsymbol{x},\boldsymbol{p})\in{\cal X}\times\Pi$ attaining
a strictly larger objective than $\mathsf{P}^{*}$, while also keeping
all the constraints inactive. Therefore, it follows that the \textit{open
orthant (an open convex set)}
\[
{\cal D}=\{(\boldsymbol{x},\boldsymbol{y},\boldsymbol{z})|(\boldsymbol{x},\boldsymbol{y},\boldsymbol{z})=(\mathsf{P}^{*}+\delta_{o}^{+},\boldsymbol{\delta}_{r}^{+},\boldsymbol{\delta}_{d}^{+}),\quad(\delta_{o}^{+},\boldsymbol{\delta}_{r}^{+},\boldsymbol{\delta}_{d}^{+})>{\bf 0}\}
\]
is disjoint from ${\cal C}$. 

Let us show that, in fact, ${\cal D}$ and $\mathrm{cl}({\cal C})$
are also disjoint. Note that ${\cal D}$ has a nonempty interior in
the standard Euclidean topology of ${\cal C}$. Therefore, for every
element of ${\cal D}$ there exists an open ball of the same dimension
as that of ${\cal C}$ entirely contained in ${\cal D}$ and having
that element of ${\cal D}$ as its center. Next, if an element of
${\cal D}$, say $\boldsymbol{d}$, was in $\mathrm{cl}({\cal C})$,
there should exist a sequence contained entirely in ${\cal C}$ converging
to $\boldsymbol{d}$. However, since $\boldsymbol{d}$ is contained
in an open ball, say ${\cal B}_{\boldsymbol{d}}$, entirely contained
in ${\cal D}$ and of the same dimension as that of ${\cal C}$, every
sequence in ${\cal C}$ converging to $\boldsymbol{d}$ must break
into that open ball. This implies that the aforementioned sequence
in ${\cal C}$ converging to $\boldsymbol{d}$ must have elements
also in ${\cal B}_{\boldsymbol{d}}$. This is absurd, since ${\cal C}\bigcap{\cal B}_{\boldsymbol{d}}=\emptyset$
(due to the fact that ${\cal C}\bigcap{\cal D}=\emptyset$).

Now suppose that $(\mathsf{P}^{*},{\bf 0},{\bf 0})$ is in the interior
of $\mathrm{cl}({\cal C})$. Then, there must exist another open ball,
say ${\cal B}_{*}$, of the same dimension as that of ${\cal C}$,
completely contained in $\mathrm{cl}({\cal C})$ and centered at $(\mathsf{P}^{*},{\bf 0},{\bf 0})$.
This means that there exists a positive perturbation vector $(\delta_{o}^{+},\boldsymbol{\delta}_{r}^{+},\boldsymbol{\delta}_{d}^{+})>{\bf 0}$
such that the point $(\mathsf{P}^{*}+\delta_{o}^{+},\boldsymbol{\delta}_{r}^{+},\boldsymbol{\delta}_{d}^{+})$
is in ${\cal B}_{*}$. This is also absurd, because $(\mathsf{P}^{*}+\delta_{o}^{+},\boldsymbol{\delta}_{r}^{+},\boldsymbol{\delta}_{d}^{+})\in{\cal D}$
and ${\cal D}\bigcap\mathrm{cl}({\cal C})=\emptyset$, as shown above.
\end{proof}
\noindent The complete argument based on an application of the supporting hyperplane
theorem now follows. 
\vspace{6pt}
\begin{theorem}[\textbf{Supporting Hyperplane Theorem}]
 Let ${\cal A}\subseteq\mathbb{R}^{n}$ be a nonempty convex set.
If $\boldsymbol{\delta}^{*}\in\mathbb{R}^{n}$ is not in the interior
of ${\cal A}$, there exists a hyperplane passing through $\boldsymbol{\delta}^{*}$
such that ${\cal A}$ is in one of its closed halfspaces. That is, there exists  $\text{\ensuremath{\boldsymbol{\lambda}}}\neq{\bf 0}$
such that, for every $\boldsymbol{\delta}\in{\cal A}$, $\langle\text{\ensuremath{\boldsymbol{\lambda}}},\boldsymbol{\delta}^{*}\rangle\ge\langle\text{\ensuremath{\boldsymbol{\lambda}}},\boldsymbol{\delta}\rangle$.
\end{theorem}
\vspace{6pt}
We start by observing that since \eqref{eq:Base} satisfies Slater's
condition, it follows that ${\cal C}$ is nonempty, and the same follows
for $\mathrm{cl}({\cal C})$. So $\mathrm{cl}({\cal C})$ is a nonempty
(closed) convex set. By Lemma \ref{lem:Shell}, we also know that
the point $(\mathsf{P}^{*},{\bf 0},{\bf 0})$ is not in the interior
of $\mathrm{cl}({\cal C})$. We can then apply the supporting hyperplane
theorem on the pair $\mathrm{cl}({\cal C})$ and $(\mathsf{P}^{*},{\bf 0},{\bf 0})$,
implying existence of a vector of multipliers $(\lambda_{o},\boldsymbol{\lambda}_{\boldsymbol{\rho}},\boldsymbol{\lambda}_{\boldsymbol{g}})\neq{\bf 0}$
such that, for every $(\delta_{o},\boldsymbol{\delta}_{r},\boldsymbol{\delta}_{d})\in\mathrm{cl}({\cal C})$,
it is true that
$
\lambda_{o}\delta_{o}+\langle\boldsymbol{\lambda}_{\boldsymbol{\rho}},\boldsymbol{\delta}_{r}\rangle +\langle\boldsymbol{\lambda}_{\boldsymbol{g}},\boldsymbol{\delta}_{d}\rangle \le\lambda_{o}\mathsf{P}^{*}.
$

It readily follows that, in fact, $(\lambda_{o},\boldsymbol{\lambda}_{\boldsymbol{\rho}},\boldsymbol{\lambda}_{\boldsymbol{g}})\ge{\bf 0}$.
Indeed, if any component of $(\lambda_{o},\boldsymbol{\lambda}_{\boldsymbol{\rho}},\boldsymbol{\lambda}_{\boldsymbol{g}})$
was negative, then we could choose $(\delta_{o},\boldsymbol{\delta}_{r},\boldsymbol{\delta}_{d})\in{\cal C}\subseteq\mathrm{cl}({\cal C})$
such that the corresponding inner product becomes arbitrarily large
(note that ${\cal C}$ is unbounded below), eventually violating the
inequality above, regardless of the sign of $\mathsf{P}^{*}$.

The second fact we may show is that $\lambda_{o}\neq0$, which implies
that $\lambda_{o}>0$. Again, if $\lambda_{o}=0$, then we would have
that
$\langle\boldsymbol{\lambda}_{\boldsymbol{\rho}},\boldsymbol{\delta}_{r}\rangle + \langle\boldsymbol{\lambda}_{\boldsymbol{g}},\boldsymbol{\delta}_{d}\rangle\le0$.
But $(\boldsymbol{\lambda}_{\boldsymbol{\rho}},\boldsymbol{\lambda}_{\boldsymbol{g}})\neq{\bf 0}$
(i.e., there is at least one nonzero component) and $(\boldsymbol{\lambda}_{\boldsymbol{\rho}},\boldsymbol{\lambda}_{\boldsymbol{g}})\ge{\bf 0}$
(shown above), and problem \eqref{eq:Base} satisfies Slater's condition,
so the preceding inequality is also absurd.

As a result, we may divide by $\lambda_{o}>0$, showing that there
is $(\boldsymbol{\lambda}_{\boldsymbol{\rho}}^{\star}\triangleq\boldsymbol{\lambda}_{\boldsymbol{\rho}}/\lambda_{o},\boldsymbol{\lambda}_{\boldsymbol{g}}^{\star}\triangleq\boldsymbol{\lambda}_{\boldsymbol{g}}/\lambda_{o})\ge{\bf 0}$,
such that, for every $(\delta_{o},\boldsymbol{\delta}_{r},\boldsymbol{\delta}_{d})\in\mathrm{cl}({\cal C})$,
it holds that $\delta_{o}+\langle\boldsymbol{\lambda}_{\boldsymbol{\rho}}^{\star},\boldsymbol{\delta}_{r}\rangle +\langle\boldsymbol{\lambda}_{\boldsymbol{g}}^{\star},\boldsymbol{\delta}_{d}\rangle\le\mathsf{P}^{*}$.
By construction of the set $\mathrm{cl}({\cal C})$ (as the closure
of ${\cal C}$), we obtain in particular that
\[
\mathsf{L}(\boldsymbol{x},\boldsymbol{p},\boldsymbol{\lambda}^{\star})=g^{o}(\boldsymbol{x})+\langle\boldsymbol{\lambda}_{\boldsymbol{g}}^{\star},\boldsymbol{g}(\boldsymbol{x})\rangle+\big\langle\boldsymbol{\lambda}_{\boldsymbol{\rho}}^{\star},-\boldsymbol{\rho}(-[\boldsymbol{F}\boldsymbol{p}]_{\boldsymbol{f}}(\boldsymbol{H}))-\boldsymbol{x}\big\rangle\le\mathsf{P}^{*},
\]
for every pair $(\boldsymbol{x},\boldsymbol{p})\in{\cal X}\times\Pi$,
since all such pairs correspond to points included in ${\cal C}\subseteq\mathrm{cl}({\cal C})$.
This further implies that we are allowed to maximize both sides over all possible $(\boldsymbol{x},\boldsymbol{p})$,
yielding
\begin{align*}
-\infty<\mathsf{D}^{*} & =\inf_{\boldsymbol{\lambda}\ge{\bf 0}}\sup_{(\boldsymbol{x},\boldsymbol{p})\in{\cal X}\times\Pi}\mathsf{L}(\boldsymbol{x},\boldsymbol{p},\boldsymbol{\lambda}) \le\sup_{(\boldsymbol{x},\boldsymbol{p})\in{\cal X}\times\Pi}\mathsf{L}(\boldsymbol{x},\boldsymbol{p},\boldsymbol{\lambda}^{\star})\le\mathsf{P}^{*},
\end{align*}
and we are done.\hfill{}$\blacksquare$

\section{Reductions, Counterexamples and Extensions}\label{Counter}
\subsection{Reduction to the Risk-Neutral Case} \label{subsec: risk-neutral reduction}
It is interesting to see how (\ref{eq:Base}) reduces to a classical
risk-neutral setting, verifying that our results generalize those
in \citep{Luo2008,Ribeiro2012,Chamon2020a,Chamon2021}, which establish strong duality for the constrained
\textit{integral-functional} problem 
\begin{equation}
\begin{array}{rl}
\underset{\boldsymbol{x},\boldsymbol{p}(\cdot)}{\mathrm{maximize}} & g^{o}(\boldsymbol{x})\\
\mathrm{subject\,to} & \boldsymbol{x}\le\mathbb{E}\{\boldsymbol{f}(\boldsymbol{p}(\boldsymbol{U}),\boldsymbol{U})\}\\
 & \boldsymbol{g}(\boldsymbol{x})\ge{\bf 0}\\
 & (\boldsymbol{x},\boldsymbol{p})\in{\cal X}\times\Pi
\end{array},\tag{{RNP}}\label{eq:RNP}
\end{equation}
posed under the assumptions of concavity of $g^{o}$, $\boldsymbol{g}$,
mere integrability of $\boldsymbol{f}$, decomposability of $\Pi$,
and nonatomicity of the distribution of $\boldsymbol{U}$. Mapping
to our setting herein, let us for simplicity set $\boldsymbol{\rho}=\mathbb{E}$
and choose $\mathrm{P}=\delta_{{\bf 0}}$. Then $[\boldsymbol{F}\boldsymbol{p}]_{\boldsymbol{f}}(\boldsymbol{H})=[\boldsymbol{F}\boldsymbol{p}]_{\boldsymbol{f}}({\bf 0})$
almost everywhere relative to $\mathrm{P}$ (i.e., the smooth operator
is essentially constant), and (\ref{eq:Base}) reduces to 
\[
\begin{array}{rl}
\underset{\boldsymbol{x},\boldsymbol{p}({\bf 0},\cdot)}{\mathrm{maximize}} & g^{o}(\boldsymbol{x})\\
\mathrm{subject\,to} & \boldsymbol{x}\le\mathbb{E}\{\boldsymbol{f}(\boldsymbol{p}({\bf 0},\boldsymbol{U}),{\bf 0},\boldsymbol{U})\}\\
 & \boldsymbol{g}(\boldsymbol{x})\ge{\bf 0}\\
 & (\boldsymbol{x},\boldsymbol{p}({\bf 0},\boldsymbol{U}))\in{\cal X}\times\Pi_{\{{\bf 0}\}\hspace{-1bp}\times{\cal U}}
\end{array},
\]
where $\Pi_{\{{\bf 0}\}\hspace{-1bp}\times{\cal U}}$ denotes the
restriction of $\Pi$ to the vertical section $\{{\bf 0}\}\hspace{-1bp}\times{\cal U}$,
which can be easily verified to be decomposable as long as $\Pi$
is. Therefore, this reduction is identical to problem (\ref{eq:RNP})
(after renaming variables), and under Assumption \ref{assu:Assumption}, in particular
nonatomicity of $\mathrm{P}_{\boldsymbol{U}|\boldsymbol{H}}=\mathrm{P}_{\boldsymbol{U}|{\bf 0}}$
almost everywhere relative to $\mathrm{P}$, we get strong duality
via Theorem \ref{thm:Main}. 

In passing, note that $\boldsymbol{\rho}$ can in fact be chosen as
any risk measure which coincides with the identity operator on equivalence
classes of constants (i.e., $\boldsymbol{\rho}(\boldsymbol{Z})=\boldsymbol{c}$
for $\boldsymbol{Z}=\boldsymbol{c}$ with probability $1$).

\subsection{Smooth Operator Strikes Back: Counterexamples}
Inspired, in a sense, from the risk-neutral setting (see above), one
might wonder at this point whether strong duality holds for the \textit{``pure
risk'' }counterpart of\textit{ }(\ref{eq:Base})
\begin{equation}
\begin{array}{rl}
\underset{\boldsymbol{x},\boldsymbol{p}(\cdot)}{\mathrm{maximize}} & g^{o}(\boldsymbol{x})\\
\mathrm{subject\,to} & \boldsymbol{x}\le-\boldsymbol{\rho}(-\boldsymbol{f}(\boldsymbol{p}(\boldsymbol{S}),\boldsymbol{S}))\\
 & \boldsymbol{g}(\boldsymbol{x})\ge{\bf 0}\\
 & (\boldsymbol{x},\boldsymbol{p})\in{\cal X}\times\Pi
\end{array},\tag{{RCP\ensuremath{\backslash\mathbb{E}}}}\label{eq:RNP-NONE}
\end{equation}
under nonatomicity of $\mathrm{P}_{\boldsymbol{S}}$ and the rest
of Assumption \ref{assu:Assumption} (possibly relaxing condition $3$). In other words,
one might ask whether the possibly somewhat nonstandard presence of
the smooth operator (conditional expectation) $[\boldsymbol{F}\boldsymbol{p}]_{\boldsymbol{f}}$
is essential in establishing strong duality (under our assumptions).
Going a bit further, even if this possibly failed, one might also
wonder whether the \textit{``pure randomization'' }counterpart of
the pure risk problem (\ref{eq:RNP-NONE}) 
\begin{equation}
\begin{array}{rl}
\underset{\boldsymbol{x},\boldsymbol{p}(\cdot)}{\mathrm{maximize}} & g^{o}(\boldsymbol{x})\\
\mathrm{subject\,to} & \boldsymbol{x}\le-\boldsymbol{\rho}(-\boldsymbol{f}(\boldsymbol{p}(\boldsymbol{H},\boldsymbol{U}),\boldsymbol{H}))\\
 & \boldsymbol{g}(\boldsymbol{x})\ge{\bf 0}\\
 & (\boldsymbol{x},\boldsymbol{p})\in{\cal X}\times\Pi
\end{array},\tag{{RCP\ensuremath{\backslash\mathbb{E}}-PureR}}\label{eq:RNP-NONE_Pure}
\end{equation}
exhibits strong duality under the same conditions, and where $U$
is standard uniform and independent of $\boldsymbol{H}$. This is
a generalization of the more specific question of whether pure policy
randomization provides convexification, in the risk-constrained nonconvex
setting above (which is generally true \textit{even at problem level}
in the risk-neutral setting, as long as the decomposable set $\Pi$
translates to a convex set of stochastic kernels; see also our ``warmup'' in 
Section \ref{warmup}). Note that \textit{both statements are true
in the risk-neutral setting}, i.e., whenever $\boldsymbol{\rho}=\mathbb{E}$
\citep{Luo2008,Ribeiro2012,Chamon2020a,Chamon2021}.

In the following, we provide a single explicit counterexample satisfying
Assumption \ref{assu:Assumption} (in fact, quite more) which constitutes an instance
of \textit{both} problems (\ref{eq:RNP-NONE}) and (\ref{eq:RNP-NONE_Pure}),
and which exhibits a non-zero duality gap, explicitly characterized.
Specifically, for $(\alpha,\beta)\in(0,1)^{2}$, we consider the nonconvex
risk-constrained problem
\begin{equation}
\begin{array}{rl}
\underset{p(\cdot)}{\mathrm{minimize}} & \mathbb{E}\{f_{0}(p(h,u))\}\\
\mathrm{subject\,to} & \mathrm{CVaR}_{\alpha}(f_{1}(p(h,u)))\le\beta\\
 & p(\cdot)\in[0,1]
\end{array},\tag{{CE}}\label{eq:CE}
\end{equation}
where $\mathrm{P}_{\boldsymbol{S}}$ is taken as the product of Lebesgue
measure on $[0,1]$ with itself, and $\boldsymbol{f}=[f_{0},f_{1}]$
with $f_{0(1)}(p) =-(+)p+\kappa p(1-p),p\in[0,1],\kappa\ge2.$
Note that $\boldsymbol{f}$ is polynomial in each dimension, thus
infinitely differentiable, and also essentially bounded relative to
$\mathrm{P}_{\boldsymbol{S}}$. Also, the set of measurable policies
taking values in the action set $[0,1]$ is trivially decomposable.
Observe that problem (\ref{eq:CE}) can be rewritten
as an instance of either (\ref{eq:RNP-NONE}) or (\ref{eq:RNP-NONE_Pure})
reading ($\boldsymbol{f}$ is state-independent)
\begin{equation}
\begin{array}{rl}
\underset{t,p(\cdot)}{\mathrm{minimize}} & t\\
\mathrm{subject\,to} & \begin{bmatrix}\hspace{-1bp}\begin{array}{l}
\mathbb{E}\{f_{0}(p(h,u))\}\\
\mathrm{CVaR}_{\alpha}(f_{1}(p(h,u))-\beta)
\end{array}\hspace{-1bp}\end{bmatrix}\le\begin{bmatrix}t\\
t'
\end{bmatrix}\\
 & (t,t',p(\cdot))\in\mathbb{R}\times\{0\}\times[0,1]
\end{array},\tag{{CE\ensuremath{'}}}
\end{equation}
whose primal \textit{and} dual values coincide with those of (\ref{eq:CE}).
It is plain to see that (\ref{eq:CE}) satisfies Asssumption \ref{assu:Assumption}; in
fact, (\ref{eq:CE}) resides in a strictly smaller class of problem
instances since, e.g., $\mathrm{P}_{\boldsymbol{S}}$ is nonatomic
in both its components (corresponding to $h$ and $u$, respectively).
We have the following result, characterizing the duality gap of problem
(\ref{eq:CE}). 
\smallskip
\begin{proposition}
\label{prop:Counter}Let $\kappa\ge2$ and choose $(\alpha,\beta)\in(0,1)^{2}$.
Then, the primal and dual values of problem (\ref{eq:CE}) are $-\beta\alpha$
and $-\beta$, respectively. In particular, the duality gap of (\ref{eq:CE})
equals $(1-\alpha)\beta>0$.
\end{proposition}
\smallskip
Proposition \ref{prop:Counter} (whose proof is omitted herein for
brevity) shows that strong duality \textit{cannot} be exhibited \textit{uniformly}
within the classes of problem instances abstracted by problems (\ref{eq:RNP-NONE})
and (\ref{eq:RNP-NONE_Pure}) satisfying Assumption \ref{assu:Assumption} (and when
condition $3$ is replaced by mere nonatomicity of $\mathrm{P}_{\boldsymbol{S}}$, which is in general weaker than condition $3$). Consequently, establishing
strong duality under Assumption \ref{assu:Assumption} for both problem settings is a provably hopeless
pursuit, effectively putting both questions formulated above peacefully
to rest.

From a structural perspective, the counterexample constructed above
demonstrates that, in the nonconvex setting, nonatomicity and decomposability
are \textit{not} enough to provide strong duality, and that nonlinearity
of (convex) functionals is the culprit (manifesting for any nontrivial
risk measure with no expectation representation); in fact, as $\alpha\nearrow1$,
it is of course the case that $\mathrm{CVaR}_{\alpha}\searrow\mathbb{E}$,
and as expected the duality gap vanishes, and strong duality
is restored. On the other hand, as we establish via Theorem \ref{thm:Main} in particular, the
presence of \textit{some} linearity, along \textit{at least one} dimension
of random observables--- \textit{is} enough to restore strong duality,
for any ensemble of convex
risk measures applied to the uncertainty that remains.

\subsection{Model Extensions} \label{subsec:theory extension}
\par Let us briefly demonstrate how our main results in (Section \ref{sec: Main results}) can be applied to certain (more or less obvious) variations or generalizations of the base model in \eqref{eq:Base}. First, it is an exercise to show that strong duality of problem \eqref{eq:Base} implies strong duality for problems in which risk components also appear in the objective. To this end, let us consider another risk vector $-\boldsymbol{\rho}^o(-[\boldsymbol{F}\boldsymbol{p}]_{\boldsymbol{f}^o}(\boldsymbol{H}))$ under the usual assumptions (either Assumption \ref{assu:Assumption} or \ref{assu:Assumption2}), where ``$o$" stands for ``objective". Second, and perhaps more trivially, one can replace $\boldsymbol{x}$ on the left-hand side of the vector risk constraint of \eqref{eq:Base} with $\boldsymbol{w}(\boldsymbol{x})$, where $\boldsymbol{w}$ is a convex function. Indeed, the proofs would follow exactly the developments of Section \ref{sec:Proof-of-Theorems}, requiring only certain  minor modifications. 
Combining these observations, we see that both Theorems \ref{thm:Main} and \ref{thm:Main2} are applicable to general \textit{risk-over-risk problems} of the form
\vspace{-1.5pt}
\begin{tagequation}{ROR}
\boxed{
\begin{array}{rl}
\underset{\boldsymbol{x},\boldsymbol{p}(\cdot)}{\mathrm{maximize}} & g^{o}(\boldsymbol{x}) +
\langle 
\mathbf{1}, -\boldsymbol{\rho}^o(-[\boldsymbol{F}\boldsymbol{p}]_{\boldsymbol{f}^o}(\boldsymbol{H})) 
\rangle
\\
\mathrm{subject\,to} & \boldsymbol{w}(\boldsymbol{x})\le
-\boldsymbol{\rho}(-[\boldsymbol{F}\boldsymbol{p}]_{\boldsymbol{f}}(\boldsymbol{H}))
\\
 & \boldsymbol{g}(\boldsymbol{x})\ge{\bf 0}\\
 & (\boldsymbol{x},\boldsymbol{p})\in{\cal X}\times\Pi
\end{array},
}
\end{tagequation}

\noindent under the appropriate conditions on the involved risk measures and service functions, as those are implied by the discussion above.

\par Lastly, we conjecture that Theorems \ref{thm:Main} and \ref{thm:Main2} can be established for more general conic formulations of \eqref{eq:Base} (that is, replacing inequality constraints by convex conic inequalities), by properly extending our proof, for instance utilizing results developed in \citep{FloresBazan2013}.
This more technical treatment is omitted here, for simplicity of exposition.

\section{Applications} \label{sec:Applications}
\par The base problem \eqref{eq:Base} is general, subsuming as special cases various useful problems arising in several applications. In this section, we discuss a small subset of such applications in more detail to highlight this fact. More specifically, we demonstrate how (risk-constrained versions of) certain problems arising in \textit{wireless systems resource allocation}, and \textit{nonconvex risk-constrained learning} can be cast into the form of \eqref{eq:Base}. In turn, this allows us to show that wide classes of such problems exhibit strong duality, substantially improving previously known results.
\par Before we proceed and as already hinted earlier, it is worth noting that in resource allocation applications it is natural to use the setting of Theorem \ref{thm:Main} (and the associated Assumption \ref{assu:Assumption}); indeed, (static) control problems arising in this context benefit from the utilization of a rich policy space ${\Pi}$, and thus the assumption that $\Pi$ is decomposable is fairly natural or even desirable in this case. On the other hand, when considering constrained learning problems, the situation is different. In this case, it may be desirable to restrict the space of policies $\Pi$. Indeed, allowing for a large $\Pi$ can result in solutions which, in practice, may massively overfit the data (in the finite-sample/empirical setting) or are hard to interpret (especially in high dimensions), and generally might not be very practical and/or computationally valuable. Thus, in this case, it may be more sensible to make use of Theorem \ref{thm:Main2} (and its associated Assumption \ref{assu:Assumption2}) instead. In general, whether one should opt for the setting of Theorem \ref{thm:Main} or that of Theorem \ref{thm:Main2} ultimately depends on the application at hand and the regularity of the resulting optimization problem. For instance, the use of continuous or smooth policies in (static) control settings might also be desirable for better interpretability of the resulting controller, which would make Theorem \ref{thm:Main2} a more suitable choice in such a case.

\subsection{Risk-Constrained Resource Allocation in Wireless Communications} \label{subsec: resource allocation examples}
\par As discussed in Section \ref{subsec:Contributions}, problem \eqref{eq:Base} is motivated by resource allocation problems from wireless systems engineering. In this setting, the random vector $\boldsymbol{S} \in \mathcal{S}$ quantifies the instantaneous quality of links between communicating or networking entities (nodes), $\boldsymbol{p}$ represents a resource allocation policy, and $\boldsymbol{f}$ measures instantaneous service levels achieved by a policy $\boldsymbol{p}$ operating on the observables $\boldsymbol{S}$. In line with the two-level structure of \eqref{eq:Base}, the observables split as $\boldsymbol{S}=(\boldsymbol{H},\boldsymbol{U})$, where $\boldsymbol{H}\in\mathcal{H}$ is the component over which risk is assessed and $\boldsymbol{U}$ is the ``inner" variable averaged out by the smooth operator. The risk of $\boldsymbol{f}$ is associated with the \textit{risk-ergodic performance} vector $\boldsymbol{x} \in \mathbb{R}^N$, and utilities $g^{o}$ and $\boldsymbol{g}$ evaluate those ergodic risks. In fact, there are numerous resource allocation tasks arising in practical settings that are readily modeled via \eqref{eq:Base}, primarily considered in the risk-neutral setting (see, e.g., \citep{Luo2008,Neely2010,WeiYu2006,Ribeiro2010,Ribeiro2010a,Ribeiro2012,Eisen2019,Kalogerias2020e,Yaylali2023}). Below we consider two related examples to showcase the expressive power of \eqref{eq:Base} in the risk-constrained setting.

\subsubsection{Multiple Access Interference Channel}\label{MAISec}
We first develop a generalization of a model studied in, e.g., \citep[Section VI.B]{Eisen2019}, \citep[Section VI]{Kalogerias2020e} and \citep{Kalogerias2020b}, where a \emph{multiple access interference channel} model is considered, in which $N_Q\equiv N_{\boldsymbol{H}}$ wireless transmitters simultaneously communicate with a central node, e.g., a base station. Here, each component of $\boldsymbol{S}$ corresponds to the strength of the communication channel between each transmitter and the base station. The signal emitted by each transmitter introduces interference to all other signals emitted by the remaining transmitters, and the goal is to optimally allocate power to each transmitter after observing the instantaneous channel state, so as to maximize a certain Quality-of-Service (QoS) network-wide utility, under a total expected power specification $p_{max}>0$.
 
Towards leveraging Theorem \ref{thm:Main}, we  exploit a common two-timescale structure of wireless channels by splitting the observables as $\boldsymbol{S}=(\boldsymbol{H},\boldsymbol{U})$. Here, $\boldsymbol{H}\in\mathbb{R}^{N_Q}_{+}$ collects \emph{large-scale} (slow fading) gains, capturing distance-dependent path loss and shadowing (which vary slowly with node geometry and propagation environment), whereas $\boldsymbol{U}\in\mathbb{R}^{N_Q}_{+}$ collects \emph{small-scale} (fast/multipath fading) gains, whose conditional law $\mathrm{P}_{\boldsymbol{U}\mid\boldsymbol{H}}$ is naturally modeled as nonatomic $\mathrm{P}$-a.e. (e.g., Rayleigh or Rician). The instantaneous gain of each link factorizes as the \emph{product} of these two effects, so that the effective (per-link) gain vector is $\boldsymbol{H}\odot\boldsymbol{U}$; this multiplicative decomposition is standard in wireless channel modeling (see, e.g., \cite[Chapters 2 and 3]{Goldsmith2005Wireless} and \citep{Tse2005}), and reflects that power received over a link equals the transmit power scaled by (slowly varying) large-scale attenuation and (rapidly varying) small-scale fading. 
Transmission power is allocated on the basis of the full split state, i.e., $\boldsymbol{p}(\boldsymbol{H},\boldsymbol{U})\ge\mathbf{0}$, and the QoS of each transmitter is measured through its \emph{fast rate}, i.e., the fast-fading average \emph{within} each slow-fading state produced by the smooth operator $[\boldsymbol{F}\boldsymbol{p}]_{\boldsymbol{f}_Q}(\boldsymbol{H})=\mathbb{E}\{\boldsymbol{f}_Q(\boldsymbol{p}(\boldsymbol{S}),\boldsymbol{S})|\boldsymbol{H}\}$, whose \emph{risk} is then assessed across the residual slow-fading uncertainty via a risk measure $\boldsymbol{\rho}_Q$. This reflects the two-level information structure in \eqref{eq:Base}: fast fading is averaged out (as coding takes place over many small-scale realizations within a coherence block), while slow-fading, too slowly-changing to average out, is the quantity over which imposing risk is meaningful.
 
More specifically, consider the hollow matrix $\boldsymbol{T}\triangleq\mathbf{1}\mathbf{1}^\top - \boldsymbol{I}$ (where $\mathbf{1}$ denotes the vector of ones of appropriate dimension), and define the service vector stacking the (instantaneous) communication rates achievable by all transmitters at the base station, i.e.,
\[\boldsymbol{f}_{Q}(\boldsymbol{p}(\boldsymbol{S}),\boldsymbol{S}) \triangleq
\log\left(\mathbf{1}+ \frac{(\boldsymbol{H}\odot\boldsymbol{U})\odot \boldsymbol{p}(\boldsymbol{H},\boldsymbol{U})}{\sigma^2\mathbf{1} 
+
\boldsymbol{T} [(\boldsymbol{H}\odot\boldsymbol{U})\odot \boldsymbol{p}(\boldsymbol{H},\boldsymbol{U})]
}\right), \quad \boldsymbol{p} \in \Pi, \]
where the $\log (\cdot)$ and division operations on vectors are interpreted component-wise, $\sigma^2 > 0$ denotes the variance of the (common) reception noise induced by the transmission to the base station, and $\Pi$ is any $\mathscr{B}(\mathcal{S})$-decomposable space of nonnegative resource policies; note that $\boldsymbol{p}$ represents power. Observe that $\boldsymbol{f}_{Q}(\cdot,\boldsymbol{S})$ is component-wise nonconcave and nonlinearly coupled, i.e., the communication rate of each transmitter depends on the power allocations of all transmitters in the network. We also define the scalar resource coupling function ${{f}_C}(\boldsymbol{p}(\boldsymbol{S}),\boldsymbol{S}) \triangleq  p_{max} - \langle \boldsymbol{p}(\boldsymbol{H},\boldsymbol{U}),\mathbf{1} \rangle,\boldsymbol{p}\in\Pi$. By choosing for simplicity $g^{o}(\boldsymbol{x}) = \langle \boldsymbol{\alpha}, \boldsymbol{x} \rangle$, i.e., a linear \textit{sumrate utility} with nonnegative weights $\boldsymbol{\alpha} \geq \mathbf{0}$ evaluating the \textit{(risk-)ergodic} service vector $\boldsymbol{x}\in\mathbb{R}^{N_Q}$, the corresponding stochastic resource allocation problem studied in \citep{Eisen2019, Kalogerias2020e, Kalogerias2020b} admits a risk-constrained generalization which reads as
\begin{tagequation}{MAI}
\begin{array}{rl}
\underset{\boldsymbol{x},\boldsymbol{p}(\cdot)}{\mathrm{maximize}} & \langle \boldsymbol{\alpha},\boldsymbol{x}\rangle\\
\mathrm{subject\,to} & \boldsymbol{x} \le
\displaystyle{
-\boldsymbol{\rho}_{Q}\bigg[\hspace{-2bp}-\hspace{-1bp}
\mathbb{E}\bigg\{\log\left(\mathbf{1}+ \frac{(\boldsymbol{H}\odot\boldsymbol{U})\odot \boldsymbol{p}(\boldsymbol{H},\boldsymbol{U})}{\sigma^2\mathbf{1}
+
\boldsymbol{T} [(\boldsymbol{H}\odot\boldsymbol{U})\odot \boldsymbol{p}(\boldsymbol{H},\boldsymbol{U})]
}\right)\hspace{-2.5pt}\bigg|\boldsymbol{H}\bigg\}
\hspace{-2bp}\bigg]}
\\
\vspace{-8bp}
\\
& \langle \mathbb{E}\{ \boldsymbol{p}(\boldsymbol{H},\boldsymbol{U}) \} , \mathbf{1} \rangle
\le p_{max}
\\
 &  (\boldsymbol{x},\boldsymbol{p})\in \mathbb{R}_{+}^{N_Q} \times\Pi
\end{array}\label{MAI},
\end{tagequation}
\hspace{-3pt}where $\boldsymbol{\rho}_Q$ is a finite-valued vector risk measure on $\mathcal{L}_p(\mathrm{P},\mathbb{R}^{N_Q})$ and compatible with the construction of \eqref{eq:Base}; for instance, $\boldsymbol{\rho}_Q$ could be a vector of $\mathrm{CVaRs}$ with potentially distinct confidence levels.
\par Evidently, problem \eqref{MAI} is a special case of \eqref{eq:Base}, since the former can be rewritten as
\begin{equation*}
\begin{array}{rl}
\underset{(\boldsymbol{x}_Q,x_C),\boldsymbol{p}(\cdot)}{\mathrm{maximize}} & g^{o}(\boldsymbol{x}_Q)\\
\mathrm{subject\,to}\, & \begin{bmatrix}\boldsymbol{x}_Q \\ x_C \end{bmatrix}
\hspace{-2bp} \le \hspace{-2bp}
\begin{bmatrix}
-\boldsymbol{\rho}_Q(-[\boldsymbol{F}\boldsymbol{p}]_{\boldsymbol{f}_Q}(\boldsymbol{H}))
\\
-\mathbb{E}\{-f_C(\boldsymbol{p}(\boldsymbol{S}),\boldsymbol{S}) \}
\end{bmatrix}\hspace{-2bp}\triangleq -\boldsymbol{\rho}(-[\boldsymbol{F}\boldsymbol{p}]_{\boldsymbol{f}}(\boldsymbol{H}))
\\
 & \vphantom{\Big|} (\boldsymbol{x}_Q,x_C,\boldsymbol{p})\in{\cal X}\times\Pi
\end{array},
\end{equation*}

\noindent with the identifications $\boldsymbol{f}=(\boldsymbol{f}_Q, f_C)$ and $\mathcal{X} = \mathbb{R}_{+}^{N_Q} \times \{0\} \subset \mathbb{R}^{N= N_Q + 1}$. Consequently, provided integrability of $\boldsymbol{f}(\boldsymbol{p}(\cdot),\cdot)$ for any $\boldsymbol{p} \in \Pi$, \eqref{MAI} exhibits strong duality under mere nonatomicity of the conditional law $\mathrm{P}_{\boldsymbol{U}\mid\boldsymbol{H}}$ (i.e., of the fast-fading given the slow-fading, which holds automatically for continuous small-scale fading) in combination with Slater's constraint qualification; those conditions suffice to ensure that Assumption \ref{assu:Assumption} holds. Notably, the slow-fading (channel) law $\mathrm{P}$ is left \emph{entirely unrestricted}; the required vector measure range convexification is supplied by the inner fast-fading variable $\boldsymbol{U}$, so that \emph{discrete or finite-state} channels (e.g., quantized slow CSI, or finitely many shadowing/LOS--NLOS states) become admissible as well. Let us note that despite using a simple linear utility $g^{o}$, \eqref{MAI} remains a highly nontrivial and challenging nonconvex problem.

\subsubsection{Frequency Division Broadcast Channel}
\par We next consider the optimal risk-aware design of a \textit{frequency division broadcast channel} model, generalizing another risk-neutral problem studied in \citep{Ribeiro2012}. In this case, an Access Point (AP) is tasked to optimally allocate limited resources in order to communicate with $N_Q=N_{\boldsymbol{S}}$ terminals, or users; this sometimes is called a \textit{broadcast} or \textit{downlink} setting, since the AP transmits information to all terminals. Each of the components of $\boldsymbol{S}$ corresponds to the channel strength between the AP and each terminal; we write $H_i$ for the link strength of terminal $i$, so that $\boldsymbol{S}=(H_1,\ldots,H_{N_Q})$. 
However, in this model the AP exploits \textit{frequency division multiplexing}, which allows simultaneous transmission to multiple terminals over different frequencies (carriers, or bands), therefore with no cross-interference. The goal of the AP is to \textit{jointly} select frequencies over which to transmit \textit{and} allocate power to each set of selected transmissions (one for each frequency selected), with \textit{both} decisions made on the basis of observing the channel vector $\boldsymbol{S}$, so as to maximize network-wide QoS, under a total expected power specification $p_{max}>0$.

Power and frequency allocation policies are denoted by $\boldsymbol{p}\in\Pi$ and $\boldsymbol{\phi}\in \Phi_{N_A}$ for some $N_A \le N_Q$, respectively, where $\Pi$ any decomposable set of nonnegative policies, and $\Phi_{N_A}$ is defined as
\[
\Phi_{N_A} = \big\{ 
\boldsymbol{\phi} \in \mathcal{L}_\infty(\mathrm{P}_{\boldsymbol{S}},\mathbb{R}^{N_Q})
\big| \boldsymbol{\phi}(\cdot)\in \{0,1\}^{N_Q},\,\, 
\langle \boldsymbol{\phi}(\cdot), \mathbf{1} \rangle \le N_A,
\,\,\mathrm{P}_{\boldsymbol{S}}\text{-a.e.}
\big\}.
\]
\noindent In other words, each element $\boldsymbol{\phi}\in \Phi_{N_A}$ indicates which transmissions are active ---these are \textit{at most} $N_A$ out of $N_Q$---, each of them realized over the corresponding frequency. It is easy to verify that for any choice of the number of active transmissions $N_A$, the set $\Phi_{N_A}$ is decomposable. By assuming practical Adaptive Modulation and Coding (AMC) with $M\in \mathbb{N}^+$ modes, the service vector stacking the communication rates achievable at each terminal takes the form \citep[Section on ``Frequency Division Broadcast Channel"]{Ribeiro2012}
\[ \boldsymbol{f}_Q([\boldsymbol{\phi}(\boldsymbol{S}),\boldsymbol{p}(\boldsymbol{S})],\boldsymbol{S}) \triangleq \sum_{m=1}^{M}
q_m
\bigg[
\boldsymbol{\phi}(\boldsymbol{S}) \odot
\mathds{1}_{[l_m,l_{m+1})}\left(
\dfrac{
\boldsymbol{S} \odot  \boldsymbol{p}(\boldsymbol{S})
}{\boldsymbol{\nu}}
\right)\hspace{-2pt}\bigg]
, \quad (\boldsymbol{\phi},\boldsymbol{p})\in\Phi_{N_A}\times\Pi,\]

\noindent where $\boldsymbol{\nu}>\mathbf{0}$ contains the variances of the reception noises (possibly) induced at each of the terminals, and where
$q_m\ge0$ is the $m$-th mode rate supported at a terminal, which is achieved whenever the corresponding received signal-to-noise ratio ---i.e., the fractional term inside the indicator function in the expression above--- is between predefined levels $l_m\ge0$ and $l_{m+1}\ge0$ (with $l_m<l_{m+1}$, of course),  for $m\in\mathbb{N}^+_M$; there are $M$ possible operation modes of this form in total. Similarly to Section \ref{MAISec}, we  define the coupling ${{f}_C}([\boldsymbol{\phi}(\boldsymbol{S}),\boldsymbol{p}(\boldsymbol{S})],\boldsymbol{S}) \triangleq  p_{max} - \langle \boldsymbol{\phi}(\boldsymbol{S}) \odot \boldsymbol{p}(\boldsymbol{S}),\mathbf{1} \rangle,(\boldsymbol{\phi},\boldsymbol{p})\in\Phi_{N_A}\times\Pi$. 
\par To recover the two-level structure of \eqref{eq:Base}, we exploit the heterogeneity of the terminals relative to the AP. Without loss of generality, let user $1$ be a terminal whose channel $H_1$ is continuously distributed ---for concreteness, the terminal closest to the AP, which enjoys a strong link ``no matter what''. We single out $H_1$ as the \emph{inner} variable over which the smooth operator averages: in the notation of \eqref{eq:Base}, the observable is the full channel $\boldsymbol{S}$, which we partition as $\boldsymbol{S}=(H_1,\boldsymbol{H}_{-1})$, setting the inner variable $\boldsymbol{U}\triangleq H_1$ and the risk coordinate $\boldsymbol{H}_{-1}\triangleq(H_2,\ldots,H_{N_Q})$ (so that the ``$\boldsymbol{H}$'' of \eqref{eq:Base} is here played by $\boldsymbol{H}_{-1}$, with law $\mathrm{P}\equiv\mathrm{P}_{\boldsymbol{H}_{-1}}$). 
Accordingly, the service of the proximal user is assessed in expectation whereas genuine risk-aversion is retained over the remaining users, i.e., $\boldsymbol{\rho}_Q=(\mathbb{E},\rho_{Q,2},\ldots,\rho_{Q,N_Q})$; this is another, purely informational (rather than operational) mechanism ---distinct from a physical fast/slow fading split of Section \ref{MAISec}--- by which the two-level form of \eqref{eq:Base} arises.

Lastly, we consider a strictly increasing concave utility $u$ evaluating each component of the (risk-)ergodic vector $\boldsymbol{x} \in \mathbb{R}^{N_Q}$ ---we succintly write $\langle u(\boldsymbol{x}),\mathbf{1}\rangle\triangleq g^o(\boldsymbol{x})$---, the latter are further constrained within the compact box $\mathcal{B}\triangleq \big\{ \boldsymbol{x}\in\mathbb{R}^{N_Q} \big| x_{min} \mathbf{1} \le \boldsymbol{x}\le x_{max} \mathbf{1} \big\}$, for fixed $0 \le x_{min} < x_{max}$. Then, the corresponding stochastic resource allocation problem may be formulated as
\begin{tagequation}{FDB}
\begin{array}{rl}
\underset{\boldsymbol{x},\boldsymbol{\phi}(\cdot),\boldsymbol{p}(\cdot)}{\mathrm{maximize}} & \langle u(\boldsymbol{x}),\mathbf{1}\rangle\\
\mathrm{subject\,to} & \boldsymbol{x} \le
\displaystyle{
-\boldsymbol{\rho}_{Q}\Bigg[\hspace{-2bp}-\hspace{-1bp}
\mathbb{E}\Bigg\{\sum_{m=1}^{M}
q_m
\Bigg[
\boldsymbol{\phi}(\boldsymbol{S}) \odot
\mathds{1}_{[l_m,l_{m+1})}\left(
\dfrac{
\boldsymbol{S} \odot  \boldsymbol{p}(\boldsymbol{S})
}{\boldsymbol{\nu}}
\right)\hspace{-2pt}\Bigg]
\Bigg|\boldsymbol{H}_{-1}\Bigg\}
\Bigg]}
\\
\vspace{-8bp}
\\
& \langle \mathbb{E}\{ 
\boldsymbol{\phi}(\boldsymbol{S})
\odot
\boldsymbol{p}(\boldsymbol{S}) \} , \mathbf{1} \rangle
\le p_{max}
\\
 &  (\boldsymbol{x}, \boldsymbol{\phi}, \boldsymbol{p})\in \mathcal{B} \times\Phi_{N_A} \times \Pi
\end{array},\label{FDB}
\end{tagequation}

\noindent where $\boldsymbol{\rho}_Q$ is a finite-valued vector risk measure on $\mathcal{L}_p(\mathrm{P}_{\boldsymbol{H}_{-1}},\mathbb{R}^{N_Q})$ ---with first component $\rho_{Q,1}=\mathbb{E}$, so that the near user's row reads $x_1\le\mathbb{E}\{f_{Q,1}\}$--- and compatible with the construction of \eqref{eq:Base}, as in Section \ref{MAISec}.
We observe that \eqref{FDB} is a highly nonconvex functional mixed-integer program; in fact, the policy $\boldsymbol{\phi}$ is integer-valued, and the (instantaneous) service function inside the risk measure $\boldsymbol{\rho}_Q$ is clearly discontinuous.
 
Again, problem \eqref{FDB} is a special case of \eqref{eq:Base}, because we can restate the former as
\begin{equation*}
\begin{array}{rl}
\underset{(\boldsymbol{x}_Q,x_C),\boldsymbol{p}'(\cdot)}{\mathrm{maximize}} & g^{o}(\boldsymbol{x}_Q)\\
\mathrm{subject\,to}\,\hspace{1.4bp} & \begin{bmatrix}\boldsymbol{x}_Q \\ x_C \end{bmatrix}
\hspace{-2bp} \le \hspace{-2bp}
\begin{bmatrix}
-\boldsymbol{\rho}_Q(-[\boldsymbol{F}\boldsymbol{p}']_{\boldsymbol{f}_Q}(\boldsymbol{H}_{-1}))
\\
-\mathbb{E}\{-f_C(\boldsymbol{p}'(\boldsymbol{S}),\boldsymbol{S}) \}
\end{bmatrix}\hspace{-2bp}\triangleq -\boldsymbol{\rho}(-[\boldsymbol{F}\boldsymbol{p}']_{\boldsymbol{f}}(\boldsymbol{H}_{-1}))
\\
 & \vphantom{\Big|} (\boldsymbol{x}_Q,x_C,\boldsymbol{p}')\in{\cal X}\times\Pi'
\end{array},
\end{equation*}
with the identifications $\boldsymbol{p}'=(\boldsymbol{\phi},\boldsymbol{p})$, $\boldsymbol{f}=(\boldsymbol{f}_Q, f_C)$, $\Pi'=\Phi_{N_A} \times \Pi$ and $\mathcal{X} = \mathcal{B} \times \{0\} \subset \mathbb{R}^{N= N_Q + 1}$, and where $[\boldsymbol{F}\boldsymbol{p}']_{\boldsymbol{f}_Q}(\boldsymbol{H}_{-1})=\mathbb{E}\{\boldsymbol{f}_Q(\boldsymbol{p}'(\boldsymbol{S}),\boldsymbol{S})\,|\,\boldsymbol{H}_{-1}\}$ is the ergodic (proximal user-averaged) service. 
Here $\mathcal{S}$ is the \emph{full} channel space (in which the observable $\boldsymbol{S}$ takes values), so the product $\Pi'$ ---already decomposable over the full channel--- is decomposable over $\mathscr{B}(\mathcal{S})$ exactly as required. 
If $\boldsymbol{f}(\boldsymbol{p}'(\cdot),\cdot)$ is integrable on $\Pi'$ (automatic here, as $\boldsymbol{f}_Q$ is bounded) and the channel of the proximal user is conditionally nonatomic given the others ---i.e., $\mathrm{P}_{H_1\mid\boldsymbol{H}_{-1}}$ is nonatomic $\mathrm{P}$-a.e., which holds whenever $H_1$ is continuously distributed--- it follows that problem \eqref{FDB} exhibits strong duality as long as it is strictly feasible.
Notably, the channels of the remaining users $\boldsymbol{H}_{-1}$ may then follow an \emph{arbitrary} (possibly discrete) law. This is true despite the discontinuity of $\boldsymbol{f}(\cdot,\boldsymbol{S})$ and the nonconvexity of the set of feasible frequency policies $\Phi_{N_A}$. 
\subsection{\label{subsec: Constrained Learning}Risk-Constrained Learning with Nonconvex Losses}
\par Next, we consider a general formulation in the context of supervised risk-constrained learning, where the associated loss functions are allowed to be nonconvex. In the risk-neutral setting, this class of problems has been recently studied in \cite{Chamon2021}. On our usual probability space $(\Omega,\mathscr{F},\mu)$, we consider random \textit{example vectors} $(\boldsymbol{X}_i,Y_i) \colon \Omega \rightarrow \mathbb{R}^d \times \mathbb{R}$, for $i \in \mathbb{N}_m$, together with their induced Borel probability distributions $\mathcal{D}_i:\mathscr{B}(\mathbb{R}^d \times \mathbb{R}) \rightarrow [0,1]$, instantiated over data pairs $(\boldsymbol{x},y)$, where $\boldsymbol{x} \in \mathbb{R}^d$ represents a realized feature or system input and $y \in \mathbb{R}$ represents a realized label or measurement.
We denote by $\mathcal{D}_{\boldsymbol{X}_i}$ the marginal distribution of feature $\boldsymbol{X}_i$, and by $\mathcal{D}_{Y_i\vert \boldsymbol{X}_i}$ the conditional distribution of label $Y_i$ given  feature $\boldsymbol{X}_i$, for $i \in \mathbb{N}_m$. 
We take as reference the (normalized) \emph{mixture of marginals}
\vspace{-7pt}
\[
\mathcal{D}_{\Sigma} \triangleq \dfrac{1}{m+1}\sum_{i=0}^{m} \mathcal{D}_{\boldsymbol{X}_i},
\]
\vspace{-7pt}

\noindent serving as a common dominating measure in our development
below (cf. \eqref{eq:RCLearning0}). By construction, $\mathcal{D}_{\boldsymbol{X}_i} \le (m+1)\,\mathcal{D}_{\Sigma}$ as
measures, and $\mathcal{D}_{\boldsymbol{X}_i} \ll \mathcal{D}_{\Sigma}$, for
every $i \in \mathbb{N}_m$. Radon--Nikodym yield densities
$w_i \triangleq \mathrm{d}\mathcal{D}_{\boldsymbol{X}_i}/\mathrm{d}\mathcal{D}_{\Sigma} \in
\mathcal{L}_1(\mathcal{D}_{\Sigma},\mathbb{R})$ with $0 \le w_i \le m+1$
$\mathcal{D}_{\Sigma}$-a.e., whose pointwise sum satisfies $\sum_{i=0}^{m} w_i(\boldsymbol{x}) = m+1$, for $\mathcal{D}_{\Sigma}\text{-a.e. } \boldsymbol{x} \in \mathbb{R}^d.$
We write $\boldsymbol{X}_{\Sigma}$ for a generic feature with law $\mathcal{D}_{\Sigma}$.
\par To obtain the two-level structure of \eqref{eq:Base}, we single out one
continuously distributed feature coordinate to serve as the inner (averaging)
variable (cf. smooth operator). Without loss of generality we take this to be the $j$-th coordinate, and
split feature vectors as $\boldsymbol{X}_i = (\boldsymbol{X}_i^{(-j)}, X_i^{(j)})$,
$i \in \mathbb{N}_m$, where $X_i^{(j)}$ is the singled-out coordinate and
$\boldsymbol{X}_i^{(-j)}$ collects the remaining $d-1$ coordinates.
\par Our formulation rests on the assumption of \emph{conditional nonatomicity of a
single feature $j$}, imposed per task $i$: for each $i \in \mathbb{N}_m$, the conditional law
$\mathcal{D}_{X_i^{(j)}\mid \boldsymbol{X}_i^{(-j)}}$ is nonatomic
($\mathcal{D}_{\boldsymbol{X}_i^{(-j)}}$-a.e.).
No assumptions are imposed on
$\boldsymbol{X}_i^{(-j)}$, whose law may be \emph{arbitrary}. The outer risk in each task acts on a conditional
mean loss over $X_i^{(j)}$; for simplicity in this treatment, for each $i \in \mathbb{N}_m$, we assume a convex, lower semicontinuous, positively homogeneous risk measure
$\rho_i^{(-j)} \colon \mathcal{L}_1(\mathcal{D}_{\boldsymbol{X}_i^{(-j)}},\mathbb{R}) \to \mathbb{R}$
such that
$\rho_i = \rho_i^{(-j)} \circ \mathbb{E}\{\cdot | \boldsymbol{X}_i^{(-j)}\}$.
Genuine
risk-aversion is thus retained over $\boldsymbol{X}_i^{(-j)}$, while the continuous coordinate
$X_i^{(j)}$ is treated in expectation. As we discuss below, conditional
nonatomicity of a single feature supplies the convexification underlying Theorems
\ref{thm:Main} and \ref{thm:Main2}; 
the risk-neutral case $\rho_i = \mathbb{E}$ may be recovered
automatically under our framework (see also Section \ref{subsec: risk-neutral reduction} for a more general argument).
\par Letting $\ell_i \colon \mathbb{R} \times \mathbb{R} \rightarrow \mathbb{R}_+$ be
given (possibly nonconvex) loss functions, for $i \in \mathbb{N}_m$, we consider the
risk-over-risk functional constrained learning problem
\begin{tagequation}{RCL}
\begin{array}{rl}
\:\underset{\boldsymbol{f}(\cdot) \in \mathcal{F}}{\mathrm{minimize}} &{\rho}_{0}\big(\widehat{{\rho}}_{0}\big( \ell_0\big(\boldsymbol{f}(\boldsymbol{X}_0),Y_0\big)\big\vert \boldsymbol{X}_0\big) \big) \\
\mathrm{subject\,to} &  {\rho}_{i}\big(\widehat{{\rho}}_{i}\big( \ell_i\big(\boldsymbol{f}(\boldsymbol{X}_i),Y_i\big)\big\vert \boldsymbol{X}_i\big) \big) \leq c_i, \quad i \in \mathbb{N}_m^+
\end{array},\label{eq:RCLearning}
\end{tagequation}

\noindent where, for each $i \in \mathbb{N}_m$, ${\rho}_{i}$ is a convex, lower
semicontinuous, and (for simplicity) positively homogeneous risk measure taking finite
values on $\mathcal{L}_1(\mathcal{D}_{\boldsymbol{X}_i},\mathbb{R})$ as described above, while
$\widehat{{\rho}}_{i}(\cdot\vert \boldsymbol{X}_i)$ is a \emph{conditional risk measure}
over the law $\mathcal{D}_{Y_i\mid\boldsymbol{X}_i}$ of $Y_i$ given
$\boldsymbol{X}_i$, defined as
$\widehat{{\rho}}_{i}(\cdot\vert \boldsymbol{X}_i) \colon \mathcal{L}_{p'}(\mathcal{D}_i,\mathbb{R}) \rightarrow \mathcal{L}_p(\mathcal{D}_{\boldsymbol{X}_i},\mathbb{R})$,
for some $p, p' \in [1,\infty]$, and such that it obeys a \textit{substitution rule} of
the form
\[ \widehat{{\rho}}_{i}( \ell_i(\boldsymbol{f}(\boldsymbol{X}_i),Y_i)\vert \boldsymbol{X}_i)  \equiv \widehat{{\rho}}_{i}( \ell_i(\boldsymbol{z},Y_i)\vert \boldsymbol{X}_i)\big\vert_{\boldsymbol{z} = \boldsymbol{f}(\boldsymbol{X}_i)} \equiv \Lambda_i\left(\boldsymbol{f}(\boldsymbol{X}_i),\boldsymbol{X}_i\right),\]
\noindent where $\Lambda_i(\boldsymbol{f}(\cdot),\cdot) \in \mathcal{L}_{p}(\mathcal{D}_{\boldsymbol{X}_i},\mathbb{R})$
for any \textit{learning representation} (i.e., policy)
$\boldsymbol{f} \colon \mathbb{R}^d \rightarrow \mathbb{R}$ lying in an appropriate
densely decomposable space (relative to some norm $\|\cdot\|$) $\mathcal{F}$ (cf.
Definition \ref{def: densely decomposable set}); also $c_i \in \mathbb{R}$ for all
$i \in \mathbb{N}_m^+$. The composed map $\Lambda_i(\boldsymbol{f}(\cdot),\cdot)$ is the
\emph{superposition} (Nemytskii) operator associated with $\boldsymbol{f}$ and  $\Lambda_i$; casting \eqref{eq:RCLearning} in terms of such operators
is what renders it an instance of \eqref{eq:Base}.
\par Next, let $\mathds{A}^{i} \subseteq \mathcal{L}_{\infty}(\mathcal{D}_{\boldsymbol{X}_i^{(-j)}},\mathbb{R})$
be the risk envelope of $\rho_i^{(-j)}$, $i \in \mathbb{N}_m$, so that,
writing out the conditional mean 
$\mathbb{E}\{Z | \boldsymbol{X}_i^{(-j)}\!=\!\boldsymbol{x}^{(-j)}\}
= \int_{\mathbb{R}} Z(\boldsymbol{x}^{(-j)},t)\,
\mathrm{d}\mathcal{D}_{X_i^{(j)}\mid \boldsymbol{X}_i^{(-j)}}(t | \boldsymbol{x}^{(-j)})$ over the singled-out coordinate,
the (nested) dual representation of $\rho_i$ reads, for $ Z \in \mathcal{L}_1(\mathcal{D}_{\boldsymbol{X}_i},\mathbb{R})$, as
\begin{equation*}
\rho_i(Z) = \sup_{\zeta \in \mathds{A}^i}
\int_{\mathbb{R}^{d-1}} \zeta(\boldsymbol{x}^{(-j)})
\left[\int_{\mathbb{R}} Z(\boldsymbol{x}^{(-j)},t)\,
\mathrm{d}\mathcal{D}_{X_i^{(j)}\mid \boldsymbol{X}_i^{(-j)}}(t | \boldsymbol{x}^{(-j)})\right]
\mathrm{d}\mathcal{D}_{\boldsymbol{X}_i^{(-j)}}(\boldsymbol{x}^{(-j)}).
\end{equation*}
\par Since each conditional law
$\mathcal{D}_{X_i^{(j)}\mid\boldsymbol{X}_i^{(-j)}}$ is nonatomic and the mixture
conditional is the corresponding convex combination of the task conditionals, the
conditional law $\mathcal{D}_{X_{\Sigma}^{(j)}\mid\boldsymbol{X}_{\Sigma}^{(-j)}}$ is
itself nonatomic\footnote{Writing $\nu_i \triangleq \mathcal{D}_{\boldsymbol{X}_i^{(-j)}}$ and
$\kappa_i \triangleq \mathcal{D}_{X_i^{(j)}\mid\boldsymbol{X}_i^{(-j)}}$, so that
$\mathcal{D}_{\boldsymbol{X}_i} = \nu_i \otimes \kappa_i$ and
$\mathcal{D}_{\Sigma} = \nu_{\Sigma}\otimes\kappa_{\Sigma}$, linearity of the pushforward gives
$\nu_{\Sigma} = \frac{1}{m+1}\sum_{i} \nu_i$, whence $\nu_i \ll \nu_{\Sigma}$; setting
$a_i \triangleq \frac{1}{m+1}\,\mathrm{d}\nu_i/\mathrm{d}\nu_{\Sigma}$ we have $a_i \ge 0$ and
$\sum_{i=0}^{m} a_i(\boldsymbol{v}) = 1$ for $\nu_{\Sigma}$-a.e. $\boldsymbol{v}$. Then, for every
$A \in \mathscr{B}(\mathbb{R}^{d-1})$ and $B \in \mathscr{B}(\mathbb{R})$,
\[
\int_A \kappa_{\Sigma}(B\mid\boldsymbol{v})\,\nu_{\Sigma}(\mathrm{d}\boldsymbol{v})
= \mathcal{D}_{\Sigma}(A\times B)
= \frac{1}{m+1}\sum_i \int_A \kappa_i(B\mid\boldsymbol{v})\,\nu_i(\mathrm{d}\boldsymbol{v})
= \int_A \sum_i a_i(\boldsymbol{v})\,\kappa_i(B\mid\boldsymbol{v})\,\nu_{\Sigma}(\mathrm{d}\boldsymbol{v}),
\]
so, $A$ being arbitrary, $\kappa_{\Sigma}(B\mid\boldsymbol{v}) = \sum_i a_i(\boldsymbol{v})
\kappa_i(B\mid\boldsymbol{v})$ for $\nu_{\Sigma}$-a.e. $\boldsymbol{v}$, for each fixed $B$, the
exceptional null set depending on $B$. Choosing $B = (-\infty,q]$, $q \in \mathbb{Q}$, and
discarding the countably many exceptional sets simultaneously, the cdfs of
$\kappa_{\Sigma}(\cdot\mid\boldsymbol{v})$ and of $\sum_i a_i(\boldsymbol{v})\kappa_i(\cdot\mid\boldsymbol{v})$
agree at all rationals, for $\nu_{\Sigma}$-a.e. $\boldsymbol{v}$ . By right-continuity the cdfs agree
everywhere, so the two measures coincide. Finally, $\nu_i(N_i) = 0$, with $N_i$ the set where
$\kappa_i$ has an atom, forces $a_i = 0$ $\nu_{\Sigma}$-a.e. on $N_i$; hence for
$\nu_{\Sigma}$-a.e. $\boldsymbol{v}$ every summand with $a_i(\boldsymbol{v}) > 0$ is nonatomic and so
is $\kappa_{\Sigma}(\cdot\mid\boldsymbol{v})$. Since $a_i > 0$ $\nu_i$-a.e., the converse also
holds, so the two conditions are equivalent.}; this is precisely the nonatomicity of $\mathrm{P}_{\boldsymbol{U}\mid\boldsymbol{H}}$
required by \eqref{eq:Base} (under the identification $\boldsymbol{S}\equiv\boldsymbol{X}_{\Sigma}$,
$\boldsymbol{U}\equiv X_{\Sigma}^{(j)}$, $\boldsymbol{H}\equiv\boldsymbol{X}_{\Sigma}^{(-j)}$
introduced below). Then, for each $i \in \mathbb{N}_m$, we construct a convex, lower semicontinuous,
positively homogeneous risk measure $\widetilde{\rho}_i$ on $\mathcal{L}_1(\mathcal{D}_{\Sigma},\mathbb{R})$
that retains the same two-level structure over the mixture reference. Specifically, for $ W \in \mathcal{L}_1(\mathcal{D}_{\Sigma},\mathbb{R})$, we have
\begin{equation*}
\widetilde{\rho}_i(W) \triangleq \sup_{\tilde{\zeta} \in \widetilde{\mathds{A}}^i}
\int_{\mathbb{R}^{d-1}} \tilde{\zeta}(\boldsymbol{x}^{(-j)})
\left[\int_{\mathbb{R}} W(\boldsymbol{x}^{(-j)},t)\,
\mathrm{d}\mathcal{D}_{X_\Sigma^{(j)}\mid \boldsymbol{X}_\Sigma^{(-j)}}(t\mid \boldsymbol{x}^{(-j)})\right]
\mathrm{d}\mathcal{D}_{\boldsymbol{X}_\Sigma^{(-j)}}(\boldsymbol{x}^{(-j)}),
\end{equation*}
with $\widetilde{\mathds{A}}^i \subseteq \mathcal{L}_{\infty}(\mathcal{D}_{\boldsymbol{X}_\Sigma^{(-j)}},\mathbb{R})$;
in operator form, $\widetilde{\rho}_i = \widetilde{\rho}_i^{(-j)} \circ
\mathbb{E}\{\,\cdot \mid \boldsymbol{X}_\Sigma^{(-j)}\}$, the inner integral being the conditional
expectation under $\mathcal{D}_{\Sigma}$.
The key identity is that, evaluated on the reweighted (Nemytskii) argument $W = Z w_i$
with $Z \in \mathcal{L}_1(\mathcal{D}_{\boldsymbol{X}_i},\mathbb{R})$, the density is absorbed by the
inner integral: setting $w_i^{(-j)} \triangleq \mathrm{d}\mathcal{D}_{\boldsymbol{X}_i^{(-j)}}/
\mathrm{d}\mathcal{D}_{\boldsymbol{X}_\Sigma^{(-j)}}$, the disintegration chain rule gives
\begin{equation*}
\begin{split}
&\int_{\mathbb{R}} Z(\boldsymbol{x}^{(-j)},t)\, w_i(\boldsymbol{x}^{(-j)},t)\,
\mathrm{d}\mathcal{D}_{X_\Sigma^{(j)}\mid \boldsymbol{X}_\Sigma^{(-j)}}(t\mid \boldsymbol{x}^{(-j)})
\\
&\quad = w_i^{(-j)}(\boldsymbol{x}^{(-j)}) \int_{\mathbb{R}} Z(\boldsymbol{x}^{(-j)},t)\,
\mathrm{d}\mathcal{D}_{X_i^{(j)}\mid \boldsymbol{X}_i^{(-j)}}(t\mid \boldsymbol{x}^{(-j)}),
\end{split}
\end{equation*}
while the outer part transports the envelope through $\tilde{\zeta}\, w_i^{(-j)}\,
\mathrm{d}\mathcal{D}_{\boldsymbol{X}_\Sigma^{(-j)}} = \tilde{\zeta}\,\mathrm{d}\mathcal{D}_{\boldsymbol{X}_i^{(-j)}}$.
Choosing $\widetilde{\mathds{A}}^i$ to agree with $\mathds{A}^i$ up to $\mathcal{D}_{\boldsymbol{X}_i^{(-j)}}$-null
sets, e.g.,
$\widetilde{\mathds{A}}^i = \{\tilde{\zeta} \in \mathcal{L}_{\infty}(\mathcal{D}_{\boldsymbol{X}_\Sigma^{(-j)}},\mathbb{R})
\mid \exists\, \zeta \in \mathds{A}^i,\ \tilde{\zeta} = \zeta \mathds{1}_{\{w_i^{(-j)}>0\}},\ \mathcal{D}_{\boldsymbol{X}_{\Sigma}^{(-j)}}\text{-a.e.}\}$, both displays coincide and
\vspace{-6.7pt}
\begin{equation}\nonumber
\widetilde{\rho}_i\big(Z\, w_i\big) = \rho_i(Z), \qquad Z \in \mathcal{L}_1(\mathcal{D}_{\boldsymbol{X}_i},\mathbb{R}),\ i \in \mathbb{N}_m.
\end{equation}
Note that, in fact, $\widetilde{\rho}_i$ is finite on all of $\mathcal{L}_1(\mathcal{D}_{\Sigma},\mathbb{R})$, for each $i \in \mathbb{N}_m$.
Therefore, from the superposition form
$G_i(\boldsymbol{f}(\cdot),\cdot) \triangleq \Lambda_i(\boldsymbol{f}(\cdot),\cdot)\, w_i
\in \mathcal{L}_1(\mathcal{D}_{\Sigma},\mathbb{R})$, it follows that, for every
$\boldsymbol{f}\in \mathcal{F}$ and every $i \in \mathbb{N}_m$,
\begin{equation*}
\begin{split}
{\rho}_{i}\big(\widehat{{\rho}}_{i}\big( \ell_i\big(\boldsymbol{f}(\boldsymbol{X}_i),Y_i\big)\big\vert \boldsymbol{X}_i\big) \big)
&= {\rho}_{i}\big( \Lambda_i(\boldsymbol{f}(\boldsymbol{X}_i),\boldsymbol{X}_i)\big)
= \widetilde{\rho}_{i}\big( G_i(\boldsymbol{f}(\boldsymbol{X}_{\Sigma}),\boldsymbol{X}_{\Sigma})\big)\\[2pt]
&= \widetilde{\rho}_i^{(-j)}\big(
\underbrace{\mathbb{E}\big\{ G_i(\boldsymbol{f}(\boldsymbol{X}_{\Sigma}),\boldsymbol{X}_{\Sigma}) \,\big\vert\, \boldsymbol{X}_{\Sigma}^{(-j)}\big\}}
_{\textstyle \triangleq\, [\boldsymbol{F}\boldsymbol{f}]_{G_i}(\boldsymbol{X}_{\Sigma}^{(-j)})}
\big).
\end{split}
\end{equation*}
\noindent The operator $[\boldsymbol{F}\boldsymbol{f}]_{G_i}$ defined above is precisely the smooth
operator of \eqref{eq:Base}: it averages the instantaneous score $G_i$ over the singled-out
continuous coordinate, conditionally on the remaining features, thus returning a function of
the latter alone. Its notation mirrors that of $[\boldsymbol{F}\boldsymbol{p}]_{\boldsymbol{f}}$
in \eqref{eq:Base}. We can then equivalently
recast \eqref{eq:RCLearning} as
\begin{tagequation}{RCL$'$}
\begin{array}{rl}
\:\underset{\boldsymbol{f}(\cdot) \in \mathcal{F}}{\mathrm{minimize}} &\widetilde{{\rho}}_{0}\left(G_0\big(\boldsymbol{f}(\boldsymbol{X}_{\Sigma}),\boldsymbol{X}_{\Sigma}\big) \right) \\
\mathrm{subject\,to} &  \widetilde{{\rho}}_{i}\left( G_i\big(\boldsymbol{f}(\boldsymbol{X}_{\Sigma}),\boldsymbol{X}_{\Sigma}\big)\right) \leq c_i, \quad i \in \mathbb{N}_m^+
\end{array}.\label{eq:RCLearning0}
\end{tagequation}

\par Altogether, \eqref{eq:RCLearning0} is readily an instance of \eqref{eq:Base} (cf. Section
\ref{subsec:theory extension}), under the identification of the observable with the (mixture)
reference feature, $\boldsymbol{S}\equiv\boldsymbol{X}_{\Sigma}$, of the inner variable with the
singled-out continuous coordinate, $\boldsymbol{U}\equiv X_{\Sigma}^{(j)}$, and of the risk
coordinate with the remaining features, $\boldsymbol{H}\equiv\boldsymbol{X}_{\Sigma}^{(-j)}$;
here the $G_i$, $i\in\mathbb{N}_m$, serve as the instantaneous performance scores, the policy
space $\Pi$ is identified with $\mathcal{F}$, and the outer risks $\widetilde{\rho}_i^{(-j)}$
act over $\boldsymbol{H}$.
Consequently, Slater's condition, alongside conditions 6$'$.(a)--(b) of Assumption
\ref{assu:Assumption2}, suffice to ensure that \eqref{eq:RCLearning0} exhibits strong duality
by Theorem \ref{thm:Main2}, and the same of course holds for \eqref{eq:RCLearning}. Following our discussion in Section
\ref{subsec: densely decomposable sets}, conditions 6$'$.(a)--(b) of Assumption
\ref{assu:Assumption2} require integrability of $G_i(\boldsymbol{f}(\cdot),\cdot)$ and
pointwise continuity of
$\widetilde{\rho}_i^{(-j)}\!\big(\mathbb{E}\{G_i(\boldsymbol{f}(\boldsymbol{X}_{\Sigma}),\boldsymbol{X}_{\Sigma})\mid\boldsymbol{X}_{\Sigma}^{(-j)}\}\big)$,
for $i \in \mathbb{N}_m$, at each $\boldsymbol{f} \in \mathcal{F}_{\|\cdot\|}$ (where, for example,
$\mathcal{F}_{\|\cdot\|}$ may be taken to be a decomposable space of which $\mathcal{F}$ is a
dense subset). These requirements are fairly minor and can be easily enforced by an
appropriate choice of $\mathcal{F}$ (again, for a detailed list of such examples, we refer
the reader to the discussion in Section \ref{subsec: densely decomposable sets}), under
also minor assumptions on the (possibly nonconvex) loss functions $\ell_i$ and on the
associated conditional risk measures $\widehat{\rho}_i(\cdot \vert \boldsymbol{X}_i)$. In
passing, it is worth noting that strong duality holds for \eqref{eq:RCLearning0} (and thus
also for \eqref{eq:RCLearning}) under merely Slater's condition using Theorem
\ref{thm:Main}, if we instead assume that $\mathcal{F}$ is decomposable to begin with,
although this result might be weaker in terms of its practical implications.

\par Having described \eqref{eq:RCLearning} in full generality, we now comment on the scope of
our assumptions and on specific instances subsumed by our framework. The absolute continuity
$\mathcal{D}_{\boldsymbol{X}_i}\ll\mathcal{D}_{\Sigma}$ noted above imposes \emph{no} restriction on
the tasks whatsoever ---in particular, the task marginals may even have disjoint supports---
thus dispensing with the common-denominator dominance requirement often imposed in the
literature (see, e.g., \citep{Chamon2021}). The only genuine data assumption is therefore the
conditional nonatomicity of a single feature, which is standard and mild: it holds, for
instance, whenever the singled-out feature admits a conditional density relative to the
(one-dimensional) Lebesgue measure in each task (see, e.g., \citep{Chamon2021}). Crucially, the remaining features are left entirely unrestricted, and may well be
discrete or categorical; this directly accommodates the mixed continuous/categorical feature
vectors common in machine learning applications (especially in practice). 
\par Next, let us discuss the generality of the conditional risk mappings $\widehat{\rho}_i(\cdot\vert \boldsymbol{X}_i)$, which are defined by enforcing the substitution rule. In fact, the substitution rule is quite general and holds beyond classes of, e.g., convex, lower semicontinuous, and/or positively homogeneous risk measures. Indeed, several compositional risk measures studied in \citep{Kalogerias2018b} satisfy the substitution rule, being in particular nonhomogeneous. More generally, a wide range of conditional risk measures constructed using conditional expectations (as building blocks) can be shown to be compatible with this definition.

\par Finally, let us specialize our result to the case where $\rho_i (\cdot )= \mathbb{E}_{\mathcal{D}_{\boldsymbol{X}_i}}\{\cdot\}$ and $\widehat{\rho}_i(\cdot\vert \boldsymbol{X}_i) = \mathbb{E}_{\mathcal{D}_{Y_i \vert \boldsymbol{X}_i}}\{\cdot\}$, for all $i \in \mathbb{N}_m$ (see also Section \ref{Counter}). This setting corresponds to risk-neutral constrained learning considered in \cite{Chamon2021}. Here conditions 6$'$.(a)--(b) of Assumption \ref{assu:Assumption2} reduce to $\mathbb{E}\left\{ \ell_i\left( \boldsymbol{f}(\boldsymbol{X}_i),Y_i\right) \vert \boldsymbol{X}_i\right\} \in \mathcal{L}_1(\mathcal{D}_{\boldsymbol{X}_i},\mathbb{R})$ for every $\boldsymbol{f} \in \mathcal{F}_{\|\cdot\|}$, so that, under our standing assumptions together with Slater's constraint qualification we recover the results of \citep[Proposition III.2 and Proposition B.1]{Chamon2021}, unifying the classification and regression regimes, while dispensing certain core assumptions. In particular, the zero duality gap result on regression problems given in \citep[Proposition B.1]{Chamon2021} is shown here to hold without the requirement postulated in \citep[Assumption 6]{Chamon2021}. More importantly, we dispense with the decomposability assumption (also employed in
\cite{Chamon2021}), allowing us to directly study the parametrized problem ---assuming,
for example, that $\mathcal{F}$ is a sufficiently expressive function class (e.g., a neural
network) yielding a densely decomposable space; see Section
\ref{subsec: densely decomposable sets} for an example of such a \emph{finite-dimensional}
universal class. Another example is a \emph{universal} RKHS. This is in contrast to the results of \cite{Chamon2021}, or any other similar result appearing in the literature, which can only guarantee an approximately zero duality gap (under additional assumptions) using ``$\varepsilon$-universal" policy parametrizations. Thus, our result strengthens substantially the alternatives available in the literature even in the risk-neutral setting, pushing the frontiers of strong duality in the context of nonconvex constrained learning in multiple ways.

\vspace{-2pt}
\section{Conclusions} \label{sec: Conclusions}
We established strong duality for a wide class of risk-constrained nonconvex functional programs. Our technical approach has exploited the Kingman-Robertson extension of Lyapunov's convexity theorem for infinite-dimensional vector-measures, and is applicable to programs involving arbitrary convex risk functionals on general $\mathcal{L}_p,p\in[1,\infty]$ spaces. We establish two distinct strong duality theorems. The former strictly generalizes existing results available for the risk-neutral setting, essentially without imposing additional assumptions. The latter goes beyond the state of the art, bypassing decomposability ---a core assumption in all existing relevant literature--- and establishing strong duality for risk-constrained problems over structured (e.g., continuous, smooth or RKHS) policies, which were previously out of reach. We further discussed the necessity of our optimization model by presenting explicitly constructed counterexamples, as well as extensions of our results covering a wider class of problems.
Lastly, we showcased the applicability of the theory on two application areas with practical interest, namely, risk-constrained wireless systems resource allocation, and risk-constrained learning with nonconvex losses. In the latter case, we recovered and extended known strong duality results for risk-neutral constrained learning under relaxed assumptions, while also emphasizing the expressive power of the risk-constrained model studied herein.
\par There are several avenues for further investigation. An example, going beyond the \emph{static} policy search setting studied herein, is the establishment of strong duality relations in \emph{dynamic} risk-constrained nonconvex stochastic programming models, under minimal and realistic assumptions.

\vspace{-3pt}
\bmhead{Acknowledgments} \blindblock{The authors would like to acknowledge support by a Microsoft gift and by the US NSF under
Grant CCF 2242215. The authors also acknowledge two anonymous reviewers for their comments and suggestions, which have significantly improved the breadth of the results.}
\vspace{-4pt}

\bibliography{library}


\begin{thebibliography}{64}
\ifx \bisbn   \undefined \def \bisbn  #1{ISBN #1}\fi
\ifx \binits  \undefined \def \binits#1{#1}\fi
\ifx \bauthor  \undefined \def \bauthor#1{#1}\fi
\ifx \batitle  \undefined \def \batitle#1{#1}\fi
\ifx \bjtitle  \undefined \def \bjtitle#1{#1}\fi
\ifx \bvolume  \undefined \def \bvolume#1{\textbf{#1}}\fi
\ifx \byear  \undefined \def \byear#1{#1}\fi
\ifx \bissue  \undefined \def \bissue#1{#1}\fi
\ifx \bfpage  \undefined \def \bfpage#1{#1}\fi
\ifx \blpage  \undefined \def \blpage #1{#1}\fi
\ifx \burl  \undefined \def \burl#1{\textsf{#1}}\fi
\ifx \doiurl  \undefined \def \doiurl#1{\url{https://doi.org/#1}}\fi
\ifx \betal  \undefined \def \betal{\textit{et al.}}\fi
\ifx \binstitute  \undefined \def \binstitute#1{#1}\fi
\ifx \binstitutionaled  \undefined \def \binstitutionaled#1{#1}\fi
\ifx \bctitle  \undefined \def \bctitle#1{#1}\fi
\ifx \beditor  \undefined \def \beditor#1{#1}\fi
\ifx \bpublisher  \undefined \def \bpublisher#1{#1}\fi
\ifx \bbtitle  \undefined \def \bbtitle#1{#1}\fi
\ifx \bedition  \undefined \def \bedition#1{#1}\fi
\ifx \bseriesno  \undefined \def \bseriesno#1{#1}\fi
\ifx \blocation  \undefined \def \blocation#1{#1}\fi
\ifx \bsertitle  \undefined \def \bsertitle#1{#1}\fi
\ifx \bsnm \undefined \def \bsnm#1{#1}\fi
\ifx \bsuffix \undefined \def \bsuffix#1{#1}\fi
\ifx \bparticle \undefined \def \bparticle#1{#1}\fi
\ifx \barticle \undefined \def \barticle#1{#1}\fi
\bibcommenthead
\ifx \bconfdate \undefined \def \bconfdate #1{#1}\fi
\ifx \botherref \undefined \def \botherref #1{#1}\fi
\ifx \url \undefined \def \url#1{\textsf{#1}}\fi
\ifx \bchapter \undefined \def \bchapter#1{#1}\fi
\ifx \bbook \undefined \def \bbook#1{#1}\fi
\ifx \bcomment \undefined \def \bcomment#1{#1}\fi
\ifx \oauthor \undefined \def \oauthor#1{#1}\fi
\ifx \citeauthoryear \undefined \def \citeauthoryear#1{#1}\fi
\ifx \endbibitem  \undefined \def \endbibitem {}\fi
\ifx \bconflocation  \undefined \def \bconflocation#1{#1}\fi
\ifx \arxivurl  \undefined \def \arxivurl#1{\textsf{#1}}\fi
\csname PreBibitemsHook\endcsname

\bibitem[\protect\citeauthoryear{Shapiro et~al.}{2021}]{ShapiroLectures_2ND}
\begin{bbook}
\bauthor{\bsnm{Shapiro}, \binits{A.}},
\bauthor{\bsnm{Dentcheva}, \binits{D.}},
\bauthor{\bsnm{Ruszczynski}, \binits{A.}}:
\bbtitle{Lectures on Stochastic Programming: Modeling and Theory, Third
  Edition},
\bedition{3rd} edn.
\bpublisher{Society for Industrial and Applied Mathematics},
\blocation{Philadelphia, PA}
(\byear{2021})
\end{bbook}
\endbibitem

\bibitem[\protect\citeauthoryear{Ogryczak and
  Ruszczy{\'n}ski}{1999}]{Ogryczak1999}
\begin{barticle}
\bauthor{\bsnm{Ogryczak}, \binits{W.}},
\bauthor{\bsnm{Ruszczy{\'n}ski}, \binits{A.}}:
\batitle{{{F}rom stochastic dominance to mean-risk models: {S}emideviations as
  risk measures}}.
\bjtitle{European Journal of Operational Research}
\bvolume{116}(\bissue{1}),
\bfpage{33}--\blpage{50}
(\byear{1999})
\end{barticle}
\endbibitem

\bibitem[\protect\citeauthoryear{Ogryczak and
  Ruszczy\'{n}ski}{2002}]{Ogryczak2002}
\begin{barticle}
\bauthor{\bsnm{Ogryczak}, \binits{W.}},
\bauthor{\bsnm{Ruszczy\'{n}ski}, \binits{A.}}:
\batitle{{{D}ual stochastic dominance and related mean-risk models}}.
\bjtitle{SIAM Journal on Optimization}
\bvolume{13}(\bissue{1}),
\bfpage{60}--\blpage{78}
(\byear{2002})
\end{barticle}
\endbibitem

\bibitem[\protect\citeauthoryear{Artzner et~al.}{1999}]{Artzner1999}
\begin{barticle}
\bauthor{\bsnm{Artzner}, \binits{P.}},
\bauthor{\bsnm{Delbaen}, \binits{F.}},
\bauthor{\bsnm{Eber}, \binits{J.-M.}},
\bauthor{\bsnm{Heath}, \binits{D.}}:
\batitle{{{C}oherent measures of risk}}.
\bjtitle{Mathematical Finance}
\bvolume{9}(\bissue{3}),
\bfpage{203}--\blpage{228}
(\byear{1999})
\end{barticle}
\endbibitem

\bibitem[\protect\citeauthoryear{Balbás et~al.}{2009}]{Balbas2009}
\begin{barticle}
\bauthor{\bsnm{Balbás}, \binits{A.}},
\bauthor{\bsnm{Garrido}, \binits{J.}},
\bauthor{\bsnm{Mayoral}, \binits{S.}}:
\batitle{{P}roperties of {distortion} {risk} {measures}}.
\bjtitle{Methodology and Computing in Applied Probability}
\bvolume{11}(\bissue{3}),
\bfpage{385}--\blpage{399}
(\byear{2009})
\end{barticle}
\endbibitem

\bibitem[\protect\citeauthoryear{Ribeiro}{2012}]{Ribeiro2012}
\begin{barticle}
\bauthor{\bsnm{Ribeiro}, \binits{A.}}:
\batitle{{{O}ptimal resource allocation in wireless communication and
  networking}}.
\bjtitle{EURASIP Journal on Wireless Communications and Networking}
\bvolume{2012}(\bissue{1}),
\bfpage{272}
(\byear{2012})
\end{barticle}
\endbibitem

\bibitem[\protect\citeauthoryear{Boyd and Vandenberghe}{2004}]{Boyd2004}
\begin{bbook}
\bauthor{\bsnm{Boyd}, \binits{S.P.}},
\bauthor{\bsnm{Vandenberghe}, \binits{L.}}:
\bbtitle{Convex Optimization},
p. \bfpage{716}.
\bpublisher{Cambridge University Press},
\blocation{Cambridge, UK}
(\byear{2004})
\end{bbook}
\endbibitem

\bibitem[\protect\citeauthoryear{Ruszczy{\'{n}}ski}{2006}]{Ruszczynski2006b}
\begin{bbook}
\bauthor{\bsnm{Ruszczy{\'{n}}ski}, \binits{A.P.}}:
\bbtitle{{N}onlinear {O}ptimization},
p. \bfpage{448}.
\bpublisher{Princeton University Press},
\blocation{Princeton, NJ, USA}
(\byear{2006})
\end{bbook}
\endbibitem

\bibitem[\protect\citeauthoryear{Bertsekas}{2009}]{Bertsekas2009}
\begin{bbook}
\bauthor{\bsnm{Bertsekas}, \binits{D.P.}}:
\bbtitle{{C}onvex {O}ptimization {T}heory},
\bedition{1st} edn.,
p. \bfpage{246}.
\bpublisher{Athena Scientific},
\blocation{Nashua NH, US}
(\byear{2009})
\end{bbook}
\endbibitem

\bibitem[\protect\citeauthoryear{Chieu et~al.}{2020}]{CHIEU2020441}
\begin{barticle}
\bauthor{\bsnm{Chieu}, \binits{N.H.}},
\bauthor{\bsnm{Jeyakumar}, \binits{V.}},
\bauthor{\bsnm{Li}, \binits{G.}}:
\batitle{{C}onvexifiability of continuous and discrete nonnegative quadratic
  programs for gap-free duality}.
\bjtitle{European Journal of Operational Research}
\bvolume{280}(\bissue{2}),
\bfpage{441}--\blpage{452}
(\byear{2020})
\end{barticle}
\endbibitem

\bibitem[\protect\citeauthoryear{Geoffrion}{1974}]{geoffrion_1974}
\begin{barticle}
\bauthor{\bsnm{Geoffrion}, \binits{A.M.}}:
\batitle{{L}agrangean relaxation for integer programming}.
\bjtitle{Mathematical Programming Study}
\bvolume{2},
\bfpage{82}--\blpage{114}
(\byear{1974})
\end{barticle}
\endbibitem

\bibitem[\protect\citeauthoryear{Hoffman and
  Kruskal}{1956}]{hoffman_kruskal_1956}
\begin{bchapter}
\bauthor{\bsnm{Hoffman}, \binits{A.J.}},
\bauthor{\bsnm{Kruskal}, \binits{J.B.}}:
\bctitle{{I}ntegral boundary points of convex polyhedra}.
In: \beditor{\bsnm{Kuhn}, \binits{H.W.}},
\beditor{\bsnm{Tucker}, \binits{A.W.}} (eds.)
\bbtitle{Linear Inequalities and Related Systems}.
\bsertitle{Annals of Mathematics Studies},
vol. \bseriesno{38},
pp. \bfpage{223}--\blpage{246}.
\bpublisher{Princeton University Press},
\blocation{Princeton, NJ}
(\byear{1956}).
\bcomment{Reprinted in various collections; often cited as Hoffman \& Kruskal
  (1956).}
\end{bchapter}
\endbibitem

\bibitem[\protect\citeauthoryear{Lemaire and Volle}{1998}]{Lemaire1998}
\begin{bbook}
\bauthor{\bsnm{Lemaire}, \binits{B.}},
\bauthor{\bsnm{Volle}, \binits{M.}}:
In: \beditor{\bsnm{Crouzeix}, \binits{J.-P.}},
\beditor{\bsnm{Martinez-Legaz}, \binits{J.-E.}},
\beditor{\bsnm{Volle}, \binits{M.}} (eds.)
\bbtitle{{D}uality in {DC} programming},
pp. \bfpage{331}--\blpage{345}.
\bpublisher{Springer},
\blocation{Boston, MA}
(\byear{1998})
\end{bbook}
\endbibitem

\bibitem[\protect\citeauthoryear{Tassiulas and
  Ephremides}{1992}]{Tassiulas1992}
\begin{barticle}
\bauthor{\bsnm{Tassiulas}, \binits{L.}},
\bauthor{\bsnm{Ephremides}, \binits{A.}}:
\batitle{{{S}tability properties of constrained queueing systems and scheduling
  policies for maximum throughput in multihop radio networks}}.
\bjtitle{IEEE Transactions on Automatic Control}
\bvolume{37}(\bissue{12}),
\bfpage{1936}--\blpage{1948}
(\byear{1992})
\end{barticle}
\endbibitem

\bibitem[\protect\citeauthoryear{Neely et~al.}{2005}]{Neely2005}
\begin{bchapter}
\bauthor{\bsnm{Neely}, \binits{M.J.}},
\bauthor{\bsnm{Modiano}, \binits{E.}},
\bauthor{\bsnm{Rohrs}, \binits{C.E.}}:
\bctitle{{{D}ynamic power allocation and routing for time-varying wireless
  networks}}.
In: \bbtitle{IEEE Journal on Selected Areas in Communications},
vol. \bseriesno{23},
pp. \bfpage{89}--\blpage{103}
(\byear{2005})
\end{bchapter}
\endbibitem

\bibitem[\protect\citeauthoryear{Georgiadis et~al.}{2006}]{Georgiadis2006a}
\begin{barticle}
\bauthor{\bsnm{Georgiadis}, \binits{L.}},
\bauthor{\bsnm{Neely}, \binits{M.J.}},
\bauthor{\bsnm{Tassiulas}, \binits{L.}}:
\batitle{{{R}esource allocation and cross-layer control in wireless networks}}.
\bjtitle{Foundations and Trends in Networking}
\bvolume{1}(\bissue{1}),
\bfpage{1}--\blpage{144}
(\byear{2006})
\end{barticle}
\endbibitem

\bibitem[\protect\citeauthoryear{Sidiropoulos et~al.}{2006}]{Sidiropoulos2006}
\begin{barticle}
\bauthor{\bsnm{Sidiropoulos}, \binits{N.D.}},
\bauthor{\bsnm{Davidson}, \binits{T.N.}},
\bauthor{\bsnm{{Zhi-Quan Luo}}}:
\batitle{{{T}ransmit beamforming for physical-layer multicasting}}.
\bjtitle{IEEE Transactions on Signal Processing}
\bvolume{54}(\bissue{6}),
\bfpage{2239}--\blpage{2251}
(\byear{2006})
\end{barticle}
\endbibitem

\bibitem[\protect\citeauthoryear{{Wei Yu} and Lui}{2006}]{WeiYu2006}
\begin{barticle}
\bauthor{\bsnm{{Wei Yu}}},
\bauthor{\bsnm{Lui}, \binits{R.}}:
\batitle{{{D}ual methods for nonconvex spectrum optimization of multicarrier
  systems}}.
\bjtitle{IEEE Transactions on Communications}
\bvolume{54}(\bissue{7}),
\bfpage{1310}--\blpage{1322}
(\byear{2006})
\end{barticle}
\endbibitem

\bibitem[\protect\citeauthoryear{Bazerque and Giannakis}{2007}]{Bazerque2007}
\begin{bchapter}
\bauthor{\bsnm{Bazerque}, \binits{J.A.}},
\bauthor{\bsnm{Giannakis}, \binits{G.B.}}:
\bctitle{{D}istributed scheduling and resource allocation for cognitive ofdma
  radios}.
In: \bbtitle{Proceedings of the 2nd International Conference on Cognitive Radio
  Oriented Wireless Networks and Communications, CrownCom},
pp. \bfpage{343}--\blpage{350}.
\bpublisher{IEEE},
\blocation{Orlando, FL, USA}
(\byear{2007})
\end{bchapter}
\endbibitem

\bibitem[\protect\citeauthoryear{Luo and Zhang}{2008}]{Luo2008}
\begin{barticle}
\bauthor{\bsnm{Luo}, \binits{Z.Q.}},
\bauthor{\bsnm{Zhang}, \binits{S.}}:
\batitle{{{D}ynamic spectrum management: {C}omplexity and duality}}.
\bjtitle{IEEE Journal on Selected Topics in Signal Processing}
\bvolume{2}(\bissue{1}),
\bfpage{57}--\blpage{73}
(\byear{2008})
\end{barticle}
\endbibitem

\bibitem[\protect\citeauthoryear{Neely}{2010}]{Neely2010}
\begin{bbook}
\bauthor{\bsnm{Neely}, \binits{M.J.}}:
\bbtitle{{S}tochastic {Network} {Optimization} with {Application} to
  {Communication} and {Queueing} {Systems}}.
\bpublisher{Morgan and Claypool Publishers},
\blocation{Switzerland}
(\byear{2010})
\end{bbook}
\endbibitem

\bibitem[\protect\citeauthoryear{Ribeiro and Giannakis}{2010}]{Ribeiro2010}
\begin{barticle}
\bauthor{\bsnm{Ribeiro}, \binits{A.}},
\bauthor{\bsnm{Giannakis}, \binits{G.B.}}:
\batitle{{Separation principles in wireless networking}}.
\bjtitle{IEEE Transactions on Information Theory}
\bvolume{56}(\bissue{9}),
\bfpage{4488}--\blpage{4505}
(\byear{2010})
\end{barticle}
\endbibitem

\bibitem[\protect\citeauthoryear{Ribeiro}{2010}]{Ribeiro2010a}
\begin{barticle}
\bauthor{\bsnm{Ribeiro}, \binits{A.}}:
\batitle{{{E}rgodic stochastic optimization algorithms for wireless
  communication and networking}}.
\bjtitle{IEEE Transactions on Signal Processing}
\bvolume{58}(\bissue{12}),
\bfpage{6369}--\blpage{6386}
(\byear{2010})
\end{barticle}
\endbibitem

\bibitem[\protect\citeauthoryear{Hu and Ribeiro}{2012}]{Hu2012}
\begin{barticle}
\bauthor{\bsnm{Hu}, \binits{Y.}},
\bauthor{\bsnm{Ribeiro}, \binits{A.}}:
\batitle{{{O}ptimal wireless networks based on local channel state
  information}}.
\bjtitle{IEEE Transactions on Signal Processing}
\bvolume{60}(\bissue{9}),
\bfpage{4913}--\blpage{4929}
(\byear{2012})
\end{barticle}
\endbibitem

\bibitem[\protect\citeauthoryear{He et~al.}{2014}]{He2014}
\begin{barticle}
\bauthor{\bsnm{He}, \binits{T.}},
\bauthor{\bsnm{Wang}, \binits{X.}},
\bauthor{\bsnm{Ni}, \binits{W.}}:
\batitle{{{O}ptimal chunk-based resource allocation for {OFDMA} systems with
  multiple {BER} requirements}}.
\bjtitle{IEEE Transactions on Vehicular Technology}
\bvolume{63}(\bissue{9}),
\bfpage{4292}--\blpage{4301}
(\byear{2014})
\end{barticle}
\endbibitem

\bibitem[\protect\citeauthoryear{Naderializadeh and
  Avestimehr}{2014}]{Naderializadeh2014}
\begin{barticle}
\bauthor{\bsnm{Naderializadeh}, \binits{N.}},
\bauthor{\bsnm{Avestimehr}, \binits{A.S.}}:
\batitle{{{ITLinQ}: {A} new approach for spectrum sharing in device-to-device
  communication systems}}.
\bjtitle{IEEE Journal on Selected Areas in Communications}
\bvolume{32}(\bissue{6}),
\bfpage{1139}--\blpage{1151}
(\byear{2014})
\end{barticle}
\endbibitem

\bibitem[\protect\citeauthoryear{Eisen et~al.}{2019a}]{Eisen2019}
\begin{barticle}
\bauthor{\bsnm{Eisen}, \binits{M.}},
\bauthor{\bsnm{Zhang}, \binits{C.}},
\bauthor{\bsnm{Chamon}, \binits{L.F.O.}},
\bauthor{\bsnm{Lee}, \binits{D.D.}},
\bauthor{\bsnm{Ribeiro}, \binits{A.}}:
\batitle{{{L}earning optimal resource allocations in wireless systems}}.
\bjtitle{IEEE Transactions on Signal Processing}
\bvolume{67}(\bissue{10}),
\bfpage{2775}--\blpage{2790}
(\byear{2019})
\end{barticle}
\endbibitem

\bibitem[\protect\citeauthoryear{Eisen et~al.}{2019b}]{Eisen2019b}
\begin{barticle}
\bauthor{\bsnm{Eisen}, \binits{M.}},
\bauthor{\bsnm{Gatsis}, \binits{K.}},
\bauthor{\bsnm{Pappas}, \binits{G.J.}},
\bauthor{\bsnm{Ribeiro}, \binits{A.}}:
\batitle{{{L}earning in wireless control systems over nonstationary channels}}.
\bjtitle{IEEE Transactions on Signal Processing}
\bvolume{67}(\bissue{5}),
\bfpage{1123}--\blpage{1137}
(\byear{2019})
\end{barticle}
\endbibitem

\bibitem[\protect\citeauthoryear{Kalogerias et~al.}{2020}]{Kalogerias2020e}
\begin{barticle}
\bauthor{\bsnm{Kalogerias}, \binits{D.S.}},
\bauthor{\bsnm{Eisen}, \binits{M.}},
\bauthor{\bsnm{Pappas}, \binits{G.J.}},
\bauthor{\bsnm{Ribeiro}, \binits{A.}}:
\batitle{{{M}odel-free learning of optimal ergodic policies in wireless
  systems}}.
\bjtitle{IEEE Transactions on Signal Processing}
\bvolume{68},
\bfpage{6272}--\blpage{6286}
(\byear{2020})
\end{barticle}
\endbibitem

\bibitem[\protect\citeauthoryear{Chamon et~al.}{2020}]{Chamon2020a}
\begin{barticle}
\bauthor{\bsnm{Chamon}, \binits{L.F.O.}},
\bauthor{\bsnm{Eldar}, \binits{Y.C.}},
\bauthor{\bsnm{Ribeiro}, \binits{A.}}:
\batitle{{{F}unctional nonlinear sparse models}}.
\bjtitle{IEEE Transactions on Signal Processing}
\bvolume{68},
\bfpage{2449}--\blpage{2463}
(\byear{2020})
\end{barticle}
\endbibitem

\bibitem[\protect\citeauthoryear{Chamon et~al.}{2022}]{Chamon2021}
\begin{botherref}
\oauthor{\bsnm{Chamon}, \binits{L.F.O.}},
\oauthor{\bsnm{Paternain}, \binits{S.}},
\oauthor{\bsnm{Calvo-Fullana}, \binits{M.}},
\oauthor{\bsnm{Ribeiro}, \binits{A.}}:
{{C}onstrained learning with non-convex losses}.
IEEE Transactions on Information Theory
\textbf{63}(3)
(2022)
\end{botherref}
\endbibitem

\bibitem[\protect\citeauthoryear{Balseiro et~al.}{2023}]{Balseiro2023}
\begin{botherref}
\oauthor{\bsnm{Balseiro}, \binits{S.R.}},
\oauthor{\bsnm{Besbes}, \binits{O.}},
\oauthor{\bsnm{Pizarro}, \binits{D.}}:
{S}urvey of {dynamic} {resource}-{constrained} {reward} {collection}
  {problems}: {Unified} {model} and {analysis}.
Operations Research
\textbf{72}(5)
(2023)
\end{botherref}
\endbibitem

\bibitem[\protect\citeauthoryear{Gatsis et~al.}{2015}]{Gatsis2015}
\begin{barticle}
\bauthor{\bsnm{Gatsis}, \binits{K.}},
\bauthor{\bsnm{Pajic}, \binits{M.}},
\bauthor{\bsnm{Ribeiro}, \binits{A.}},
\bauthor{\bsnm{Pappas}, \binits{G.J.}}:
\batitle{{{O}pportunistic control over shared wireless channels}}.
\bjtitle{IEEE Transactions on Automatic Control}
\bvolume{60}(\bissue{12}),
\bfpage{3140}--\blpage{3155}
(\byear{2015})
\end{barticle}
\endbibitem

\bibitem[\protect\citeauthoryear{Gatsis et~al.}{2018}]{Gatsis2018}
\begin{barticle}
\bauthor{\bsnm{Gatsis}, \binits{K.}},
\bauthor{\bsnm{Ribeiro}, \binits{A.}},
\bauthor{\bsnm{Pappas}, \binits{G.J.}}:
\batitle{{Random Access Design for Wireless Control Systems}}.
\bjtitle{Automatica}
\bvolume{91},
\bfpage{1}--\blpage{9}
(\byear{2018})
\end{barticle}
\endbibitem

\bibitem[\protect\citeauthoryear{Berliocchi and
  Lasry}{1973}]{BerliocchiLasry1973}
\begin{barticle}
\bauthor{\bsnm{Berliocchi}, \binits{H.}},
\bauthor{\bsnm{Lasry}, \binits{J.-M.}}:
\batitle{{I}nt\'egrandes normales et mesures param\'etr\'ees en calcul des
  variations}.
\bjtitle{Bulletin de la Soci\'et\'e Math\'ematique de France}
\bvolume{101}(\bissue{2}),
\bfpage{129}--\blpage{184}
(\byear{1973})
\end{barticle}
\endbibitem

\bibitem[\protect\citeauthoryear{Aubin and Ekeland}{1976}]{AubinEkeland1976}
\begin{barticle}
\bauthor{\bsnm{Aubin}, \binits{J.-P.}},
\bauthor{\bsnm{Ekeland}, \binits{I.}}:
\batitle{{E}stimates of the duality gap in non-convex optimization}.
\bjtitle{Mathematics of Operations Research}
\bvolume{1}(\bissue{3}),
\bfpage{225}--\blpage{245}
(\byear{1976})
\end{barticle}
\endbibitem

\bibitem[\protect\citeauthoryear{Aumann and Perles}{1965}]{AumannPerles1965}
\begin{barticle}
\bauthor{\bsnm{Aumann}, \binits{R.J.}},
\bauthor{\bsnm{Perles}, \binits{M.}}:
\batitle{{A} variational problem arising in economics}.
\bjtitle{Journal of Mathematical Analysis and Applications}
\bvolume{11}(\bissue{1}),
\bfpage{488}--\blpage{503}
(\byear{1965})
\end{barticle}
\endbibitem

\bibitem[\protect\citeauthoryear{Ekeland and Témam}{1999}]{EkelandTemam1999}
\begin{bbook}
\bauthor{\bsnm{Ekeland}, \binits{I.}},
\bauthor{\bsnm{Témam}, \binits{R.}}:
\bbtitle{{C}onvex {A}nalysis and {V}ariational {P}roblems}.
\bsertitle{Classics in Applied Mathematics},
vol. \bseriesno{28},
p. \bfpage{402}.
\bpublisher{Society for Industrial and Applied Mathematics (SIAM)},
\blocation{Philadelphia, PA}
(\byear{1999}).
\bcomment{Translated from the French; first English ed. 1976}
\end{bbook}
\endbibitem

\bibitem[\protect\citeauthoryear{Giannessi}{2005}]{Giannessi2005}
\begin{bbook}
\bauthor{\bsnm{Giannessi}, \binits{F.}}:
\bbtitle{Constrained Optimization and Image Space Analysis}.
\bpublisher{Springer},
\blocation{New York, NY}
(\byear{2005})
\end{bbook}
\endbibitem

\bibitem[\protect\citeauthoryear{Diestel and Uhl}{1977}]{Diestel1977}
\begin{bbook}
\bauthor{\bsnm{Diestel}, \binits{J.}},
\bauthor{\bsnm{Uhl}, \binits{J.R.}}:
\bbtitle{Vector Measures},
\bedition{Math. surv} edn.
\bpublisher{American Mathematical Society},
\blocation{Providence}
(\byear{1977})
\end{bbook}
\endbibitem

\bibitem[\protect\citeauthoryear{Lyapunov}{1940}]{Lyapunov_original}
\begin{barticle}
\bauthor{\bsnm{Lyapunov}, \binits{A.A.}}:
\batitle{{S}ur les fonctions-vecteurs compl{\`e}tement additives}.
\bjtitle{Izvestiya Akademii Nauk SSSR. Seriya Matematicheskaya}
\bvolume{4},
\bfpage{465}--\blpage{478}
(\byear{1940})
\end{barticle}
\endbibitem

\bibitem[\protect\citeauthoryear{Aumann}{1965}]{AUMANN19651}
\begin{barticle}
\bauthor{\bsnm{Aumann}, \binits{R.J.}}:
\batitle{{I}ntegrals of set-valued functions}.
\bjtitle{Journal of Mathematical Analysis and Applications}
\bvolume{12}(\bissue{1}),
\bfpage{1}--\blpage{12}
(\byear{1965})
\end{barticle}
\endbibitem

\bibitem[\protect\citeauthoryear{Flores-Bazán and
  Mastroeni}{2013}]{FloresBazan2013}
\begin{barticle}
\bauthor{\bsnm{Flores-Bazán}, \binits{F.}},
\bauthor{\bsnm{Mastroeni}, \binits{G.}}:
\batitle{{S}trong {duality} in {cone} {constrained} {nonconvex}
  {optimization}}.
\bjtitle{SIAM Journal on Optimization}
\bvolume{23}(\bissue{1}),
\bfpage{153}--\blpage{169}
(\byear{2013})
\end{barticle}
\endbibitem

\bibitem[\protect\citeauthoryear{Luenberger}{1968}]{Luenberger1968}
\begin{bbook}
\bauthor{\bsnm{Luenberger}, \binits{D.G.}}:
\bbtitle{{Optimization by Vector Space Methods}}.
\bpublisher{Wiley},
\blocation{New York, NY}
(\byear{1968})
\end{bbook}
\endbibitem

\bibitem[\protect\citeauthoryear{Yaylali and Kalogerias}{2023}]{Yaylali2023}
\begin{bchapter}
\bauthor{\bsnm{Yaylali}, \binits{G.}},
\bauthor{\bsnm{Kalogerias}, \binits{D.}}:
\bctitle{{R}obust and reliable stochastic resource allocation via tail
  waterfilling}.
In: \bbtitle{2023 IEEE 24th International Workshop on Signal Processing
  Advances in Wireless Communications (SPAWC)},
pp. \bfpage{256}--\blpage{260}
(\byear{2023})
\end{bchapter}
\endbibitem

\bibitem[\protect\citeauthoryear{Yaylali and Kalogerias}{2024}]{Yaylali2024}
\begin{bchapter}
\bauthor{\bsnm{Yaylali}, \binits{G.}},
\bauthor{\bsnm{Kalogerias}, \binits{D.}}:
\bctitle{Distributionally robust power policies for wireless systems under
  power fluctuation risk}.
In: \bbtitle{2024 58th Asilomar Conference on Signals, Systems, and Computers},
pp. \bfpage{1005}--\blpage{1012}
(\byear{2024})
\end{bchapter}
\endbibitem

\bibitem[\protect\citeauthoryear{Yaylali et~al.}{2025}]{Yaylali2025}
\begin{bchapter}
\bauthor{\bsnm{Yaylali}, \binits{G.}},
\bauthor{\bsnm{Khan}, \binits{A.A.}},
\bauthor{\bsnm{Kalogerias}, \binits{D.}}:
\bctitle{Waterfilling at the edge: Optimal percentile resource allocation via
  risk-averse reduction}.
In: \bbtitle{2025 IEEE 10th International Workshop on Computational Advances in
  Multi-Sensor Adaptive Processing (CAMSAP)},
pp. \bfpage{181}--\blpage{185}
(\byear{2025})
\end{bchapter}
\endbibitem

\bibitem[\protect\citeauthoryear{Yaylali and Kalogerias}{2026}]{Yaylali2026}
\begin{barticle}
\bauthor{\bsnm{Yaylali}, \binits{G.}},
\bauthor{\bsnm{Kalogerias}, \binits{D.}}:
\batitle{Distributionally robust ultra-reliable resource allocation via double
  tail waterfilling under fading risk}.
\bjtitle{IEEE Open Journal of Signal Processing}
\bvolume{7},
\bfpage{154}--\blpage{164}
(\byear{2026})
\end{barticle}
\endbibitem

\bibitem[\protect\citeauthoryear{Rockafellar and
  Uryasev}{2000}]{Rockafellar1997}
\begin{barticle}
\bauthor{\bsnm{Rockafellar}, \binits{R.T.}},
\bauthor{\bsnm{Uryasev}, \binits{S.}}:
\batitle{{{O}ptimization of conditional value-at-risk}}.
\bjtitle{Journal of Risk}
\bvolume{2},
\bfpage{21}--\blpage{41}
(\byear{2000})
\end{barticle}
\endbibitem

\bibitem[\protect\citeauthoryear{Rockafellar}{1971}]{Rockafellar1971Integrals}
\begin{barticle}
\bauthor{\bsnm{Rockafellar}, \binits{R.T.}}:
\batitle{{I}ntegrals which are convex functionals}.
\bjtitle{Pacific Journal of Mathematics}
\bvolume{24}(\bissue{3}),
\bfpage{525}--\blpage{539}
(\byear{1971})
\end{barticle}
\endbibitem

\bibitem[\protect\citeauthoryear{Hiai and Umegaki}{1977}]{HIAI1977149}
\begin{barticle}
\bauthor{\bsnm{Hiai}, \binits{F.}},
\bauthor{\bsnm{Umegaki}, \binits{H.}}:
\batitle{{I}ntegrals, conditional expectations, and martingales of multivalued
  functions}.
\bjtitle{Journal of Multivariate Analysis}
\bvolume{7}(\bissue{1}),
\bfpage{149}--\blpage{182}
(\byear{1977})
\end{barticle}
\endbibitem

\bibitem[\protect\citeauthoryear{Pennanen and
  Perkkiö}{2024}]{PennanenPerkkio2024}
\begin{bbook}
\bauthor{\bsnm{Pennanen}, \binits{T.}},
\bauthor{\bsnm{Perkkiö}, \binits{A.-P.}}:
\bbtitle{{C}onvex {S}tochastic {O}ptimization: {D}ynamic {P}rogramming and
  {D}uality in {D}iscrete {T}ime}.
\bsertitle{Probability Theory and Stochastic Modelling},
vol. \bseriesno{107}.
\bpublisher{Springer},
\blocation{Cham, Switzerland}
(\byear{2024}).
\bcomment{1st edition}
\end{bbook}
\endbibitem

\bibitem[\protect\citeauthoryear{Rockafellar and
  Wets}{2004}]{Rockafellar2009VarAn}
\begin{bbook}
\bauthor{\bsnm{Rockafellar}, \binits{R.T.}},
\bauthor{\bsnm{Wets}, \binits{R.J.-B.}}:
\bbtitle{Variational Analysis}
vol. \bseriesno{317},
p. \bfpage{734}.
\bpublisher{Springer},
\blocation{Heidelberg, DE}
(\byear{2004})
\end{bbook}
\endbibitem

\bibitem[\protect\citeauthoryear{Kingman and
  Robertson}{1968}]{KingmanRobertson}
\begin{barticle}
\bauthor{\bsnm{Kingman}, \binits{J.F.C.}},
\bauthor{\bsnm{Robertson}, \binits{A.P.}}:
\batitle{On a theorem of lyapunov}.
\bjtitle{Journal of the London Mathematical Society}
\bvolume{s1-43}(\bissue{1}),
\bfpage{347}--\blpage{351}
(\byear{1968})
\end{barticle}
\endbibitem

\bibitem[\protect\citeauthoryear{Knowles}{1975}]{SIAMCon:Knowles}
\begin{barticle}
\bauthor{\bsnm{Knowles}, \binits{G.}}:
\batitle{{L}yapunov vector measures}.
\bjtitle{SIAM Journal on Control}
\bvolume{13}(\bissue{2}),
\bfpage{294}--\blpage{303}
(\byear{1975})
\end{barticle}
\endbibitem

\bibitem[\protect\citeauthoryear{Uhl}{1969}]{Uhl1969}
\begin{barticle}
\bauthor{\bsnm{Uhl}, \binits{J.J.}}:
\batitle{{{T}he range of a vector-valued measure}}.
\bjtitle{Proceedings of the American Mathematical Society}
\bvolume{23}(\bissue{1}),
\bfpage{158}--\blpage{163}
(\byear{1969})
\end{barticle}
\endbibitem

\bibitem[\protect\citeauthoryear{Kadets and Schechtman}{1992}]{Kadets1992}
\begin{barticle}
\bauthor{\bsnm{Kadets}, \binits{V.M.}},
\bauthor{\bsnm{Schechtman}, \binits{G.}}:
\batitle{Lyapunov’s theorem for $\ell_p$-valued measures}.
\bjtitle{Algebra i Analiz}
\bvolume{4}(\bissue{5}),
\bfpage{148}--\blpage{154}
(\byear{1992}).
\bcomment{Translation in St. Petersburg Math. J. 4(5), 961–966 (1993)}
\end{barticle}
\endbibitem

\bibitem[\protect\citeauthoryear{Bertsekas}{2016}]{Bertsekas2016Nonlinear}
\begin{bbook}
\bauthor{\bsnm{Bertsekas}, \binits{D.P.}}:
\bbtitle{{N}onlinear {P}rogramming},
\bedition{3rd} edn.
\bpublisher{Athena Scientific},
\blocation{Belmont, MA}
(\byear{2016})
\end{bbook}
\endbibitem

\bibitem[\protect\citeauthoryear{Steinwart and
  Christmann}{2008}]{steinwart2008support}
\begin{bbook}
\bauthor{\bsnm{Steinwart}, \binits{I.}},
\bauthor{\bsnm{Christmann}, \binits{A.}}:
\bbtitle{{S}upport {V}ector {M}achines}.
\bsertitle{Information Science and Statistics}.
\bpublisher{Springer},
\blocation{New York, NY}
(\byear{2008})
\end{bbook}
\endbibitem

\bibitem[\protect\citeauthoryear{Guliyev and Ismailov}{2018}]{GULIYEV2018262}
\begin{barticle}
\bauthor{\bsnm{Guliyev}, \binits{N.J.}},
\bauthor{\bsnm{Ismailov}, \binits{V.E.}}:
\batitle{{A}pproximation capability of two hidden layer feedforward neural
  networks with fixed weights}.
\bjtitle{Neurocomputing}
\bvolume{316},
\bfpage{262}--\blpage{269}
(\byear{2018})
\end{barticle}
\endbibitem

\bibitem[\protect\citeauthoryear{Kalogerias et~al.}{2020}]{Kalogerias2020b}
\begin{bchapter}
\bauthor{\bsnm{Kalogerias}, \binits{D.S.}},
\bauthor{\bsnm{Eisen}, \binits{M.}},
\bauthor{\bsnm{Pappas}, \binits{G.J.}},
\bauthor{\bsnm{Ribeiro}, \binits{A.}}:
\bctitle{{A} zeroth-order learning algorithm for ergodic optimization of
  wireless systems with no models and no gradients}.
In: \bbtitle{ICASSP, IEEE International Conference on Acoustics, Speech and
  Signal Processing - Proceedings},
pp. \bfpage{5045}--\blpage{5049}.
\bpublisher{Institute of Electrical and Electronics Engineers Inc.},
\blocation{Barcelona, ES}
(\byear{2020})
\end{bchapter}
\endbibitem

\bibitem[\protect\citeauthoryear{Goldsmith}{2005}]{Goldsmith2005Wireless}
\begin{bbook}
\bauthor{\bsnm{Goldsmith}, \binits{A.}}:
\bbtitle{{Wireless {C}ommunications}}.
\bpublisher{Cambridge university press}, \blocation{???}
(\byear{2005})
\end{bbook}
\endbibitem

\bibitem[\protect\citeauthoryear{Tse and Viswanath}{2005}]{Tse2005}
\begin{bbook}
\bauthor{\bsnm{Tse}, \binits{D.}},
\bauthor{\bsnm{Viswanath}, \binits{P.}}:
\bbtitle{{F}undamentals of {{W}ireless} {{C}ommunication}},
\bedition{Illustrated edition} edn.
\bpublisher{Cambridge University Press},
\blocation{Cambridge, UK ; New York}
(\byear{2005})
\end{bbook}
\endbibitem

\bibitem[\protect\citeauthoryear{Kalogerias and Powell}{2018}]{Kalogerias2018b}
\begin{botherref}
\oauthor{\bsnm{Kalogerias}, \binits{D.S.}},
\oauthor{\bsnm{Powell}, \binits{W.B.}}:
{{R}ecursive optimization of convex risk measures: {M}ean-semideviation
  models}.
arXiv:1804.00636
(2018)
\end{botherref}
\endbibitem

\end{thebibliography}

\end{document}